\title{$A_2$-Planar Algebras I}
\author{
        David E. Evans and Mathew Pugh \\ \\
        School of Mathematics, \\
        Cardiff University, \\
        Senghennydd Road, \\
        Cardiff, CF24 4AG, \\
        Wales, U.K.
}
\date{\today}

\documentclass[12pt]{article}
\usepackage{amssymb}
\usepackage{graphicx}
\usepackage[all]{xy}

\newtheorem{Def}{Definition}[section]
\newtheorem{Prop}[Def]{Proposition}
\newtheorem{Lemma}[Def]{Lemma}
\newtheorem{Cor}[Def]{Corollary}
\newtheorem{Thm}[Def]{Theorem}

\begin{document}
\maketitle

\begin{abstract}
We give a diagrammatic presentation of the $A_2$-Temperley-Lieb algebra.
Generalizing Jones' notion of a planar algebra, we formulate an $A_2$-planar algebra motivated by Kuperberg's $A_2$-spider.
This $A_2$-planar algebra contains a subfamily of vector spaces which will capture the double complex structure pertaining to the subfactor for a finite $SU(3)$ $\mathcal{ADE}$ graph with a flat cell system, including both the periodicity three coming from the $A_2$-Temperley-Lieb algebra as well as the periodicity two coming from the subfactor basic construction.
We use an $A_2$-planar algebra to obtain a description of the (Jones) planar algebra for the Wenzl subfactor in terms of generators and relations. 
\end{abstract}

\section{Introduction}

A braided inclusion $N \subset M$, where there is a braided system of $SU(3)_k$ endomorphisms ${}_N \mathcal{X}_N$ on the factor $N$, yields a nimrep (non-negative integer matrix representation) of the right action of the $N$-$N$ sectors ${}_N \mathcal{X}_N$ on the $M$-$N$ sectors ${}_M \mathcal{X}_N$ via the theory of $\alpha$-induction \cite{bockenhauer/evans:1998, bockenhauer/evans:1999i, bockenhauer/evans:1999ii}. This nimrep defines a classifying graph $\mathcal{G}$, which is of $\mathcal{ADE}$ type.
One can build an Ocneanu cell system $W$ on an $\mathcal{ADE}$ graph $\mathcal{G}$ \cite{ocneanu:2000ii}, which attaches a complex number to each closed path of length three on the edges of $\mathcal{G}$. A cell system $W$ naturally gives rise to a representation of the Hecke algebra, or more precisely, of the $A_2$-Temperley-Lieb algebra \cite{evans/pugh:2009i, evans/pugh:2009ii}, which is a quotient of the Hecke algebra given by the fixed point algebra of $\bigotimes_{\mathbb{N}} M_3$ under the action of $SU(3)_k$.
This $A_2$-Temperley-Lieb algebra has an inherent periodicity of three coming from the representation theory of $SU(3)$.

To each pair $(\mathcal{G},W)$, consisting of an $SU(3)$ $\mathcal{ADE}$ graph $\mathcal{G}$ and a cell system $W$ on $\mathcal{G}$, there is associated a subfactor $N \subset M$, or rather, a subfactor double complex (c.f. \cite{evans/kawahigashi:1994}), which has a periodicity of three in the horizontal direction, coming from the $A_2$-Temperley-Lieb algebra, and a periodicity of two in the vertical direction, coming from the subfactor basic construction of Jones \cite{jones:1983}, or equivalently, from the (usual) Temperley-Lieb algebra.
The subfactor double complex contains the tower of higher relative commutants $N' \subset M_i$ as its initial column, where $N \subset M \subset M_1 \subset M_2 \ldots \;$ is the tower obtained by iterating the basic construction. However, it also contains the $SU(3)$ structure captured by the $A_2$-Temperley-Lieb operators, which is lost in the tower of higher relative commutants, or indeed in the standard invariant.

The main goal of this paper is to provide a framework for an $A_2$ version of a planar algebra which describes the subfactor double complex. We begin by giving a diagrammatic presentation of the $A_2$-Temperley-Lieb algebra, consisting of $A_2$-tangles which are a special class of Kuperberg's $A_2$ webs \cite{kuperberg:1996}. The $A_2$-Temperley-Lieb algebra is the underlying algebra in our $A_2$-planar algebra, which is a family of vector spaces which carry an action of the $A_2$-tangles.

The main result of the paper is Theorem \ref{thm:type_I_subfactors_give_flat_planar_algebras}, where for any pair $(\mathcal{G},W)$ we explicitly associate to the corresponding subfactor double complex an $A_2$-planar algebra, that is, there is an action of the $A_2$-tangles on each finite-dimensional vector space in the subfactor double complex. As an immediate corollary we obtain a description of the (usual) planar algebra for Wenzl's Hecke subfactor in terms of generators and relations \cite{wenzl:1988}.
This work provides a framework for studying subfactor double complexes, even in the continuous $SU(3)$ regime beyond index nine.

\section{Preliminaries}

A subfactor encodes symmetries. These can be understood and studied from a number of vantage points and directions which have interlocking ideas. In a subfactor's most fundamental setting, these symmetries may arise from a group, a group dual or a Hopf algebra, and their actions on a von Neumann algebra $M$, but subfactor symmetries go far beyond this, and beyond quantum groups. The symmetries of a group $G$ and group dual may be recovered from the position of the fixed point algebra $M^G$ in the ambient algebra $M$ and the position of $M$ in the crossed product $M \rtimes G$.
More generally, the symmetry or quantum symmetry is encoded by the position of a von Neumann algebra in another. Subfactors encode data, algebraic, combinatorial and analytic, and the question arises as to how to recover the data from the subfactor $N \subset M$ and vice versa.

Iterating the basic construction of Jones \cite{jones:1983} in the type $\mathrm{II}_1$ setting, one obtains a tower $N \subset M \subset M_1 \subset M_2 \ldots \;$. The standard invariant is obtained by considering the tower of relative commutants $M_i' \cap M_j$, which are finite dimensional in the case of finite index. Different axiomatizations of the standard invariant are given by Ocneanu with paragroups \cite{ocneanu:1988}, emphasising connections and their flatness, and by Popa with $\lambda$-lattices and a more probabilistic language which permit reconstruction of the (extremal finite index) subfactor under certain amenable conditions \cite{popa:1995}.
Jones \cite{jones:planar} produced another formulation using planar algebras, a diagrammatic incarnation of the relative commutants, closed under planar contractions or carrying operations indexed by certain planar diagrams, such that any extremal subfactor gives a planar algebra. Conversely, using the work of Popa on $\lambda$-lattices, every planar algebra, with suitable positivity properties, produces an extremal finite index subfactor. Recent work of \cite{guionnet/jones/shlyakhtenko:2007}, this time with a free probabilistic input of ideas, has recovered the characterisation of Popa.

The most fundamental symmetry of a subfactor is through the Temperley Lieb algebra \cite{temperley/lieb:1971}.
The Jones basic construction $M_{i-1} \subset M_i \subset M_{i+1}$ is through adjoining an extra projection $e_i$ arising from the projection or conditional expectation of $M_i$ onto $M_{i-1}$. These projections satisfy the Temperley-Lieb relations of integrable statistical mechanics. They are contained in the tower of relative commutants of any finite index subfactor and are in some sense the minimal symmetries. The planar algebra of a subfactor also has to encode what else is there, but in the case of the Temperley-Lieb algebra its planar algebra corresponds to Kauffman's diagrammatic presentation of the Temperley-Lieb algebra. The Temperley-Lieb algebra has a realization from $SU(2)$, from the fixed point algebras of quantum $SU(2)$ on the Pauli algebra and special representations of Hecke algebras of type $A$. These $SU(2)$ subfactors generalize to $SU(3)$ (and beyond \cite{wenzl:1988, wassermann:1998}). These subfactors can be used to understand $SU(3)$ orbifold subfactors, conformal embeddings and modular invariants \cite{evans/kawahigashi:1994, xu:1998, bockenhauer/evans:1998, bockenhauer/evans:1999i, bockenhauer/evans:1999ii, evans/pugh:2009ii}.

Here we give a planar study of subfactors which encodes the representation theory of quantum $SU(3)$ diagrammatically. The Temperley-Lieb algebra is then generalized to the following.
The Hecke algebra $H_n(q)$, $q \in \mathbb{C}$, is the algebra generated by invertible operators $g_j$, $j=1,2,\ldots,n-1$, satisfying the relations
\begin{eqnarray}
(q^{-1} - g_j)(q + g_j) & = & 0, \\
g_i g_j & = & g_j g_i, \quad |i-j|>1,\\
g_i g_{i+1} g_i & = & g_{i+1} g_i g_{i+1}. \label{eqn:braiding_relation}
\end{eqnarray}
When $q=1$, the first relation becomes $g_j^2 = 1$, so that $H_n(1)$ reduces to the group ring of the symmetric, or permutation, group $S_n$, where $g_j$ represents a transposition $(j,j+1)$.
Writing $g_j = q^{-1} - U_j$ where $|q|=1$, and setting $\delta = q+q^{-1}$, these generators and relations lead to self-adjoint operators $\mathbf{1}, U_1, U_2, \ldots, U_{n-1}$ and relations
\begin{center}
\begin{minipage}[b]{12.5cm}
 \begin{minipage}[b]{1.3cm}
  \begin{eqnarray*}
  \textrm{H1:}\\
  \textrm{H2:}\\
  \textrm{H3:}
  \end{eqnarray*}
 \end{minipage}
 \hspace{2cm}
 \begin{minipage}[b]{7cm}
  \begin{eqnarray*}
  U_i^2 & = & \delta U_i,\\
  U_i U_j & = & U_j U_i, \quad |i-j|>1,\\
  U_i U_{i+1} U_i - U_i & = & U_{i+1} U_i U_{i+1} - U_{i+1},
  \end{eqnarray*}
 \end{minipage}
\end{minipage}
\end{center}
where $\delta = q+q^{-1}$.

To any $\sigma$ in the permutation group $S_n$, decomposed into transpositions of nearest neighbours $\sigma = \prod_{i \in I_{\sigma}} \tau_{i,i+1}$, we associate the operator $g_{\sigma} = \prod_{i \in I_{\sigma}} g_i$, which is well defined because of the braiding relation (\ref{eqn:braiding_relation}).
Then the commutant of the quantum group $SU(N)_q$ is obtained from the Hecke algebra by imposing an extra condition, which is the vanishing of the $q$-antisymmetrizer \cite{di_francesco/zuber:1990}
\begin{equation}
\label{SU(N)q condition}
\sum_{\sigma \in S_{n}} (-q)^{|I_{\sigma}|} g_{\sigma} = 0.
\end{equation}
For $SU(2)$ it reduces to the Temperley-Lieb condition $U_i U_{i \pm 1} U_i - U_i = 0$, whilst for $SU(3)$ it is
\begin{equation}
\label{eqn:SU(3)q_condition}
\left( U_i - U_{i+2} U_{i+1} U_i + U_{i+1} \right) \left( U_{i+1} U_{i+2} U_{i+1} - U_{i+1} \right) = 0.
\end{equation}
The $A_2$-Temperley-Lieb algebra will be the algebra generated by a family $\{ U_n \}$ of self-adjoint operators which satisfy the Hecke relations H1-H3 and the extra condition (\ref{eqn:SU(3)q_condition}) (c.f. \cite{martin:1991, brzezinski/katriel:1995}).
The $A_2$-Temperley-Lieb algebra is the fixed point algebra of $\bigotimes_{\mathbb{N}} M_3$ under the product action of $\bigotimes_{\mathbb{N}} \mathrm{Ad}({\rho})$ of $SU(3)$ or its quantum version $SU(3)_q$. There is an inherent periodicity three which comes from the representation theory of $SU(3)$, which is reflected in the Bratteli diagram of the McKay graph of the fusion of the fundamental representation. For $q$ a $k^{\mathrm{th}}$ root of unity or $q=1$, the $A_2$-Temperley-Lieb algebra is isomorphic to the path algebra of the $SU(3)$ graph $\mathcal{A}^{(k+3)}$, where $k=\infty$ for $q=1$, which is tripartite, or three-colourable, so that all closed paths on $\mathcal{A}^{(k+3)}$ have lengths which are multiples of three.

\subsection{Background on Jones' planar algebras}

Jones introduced the notion of a planar algebra in \cite{jones:planar} to study subfactors.
Let us briefly review the essential construction of Jones' planar algebras. A planar $k$-tangle consists of a disc $D$ in the plane with $2k$ vertices on its boundary, $k \geq 0$, and $n \geq 0$ internal discs $D_j$, $j=1,\ldots,n$, where the disc $D_j$ has $2k_j$ vertices on its boundary, $k_j \geq 0$. One vertex on the boundary of each disc (including the outer disc $D$) is chosen as a marked vertex, and the segment of the boundary of each disc between the marked vertex and the vertex immediately adjacent to it as we move around the boundary in an anti-clockwise direction is labelled either $+$ or $-$. For a disc which has no vertices on its boundary, we label its entire boundary by $+$ or $-$. Inside $D$ we have a collection of disjoint smooth curves, called strings, where any string is either a closed loop, or else has as its endpoints the vertices on the discs, and such that every vertex is the endpoint of exactly one string. Any tangle must also allow a checkerboard colouring of the regions inside $D$, which are bounded by the strings and the boundaries of the discs, where every region is coloured black or white such that any two regions which share a common boundary are not coloured the same, and any region which meets the boundary of a disc at the segment marked $+$, $-$ is coloured black, white respectively.

A planar $k$-tangle with an internal disc $D_j$ with $2k_j$ vertices on its boundary can be composed with a $k_j$-tangle $S$, giving a new $k$-tangle $T \circ_j S$, by inserting the tangle $S$ inside the inner disc $D_j$ of $T$ such that the vertices on the outer disc of $S$ coincide with those on the disc $D_j$, and in particular the two marked vertices must coincide. The boundary of the disc $D_j$ is then removed, and the strings are smoothed if necessary. The collection of all diffeomorphism classes of such planar tangles, with composition defined as above, is called the planar operad.

A planar algebra $P$ is then defined to be an algebra over this operad, i.e. a family $P = (P_k^+, P_k^-; k \geq 0)$ of vector spaces with $P_k^{\pm} \subset P_{k'}^{\pm}$ for $k < k'$, and with the following property. For every $k$-tangle $T$ with $n$ internal discs $D_j$ labelled by elements $x_j \in P_{k_j}$, $j=1,\ldots,n$, there is an associated linear map $Z(T):\otimes_{j=1}^n P_{k_j} \rightarrow P_k$, which is compatible with the composition of tangles and re-ordering of internal discs.

These planar algebras gave a topological reformulation of the standard invariant, described in terms of relative commutants in the standard tower of a subfactor. More precisely, the standard invariant of an extremal subfactor $N \subset M$ is a (subfactor) planar algebra $P = (P_k)_{k \geq 0}$ with $P_k = N' \cap M_{k-1}$. Conversely, every planar algebra can be realised by a subfactor \cite{popa:1995, jones:planar} (see also \cite{guionnet/jones/shlyakhtenko:2007, jones/shlyakhtenko/walker:2008, kodiyalam/sunder:2008}).
The index \cite{jones:1983} is a crude measure of the complexity of a subfactor -- those subfactors with index $<4$ being the simplest. Since every relative commutant contains the Temperley-Lieb algebra, another notion of complexity is the number of non-Temperley-Lieb elements that are required to generate the relative commutants. In the planar algebra set-up, planar algebras $P$ generated by a single element, for which the dimension of $P_3$ is at most 13, were classified in \cite{bisch/jones:2000+2003}. In the recent work of \cite{kodiyalam/tupurani:2010i+ii} it was shown that any  subfactor planar algebra $P$ of depth $k$ is generated by a single element in $P_t$, for some $t \leq k+1$.

In \cite{jones:2001} Jones studied \emph{annular} tangles, that is, tangles with a distinguished internal disc. He introduced the notion of modules over a planar algebra, which are modules over an annular category whose morphisms are given by such annual tangles, and gave a description of all irreducible Temperley-Lieb modules. A more general planar algebra is the graph planar algebra of a bipartite graph \cite{jones:2000}. Jones and Reznikoff obtained the decomposition of the graph planar algebras for the $ADE$ graphs into irreducible Temperley-Lieb modules \cite{jones:2001, reznikoff:2005}. A similar notion to an tangle is that of an \emph{affine} tangle. Affine Temperley-Lieb algebras were studied in \cite{jones/reznikoff:2006, reznikoff:2008}.

One way to construct planar algebras is by generators and relations. One problem that arises with this method is to determine whether or not a set of generators and relations will produce a finite dimensional planar algebra, that is, a planar algebra $P$ where each $P_k$, $k>0$, is finite-dimensional. Landau \cite{landau:2002} obtained a condition called an \emph{exchange relation}, which guarantees that a planar algebra is in fact finite dimensional, and this condition was extended and generalized in \cite{ghosh:2003}.
A bigger problem is to show whether or not the trace defined on the planar algebra is positive definite. The graph planar algebras have a positive definite trace.
A recently published result in \cite[Corollary 4.2]{jones/penneys:2010} says that every finite-depth subfactor planar algebra is a planar subalgebra of the graph planar algebra of its principal graph. If a planar algebra can be found as a planar subalgebra of a graph planar algebra then the trace it inherits from the graph planar algebra will be automatically positive definite.
This motivated the construction of the planar algebra for the $ADE$ subfactors in terms of generators and relations \cite{bigelow:2010, morrison/peters/snyder:2008}, and more recently for the Haagerup subfactor \cite{peters:2009}, and the extended Haagerup subfactor \cite{bigelow/morrison/peters/snyder:2009} where planar algebras were used to show the existence of the extended Haagerup subfactor for the first time.

The planar algebras associated to different constructions of subfactors have been described: the planar algebra associated to subfactors arising from the outer actions on a factor by a finite-dimensional Kac algebra \cite{kodiyalam/landau/sunder:2003}, by a semisimple and cosemisimple Hopf algebra \cite{kodiyalam/sunder:2006} and more recently by the actions of finite groups or finitely generated, countable, discrete groups \cite{ghosh:2006, gupta:2008, bisch/das/ghosh:2008i, bisch/das/ghosh:2008ii}. Planar algebras associated to the action of compact quantum groups on finite quantum spaces were studied in \cite{banica:2005}.

\section{Taking Jones' planar algebras to the $A_2$ setting} \label{sect:A2tangles}

Our planar description naturally begins in this section with the spiders of Kuperberg \cite{kuperberg:1996} who developed some of the basic diagrammatics of the representation theory of $A_2$ and other rank two Lie algebras. Here we give a diagrammatic presentation of the $A_2$-Temperley-Lieb algebra using Kuperberg's $A_2$ spider, and show that the $A_2$-Temperley-Lieb algebra is isomorphic to Wenzl's quotient of the Hecke algebra \cite{wenzl:1988}.
In Section \ref{sec:A_2-planar_algebras} we introduce and study the notion of a general $A_2$-planar algebra and in Section \ref{sect:A2-planar_algebras+flatness} the notion of an $A_2$-planar algebra and the notion of flatness. In Section \ref{sect:i,j-tangles} we describe particular subspaces that we are interested in, which will correspond exactly to the double complex associated to the $SU(3)$-subfactors.

The $SU(3)$ $\mathcal{ADE}$ graphs appear as nimreps for the $SU(3)$ modular invariants \cite{evans/pugh:2009i, evans/pugh:2009ii}. For each graph there is a construction of a subfactor via a double complex of finite-dimensional algebras (cf. $\lambda$-lattice in what one could call the $SU(2)$ setting) which relies on the existence of a cell system which defines a connection or Boltzmann weight. The series of the commuting squares in these double complexes are not canonical in the sense of Popa, because although these double complexes have period 2 vertically (coming from the subfactor basic construction) they have period 3 horizontally (coming from the underlying $A_2$-Temperley-Lieb algebraic structure). These double complexes were used by Evans and Kawahigashi \cite{evans/kawahigashi:1994} to understand the Wenzl subfactors and their orbifolds, and in particular to compute their principal graphs.
The main result of the paper is Theorem \ref{thm:type_I_subfactors_give_flat_planar_algebras} in Section \ref{sect:Planar_algebras_give_subfactors}, where we show how the subfactor, or associated double complex, for a finite $\mathcal{ADE}$ graph with a flat cell system diagrammatically gives rise to a flat $A_2$-$C^{\ast}$-planar algebra.
Jones' ($A_1$-)planar algebra is contained in the $A_2$-planar algebra, as the algebra over a certain suboperad of our $A_2$-planar operad.
In Section \ref{sect:planar_alg_for_A=PTL} we obtain an $A_2$-planar algebra description of the Wenzl subfactor, and as a corollary we have a construction of Jones' planar algebra for the Wenzl subfactor in terms of generators and relations which come from the $A_2$-planar algebra.

In \cite{evans/pugh:2009i} we computed the numerical values of the Ocneanu cells, announced by Ocneanu (e.g. \cite{ocneanu:2000ii, ocneanu:2002}), and consequently representations of the Hecke algebra, for the $SU(3)$ $\mathcal{ADE}$ graphs. These cells assign a numerical weight to Kuperberg's diagram of trivalent vertices -- corresponding to the fact that the trivial representation is contained in the triple product of the fundamental representation of $SU(3)$ through the determinant. They will yield, in a natural way, representations of an $A_2$-Temperley-Lieb or Hecke algebra.
For bipartite graphs, the corresponding weights (associated to the diagrams of cups or caps), arise in a more straightforward fashion from a Perron-Frobenius eigenvector, giving a natural representation of the Temperley-Lieb algebra or Hecke algebra.

In the sequel \cite{evans/pugh:2009iv} we introduce the notion of modules over an $A_2$-planar algebra, and describe certain irreducible Hilbert $A_2$-$TL$-modules. A partial decomposition of graph $A_2$-planar algebras for the $\mathcal{ADE}$ graphs is achieved. The graph $A_2$-planar algebra $P^{\mathcal{G}}$ of an $\mathcal{ADE}$ graph is an $A_2$-$C^{\ast}$-planar algebra with $\mathrm{dim}(P^{\mathcal{G}}_0) > 1$, which is a generalization of the bipartite graph planar algebra to the $A_2$ setting. These graph $A_2$-planar algebras are diagrammatic representations of another double complex of finite dimensional algebras, where now the initial space in the double complex is $\mathbb{C}^n$ where $n>1$ (note that $n=1$ for the initial space in the double complex associated to an $SU(3)$-subfactor).

The bipartite theory of the $SU(2)$ setting has to some degree become a three-colourable theory in our $SU(3)$ setting. This theory is not completely three-colourable since some of the graphs are not three-colourable -- namely the graphs $\mathcal{A}^{(n)\ast}$ associated to the conjugate modular invariants, $n \geq 4$, $\mathcal{D}^{(n)}$ associated to the orbifold modular invariants, $n \neq 0 \textrm{ mod } 3$, and the exceptional graph $\mathcal{E}^{(8)\ast}$. The figures for the complete list of the $\mathcal{ADE}$ graphs are given in \cite{behrend/pearce/petkova/zuber:2000, evans/pugh:2009i}.

We have laid the foundations for a planar algebra formulation of an $SU(3)$ theory which may help resolve some of the unanswered questions left open in the programme which we set out on in \cite{evans/pugh:2009i, evans/pugh:2009ii} to understand $SU(3)$ modular invariants and their representation by braided subfactors. We realised all $SU(3)$ modular invariants by braided $SU(3)$ subfactors \cite{evans/pugh:2009ii} but did not classify their associated nimreps or claim that the known list is exhaustive. In the case of one of the exceptional modular invariants, we could not identify the nimrep. We verified that all known candidate nimrep graphs carried Ocneanu cell systems \cite{evans/pugh:2009i}, apart from one exceptional graph $\mathcal{E}_4^{(12)}$. However, we did not determine when such a cell system yields a local braided subfactor, but speculated that this should correspond to type I cell systems, that is, cell systems such that the connection defined by equations (\ref{eqn:connection}), (\ref{eqn:inverse_graph_connection}) in the present paper is flat. This is only known for the $\mathcal{A}$ and $\mathcal{D}$ graphs at present \cite{evans/kawahigashi:1994}.

The question of whether all nimreps have been realised is open.
There are some nimreps which do not have braided subfactors. We also want to go beyond the $\mathcal{ADE}$ classification to study subfactors for more exotic graphs which support a cell system, just as Jones' planar algebras facilitated the study of the Haagerup and extended Haagerup subfactors. The tools being drawn up in this paper may aid these further studies.

\subsection{Orbifolds}

The orbifold construction is a standard procedure in operator algebras, in $C^{\ast}$-algebras and subfactor theory in von Neumann algebras, as well as in integrable statistical mechanics and conformal field theory. A finite abelian group action on the underlying structure can bring about an orbifold, by suitably dividing out by the group elements (usually called simple currents in conformal field theory) which may or may not describe completely different theory from the original one.
This usually depends on having fixed points, and understanding their role or the resolution of these singularities is the key.

For example, in the theory of $C^{\ast}$-algebras, the fixed point algebra of the irrational rotation algebra by a flip on the generators or the underlying two dimensional torus has an AF fixed point algebra, and so has a completely different character to the ambient noncommutative torus which has non trivial $K_1$.
This is reviewed with full references in \cite[notes to Ch.3, pp 125--146]{evans/kawahigashi:1998}.
The invariants involved in
understanding  or comparing orbifolds, the fixed point algebras or crossed products,  with the original algebras being $K$-theory or equivariant $K$-theory. Partly motivated by this, orbifold methods were introduced into subfactor theory \cite{evans/kawahigashi:1994}, but first we digress to the underlying statistical mechanics and conformal field theories.

In statistical mechanics, Date, Jimbo, Miwa and Okado \cite{date/jimbo/miwa/okado:1989} introduced integrable models associated with the level $k$-integrable models of the Kac-Moody algebra of $SU(n)$. The Boltzmann weights lie in the fixed point algebra of the infinite tensor product of $M_n$ under the action of $SU(n)_k$.

The notion of an orbifold of such a model by dividing out by a subgroup $Z$ of the centre of $SU(n)$ were introduced by Pasquier \cite{pasquier:1987}, Fendley and Ginsparg \cite{fendley/ginsparg:1989} for $n=2$ and by Di-Franceso and Zuber \cite{di_francesco/zuber:1990} for $n=3$, borrowing from an orbifold notion in conformal field theory \cite{dixon/harvey/vafa/witten:1985+1986}. In the Wess-Zumino-Witten model, a two dimensional  conformal field theory arises from classical fields taking values in the target $SU(n)$ models and their orbifolds by $Z$ are meant to be those living in the quotient $SU(n)/Z$.

With all this mind, the orbifold construction was introduced in subfactor theory in \cite{evans/kawahigashi:1994}, with the Boltzmann weights being in the relevant fixed point algebras and hence naturally satisfy the Yang-Baxter equation, and the subfactors introduced through the action of the subgroup $Z$ of the centre as a group of automorphisms and crossed products. It is still a question whether one is really finding a new subfactor, as in the $N=2$ case, one cannot simply take the orbifold of the $A_{4m-1}$-principal graph which would be $D_{2m+1}$, as only $D_m$ for $m$ even can arise as a principal graph of a subfactor. For the $SU(3)$ subfactor the action of the center $\mathbb{Z}_3$ of $SU(3)$ introduces an action for each integer level $k$ on ${}_N \mathcal{X}_N$, a system of endomorphisms of a type III factor $N$ represented by the vertices
of the truncated diagram $\mathcal{A}^{(k+3)}$.

These orbifolds are best understood through $\alpha$-induction in subfactor theory \cite[Section 3]{bockenhauer/evans:1999i}, \cite[Section 6.2]{bockenhauer/evans:1999ii}, \cite[Section 8]{bockenhauer/evans:2002}, which we summarize here.

Simple currents \cite{schellekens/yankielowicz:1989} are primary fields with unit quantum dimension and appear in the subfactor framework as automorphisms in the system ${}_N \mathcal{X}_N$. They form a closed abelian group under fusion. Simple currents give rise to modular invariants, and all such invariants have been classified \cite{gato-rivera/schellekens:1992, kreuzer/schellekens:1994}.
We are focussing on $SU(n)$ here for $n=2,3$, and so will only consider cyclic simple current groups $\mathbb{Z}_n$.

By taking a generator $[\sigma]$ of the cyclic simple current group $\mathbb{Z}_n$ we can construct the crossed product subfactor $N \subset M = N \rtimes \mathbb{Z}_n$ whenever we can choose a representative $\sigma$ in each such simple current sector such that we have exact cyclicity $\sigma^n= 1$ (and not only as sectors).
Rehren's lemma \cite{rehren:1990} states that such a choice is possible if and only if the statistics phase $\omega_{\sigma}$ is an $n$-th root of unity, i.e. if and only if the conformal weight $h_{\sigma}$ is an integer multiple of $1/n$.
This construction gives rise to a non-trivial subfactor and in turn to a modular invariant.
For $SU(n)_k$ the simple current group $\mathbb{Z}_n$ corresponds to weights $k\Lambda_{(j)}$, $j=0,1,...,n-1$.
The conformal dimensions are $h_{k\Lambda_{(j)}}=kj(n-j)/2n$, which by Rehren's Lemma \cite{rehren:1990} allow for full $\mathbb{Z}_n$ extensions except when $n$ is even and $k$ is odd in which case the maximal extension is $N \subset M = N \rtimes \mathbb{Z}_{n/2}$ because we can only use the even labels $j$.
(This reflects the fact that e.g. for $SU(2)$ there are no D-invariants at odd levels.)
Thus Rehren's lemma has told us that extensions are labelled by all the divisors of $n$ unless $n$ is even and $k$ is odd in which case they are labelled by the divisors of $n/2$.
This matches exactly the simple current modular invariant classification of \cite{gato-rivera/schellekens:1992, kreuzer/schellekens:1994}.
An extension by a simple current subgroup $\mathbb{Z}_m$, with $m$ a divisor of $n$ or $n/2$, is moreover local, if the generating current (and hence all in the $\mathbb{Z}_m$ subgroup) has integer conformal weight, $h_{k\Lambda_{(q)}} \in \mathbb{Z}$, where $n=mq$.
This happens exactly if $kq \in 2m\mathbb{Z}$ if $n$ is even, or $kq \in m\mathbb{Z}$ if $n$ is odd \cite{bockenhauer/evans:1999ii}.
For $SU(2)$ this corresponds to the $\mathrm{D}_{\mathrm{even}}$ series whereas the $\mathrm{D}_{\mathrm{odd}}$ series are non-local extensions. For $SU(3)$, there is a simple current extension at each level, but only those at $k \in 3\mathbb{Z}$ are local.
For the case of $SU(3)$ at level $3p$, the crossed product $N \subset N \rtimes \mathbb{Z}_3$ with canonical endomorphism $[\theta]  = [\lambda_{(0,0)}] \oplus [\lambda_{(3p,0)}] \oplus [\lambda_{(0,3p)}]$, the procedure of alpha induction \cite[p.89]{bockenhauer/evans:1999i} yields from
$\langle \alpha_\lambda, \alpha_\mu \rangle = \langle \theta \lambda, \mu \rangle$
that at the fixed point $f = (p,p)$, $[\alpha_f] = [\alpha_f^{(1)}] \oplus [\alpha_f^{(2)}] \oplus [\alpha_f^{(3)}]$ splits into three irreducibles whilst otherwise $[\alpha_{\lambda}]$ is irreducible and identified with $[\alpha_{\sigma \lambda}]$, $\sigma \in \mathbb{Z}_3$, under the action of the centre $\mathbb{Z}_3$ or simple currents. Thus under alpha induction, the Verlinde algebra or the tensor category of $SU(3)$ at level $3p$, represented by a system of endomorphisms ${}_N \mathcal{X}_N$ is taken to its orbifold ${}_N \mathcal{X}_N^{\pm}$, and taking the dual action reverses this procedure. The principal graphs (the fusion graphs of $[\alpha_{(1,0)}]$) are the orbifold graphs $\mathcal{D}^{3p+3}$.

M{\"u}ger \cite{muger:2000} and Bruguieres \cite{bruguieres:2000} have subsequently introduced an orbifold procedure which can handle non abelian groups, and this procedure is sometimes described as equivariantization/deequivariantization in the category oriented literature.
We pointed out in \cite{evans/pugh:2009i} recent work in condensed matter physics \cite{bais/slingerland:2008} where we see that $\alpha$-induction is playing a key role.
For example, the computation of pages 8--9 is $\alpha$-induction for an orbifold embedding of $SU(2)_4$ which gives fusion graph $D_4$. Other examples are the conformal embedding of $SU(2)_4 \subset SU(3)_1$ on pages 14--15, which again gives fusion graph $D_4$, and the conformal embedding $SU(2)_{10} \subset SO(5)_1$ on pages 15--16, which gives fusion graph $E_6$.

\subsection{$A_2$-tangles} \label{sec:A_2-tangles}

In \cite{kuperberg:1996}, Kuperberg defined the notion of a spider, which is an axiomatization of the representation theory of groups and other group-like objects. The invariant spaces have bases given by certain planar graphs. These graphs are called webs, hence the term spider. In \cite{kuperberg:1996} certain spiders were defined in terms of generators and relations, isomorphic to the representation theories of rank two Lie algebras and the quantum deformations of these representation theories. This formulation generalized a well-known construction for $A_1 = \textrm{su}(2)$ by Kauffman \cite{kauffman:1987}.

\begin{figure}[bt]
\begin{center}
  \includegraphics[width=40mm]{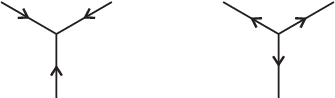}\\
 \caption{$A_2$ webs}\label{fig:A2-webs}
\end{center}
\end{figure}

For the $A_2 = \textrm{su}(3)$ case, we have the $A_2$ webs, illustrated in Figure \ref{fig:A2-webs}. We will call these webs \textbf{incoming and outgoing trivalent vertices} respectively. We call the oriented lines \textbf{strings}. We may join the $A_2$ webs together by attaching free ends of outgoing trivalent vertices to free ends of incoming trivalent vertices, and isotoping the strings if needed so that they are smooth.

We are now going to systematically define an algebra of web tangles, and express this in terms of generators and relations.

\begin{Def}
An \textbf{$A_2$-tangle} will be a connected collection of strings joined together at incoming or outgoing trivalent vertices (see Figure \ref{fig:A2-webs}), possibly with some free ends, such that the orientations of the individual strings are consistent with the orientations of the trivalent vertices.
\end{Def}

\begin{Def}
We call a vertex a \textbf{source} vertex if the string attached to it has orientation away from the vertex. Similarly, a \textbf{sink} vertex will be a vertex where the string attached has orientation towards the vertex.
\end{Def}

\begin{Def}

For $m,n \geq 0$, an \textbf{$A_2$-$(m,n)$-tangle} will be an $A_2$-tangle $T$ on a rectangle, where $T$ has $m+n$ free ends attached to $m$ source vertices along the top of the rectangle and $n$ sink vertices along the bottom such that the orientation of the strings is respected. If $m=n$ we call $T$ simply an $A_2$-$m$-tangle, and we position the vertices so that for every vertex along the top there is a corresponding vertex directly beneath it along the bottom.

Two $A_2$-$(m,n)$-tangles are equivalent if one can be obtained from the other by an isotopy which moves the strings and trivalent vertices, but leaves the boundary vertices unchanged. We define $\mathcal{T}^{A_2}_{m,n}$ to be the set of all (equivalence classes of) $A_2$-$(m,n)$-tangles.
\end{Def}

The composition $TS \in \mathcal{T}^{A_2}_{m,k}$ of an $A_2$-$(m,n)$-tangle $T$ and an $A_2$-$(n,k)$-tangle $S$ is given by gluing $S$ vertically below $T$ such that the vertices at the bottom of $T$ and the top of $S$ coincide, removing these vertices, and isotoping the glued strings if necessary to make them smooth. The composition is clearly associative.

\begin{Def}
We define the vector space $\mathcal{V}^{A_2}_{m,n}$ to be the free vector space over $\mathbb{C}$ with basis $\mathcal{T}^{A_2}_{m,n}$. Then $\mathcal{V}^{A_2}_{m,n}$  has an algebraic structure with multiplication given by composition of tangles. In particular, we will write $\mathcal{V}^{A_2}_{m}$ for $\mathcal{V}^{A_2}_{m,m}$, and $\mathcal{V}^{A_2} = \bigcup_{m \geq 0} \mathcal{V}^{A_2}_{m}$. For $n<m$ we have $\mathcal{V}^{A_2}_{n} \subset \mathcal{V}^{A_2}_{m}$, with the inclusion of an $n$-tangle $\mathcal{T} \in \mathcal{T}^{A_2}_{n}$ in $\mathcal{T}^{A_2}_{m}$ given by adding $m-n$ vertices along the top and bottom of the rectangle after the rightmost vertex, with $m-n$ downwards oriented vertical strings connecting the extra vertices along the top to those along the bottom. The inclusion for $\mathcal{V}^{A_2}_{n}$ in $\mathcal{V}^{A_2}_{m}$ is the linear extension of this map.
\end{Def}

Note that $\mathcal{T}^{A_2}_{m,n}$ is infinite, and thus the vector space $\mathcal{V}^{A_2}_{m,n}$ is infinite dimensional. However, we will take a quotient of $\mathcal{V}^{A_2}_{m,n}$ which will turn out to be finite dimensional.
Let K1-K3 denote the following relations on local parts of tangles, for $\alpha, \delta \in \mathbb{C}$ \cite{kuperberg:1996}:
\\
\vspace{-12mm}
\begin{center}
\begin{minipage}[b]{11.5cm}
 \begin{minipage}[t]{3cm}
  \parbox[t]{2cm}{\begin{eqnarray*}\textrm{K1:}\end{eqnarray*}}
 \end{minipage}
 \begin{minipage}[t]{5.5cm}
  \begin{center}
  \mbox{} \\
 \includegraphics[width=20mm]{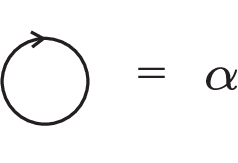}
  \end{center}
 \end{minipage}
 \begin{minipage}[t]{2cm}
  \mbox{} \\
  \parbox[t]{1cm}{}
 \end{minipage}
\end{minipage}
\begin{minipage}[b]{11.5cm}
 \begin{minipage}[t]{3cm}
  \parbox[t]{2cm}{\begin{eqnarray*}\textrm{K2:}\end{eqnarray*}}
 \end{minipage}
 \begin{minipage}[t]{5.5cm}
  \begin{center}
  \mbox{} \\
 \includegraphics[width=23mm]{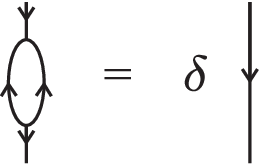}
  \end{center}
 \end{minipage}
 \begin{minipage}[t]{2cm}
  \mbox{} \\
  \parbox[t]{1cm}{}
 \end{minipage}
\end{minipage}
\begin{minipage}[b]{11.5cm}
 \begin{minipage}[t]{3cm}
  \parbox[t]{2cm}{\begin{eqnarray*}\textrm{K3:}\end{eqnarray*}}
 \end{minipage}
 \begin{minipage}[t]{5.5cm}
  \begin{center}
  \mbox{} \\
 \includegraphics[width=55mm]{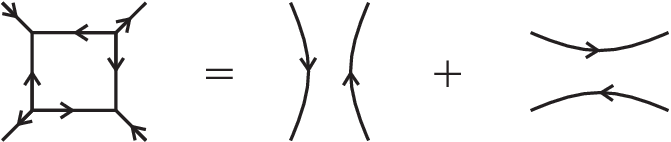}
  \end{center}
 \end{minipage}
 \begin{minipage}[t]{2cm}
  \mbox{} \\
  \parbox[t]{1cm}{}
 \end{minipage}
\end{minipage}
\end{center}

\begin{Def}
We define $I_{m,n} \subset \mathcal{V}^{A_2}_{m,n}$ to be the ideal of $\mathcal{V}^{A_2}_{m,n}$ which is the linear span of the relations K1-K3.
\end{Def}

By the linear span of the relations K1-K3 is meant the linear span of the differences of the left hand side and the right hand side of each of the relations, as local parts of the tangles, where the rest of the tangle is identical in each term in the difference. We will denote $I_{m,m}$ by $I_m$. Note that $I_m \subset I_{m+1}$.

\begin{Def}
The algebra $V^{A_2}_{m}$ is defined to be the quotient of the space $\mathcal{V}^{A_2}_{m}$ by the ideal $I_m$, and $V^{A_2} = \bigcup_{m \geq 0} V^{A_2}_{m}$.
\end{Def}

A basis of $V^{A_2}_{m}$ is given by all $A_2$-$m$-tangles which do not contain the local pictures which appear on the left hand side of K1-K3 (which Kuperberg calls elliptic faces). We will call the local picture $\,$ \includegraphics[width=10mm]{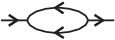} $\,$ a digon, and $\,$ \includegraphics[width=8mm]{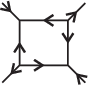} $\,$ an embedded square.
We could replace the Kuperberg relation K1 by the more general relations:\\

\begin{center}
\begin{minipage}[b]{11.5cm}
 \begin{minipage}[t]{2.5cm}
  \parbox[t]{2cm}{K1':}
 \end{minipage}
 \begin{minipage}[t]{6.5cm}
   \includegraphics[width=60mm]{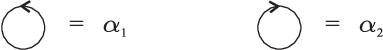}
 \end{minipage}
 \begin{minipage}[t]{2cm}
  \parbox[t]{1cm}{}
 \end{minipage}
\end{minipage}
\end{center}

Although it now appears that we have three independent parameters $\alpha_1, \alpha_1, \delta$, we actually have only one, as shown in the following Lemma:

\begin{Lemma}\label{Lemma:delta^2=alpha+1}
For a fixed complex number $\delta \neq 0$ we must have either $\alpha_1 = \alpha_2 = \delta^2 - 1$ or $\alpha_1 = \alpha_2 = 0$.
\end{Lemma}
\emph{Proof:}
\begin{figure}[tb]
\begin{center}
  \includegraphics[width=140mm]{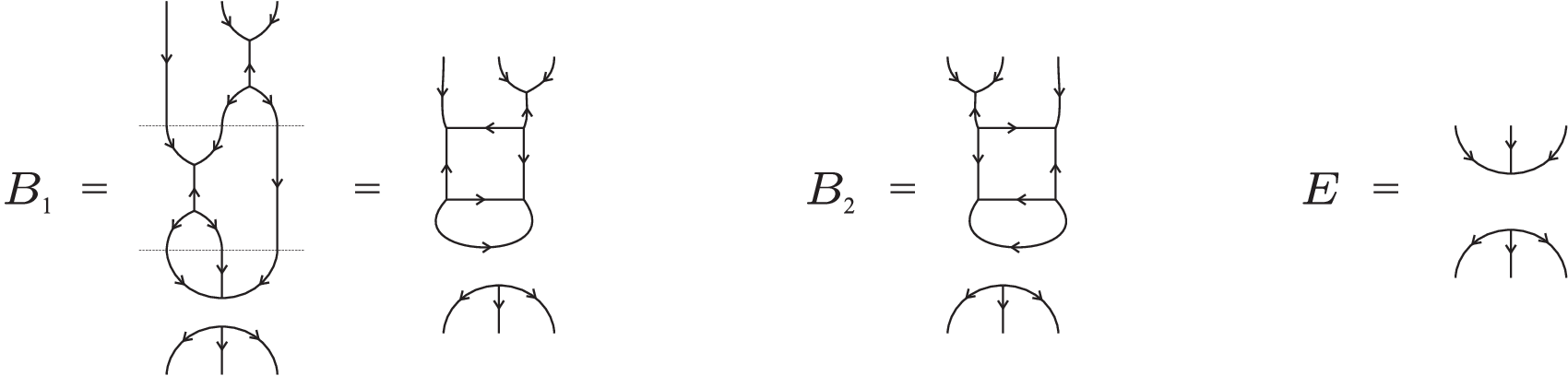}\\
 \caption{3-tangles $B_1$, $B_2$, $E$}\label{fig:3-tangles_B1,B2,E}
\end{center}
\end{figure}
Let $B_1$ be the 3-tangle illustrated in Figure \ref{fig:3-tangles_B1,B2,E}, which is the composition of three basis tangles in $V^{A_2}_{3}$. Let $B_2$ be a 3-tangle which comes from a similar composition, and $E$ a basis tangle in $V^{A_2}_{3}$, both also illustrated in Figure \ref{fig:3-tangles_B1,B2,E}.
Reducing $B_1$ using K2 twice, we get $B_1 = \delta^2 E$. On the other hand, if we reduce $B_1$ using K3, we get an anticlockwise oriented closed loop, which by K1' contributes a scalar factor $\alpha_1$. Then we also have $B_1 = E + \alpha_1 E$.
If $E \neq 0$, then $\delta^2 = 1 + \alpha_1$, and by the same argument on $B_2$ we also obtain $\delta^2 = 1 + \alpha_2$.
Suppose now that $E = 0$.
Let $\widehat{E}$ be the tangle given by composition of $E$ (embedded in $V^{A_2}_{6}$) with three nested caps above and three nested cups below, i.e. $\widehat{E}$ is the tangle
\begin{center}
\includegraphics[width=45mm]{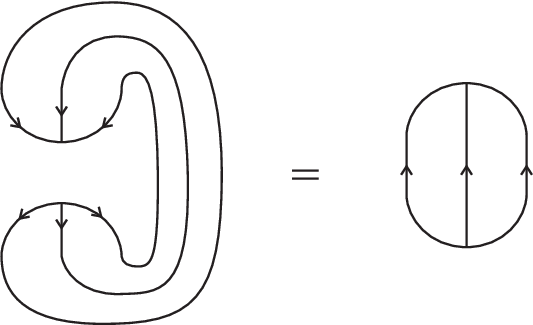}
\end{center}
If we use K2 to remove the left digon, we obtain an anticlockwise oriented loop, and so the diagram counts as the scalar $\alpha_1 \delta$. If instead we used K2 to remove the right digon we would obtain the scalar $\alpha_2 \delta$.
Since $\widehat{E} = 0$ and $\delta \neq 0$, we have $\alpha_1 = \alpha_2 = 0$.
\hfill
$\Box$

For $m \in \mathbb{Z}$, we define the quantum integer $[m]_q$ by $[m]_q = (q^m-q^{-m})/(q-q^{-1})$, where $q \in \mathbb{C}$. Note that if $\delta = [2]_q$, then by Lemma \ref{Lemma:delta^2=alpha+1} $\alpha = \delta^2 - 1 = [3]_q$ (or zero). When $q$ is an $n^{\mathrm{th}}$ root of unity, $q = e^{2\pi i/n}$, we will usually write $[m]$ for $[m]_q$.

There is a braiding on $V^{A_2}$, defined locally by the following linear combinations of local diagrams in $V^{A_2}$, for a choice of third root $q^{1/3}$, $q \in \mathbb{C}$ (see \cite{kuperberg:1996, suciu:1997}):
\begin{center}
\begin{minipage}[b]{16cm}
 \begin{minipage}[t]{4.5cm}
  \mbox{} \\
  \parbox[t]{1cm}{}
 \end{minipage}
 \begin{minipage}[t]{5.5cm}
  \mbox{} \\
   \includegraphics[width=55mm]{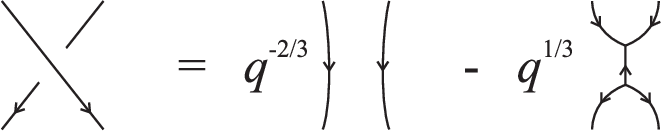}
 \end{minipage}
 \begin{minipage}[t]{5.5cm}
  \hfill
  \parbox[t]{2cm}{\begin{eqnarray}\label{braiding1}\end{eqnarray}}
 \end{minipage}
\end{minipage} \\
\mbox{} \\
\begin{minipage}[b]{16cm}
 \begin{minipage}[t]{4.5cm}
  \mbox{} \\
  \parbox[t]{1cm}{}
 \end{minipage}
 \begin{minipage}[t]{5.5cm}
  \mbox{} \\
   \includegraphics[width=55mm]{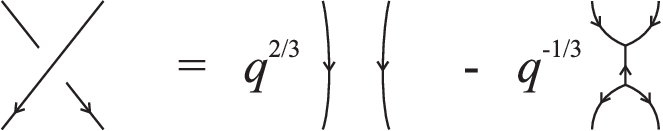}
 \end{minipage}
 \begin{minipage}[t]{5.5cm}
  \hfill
  \parbox[t]{2cm}{\begin{eqnarray}\label{braiding2}\end{eqnarray}}
 \end{minipage}
\end{minipage}
\end{center}

\noindent
The braiding satisfies the following properties locally, provided $\delta = [2]_q$ and $\alpha = [3]_q$:
\begin{center}
\begin{minipage}[b]{16cm}
 \begin{minipage}[t]{1cm}
  \hfill
  \parbox[t]{0.5cm}{}
 \end{minipage}
 \begin{minipage}[t]{13cm}
  \begin{center}
  \mbox{} \\
   \includegraphics[width=110mm]{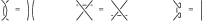}
  \end{center}
 \end{minipage}
 \begin{minipage}[t]{1.5cm}
  \hfill
  \parbox[t]{1.5cm}{\begin{eqnarray}\label{braiding_relations1-3}\end{eqnarray}}
 \end{minipage}
\end{minipage} \\
\mbox{} \\
\begin{minipage}[b]{16cm}
 \begin{minipage}[t]{2.5cm}
  \mbox{} \\
  \parbox[t]{1cm}{}
 \end{minipage}
 \begin{minipage}[t]{10cm}
  \begin{center}
  \mbox{} \\
   \includegraphics[width=90mm]{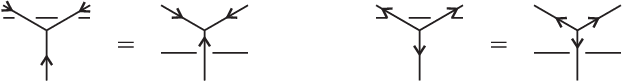}
  \end{center}
 \end{minipage}
 \begin{minipage}[t]{3cm}
  \hfill
  \parbox[t]{2cm}{\begin{eqnarray}\label{braiding_relation4}\end{eqnarray}}
 \end{minipage}
\end{minipage}
\end{center}
where we also have relation (\ref{braiding_relation4}) with the crossings all reversed.

\begin{figure}[h]
\begin{center}
  \includegraphics[width=90mm]{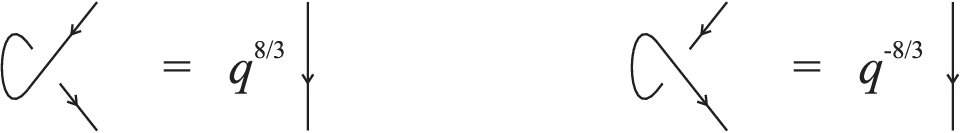}\\
 \caption{Removing kinks}\label{fig:expanding_kinks}
\end{center}
\end{figure}

We call the local pictures illustrated on the left hand sides of relations (\ref{braiding1}), (\ref{braiding2}) respectively a negative, positive crossing respectively. With this braiding, kinks (or twists) contribute a scalar factor of $q^{8/3}$ for those involving a positive crossing, and $q^{-8/3}$ for those involving a negative crossing, as shown in Figure \ref{fig:expanding_kinks}.

We now define a $\ast$-operation on $\mathcal{V}^{A_2}_{m}$, which is an involutive conjugate linear map. For an $m$-tangle $T \in \mathcal{T}^{A_2}_{m}$, $T^{\ast}$ is the $m$-tangle obtained by reflecting $T$ about a horizontal line halfway between the top and bottom vertices of the tangle, and reversing the orientations on every string. Then $\ast$ on $\mathcal{V}^{A_2}_{m}$ is the conjugate linear extension of $\ast$ on $\mathcal{T}^{A_2}_{m}$. Note that the $\ast$-operation leaves the relation K2 invariant if and only if $\delta \in \mathbb{R}$. For $\delta \in \mathbb{R}$, the $\ast$-operation leaves the ideal $I_m$ invariant due to the symmetry of the relations K1-K3. Then $\ast$ passes to $V^{A_2}_{m}$, and is an involutive conjugate linear anti-automorphism.

\subsection{Diagrammatic presentation of the $A_2$-Temperley-Lieb algebra} \label{sec:diagrammatic_A2-TL}

\begin{figure}[bt]
\begin{center}
  \includegraphics[width=40mm]{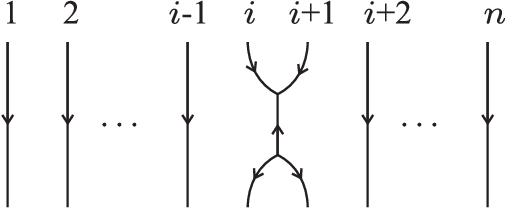}\\
 \caption{The $n$-tangle $W_i$, $i=1,\ldots, n-1$.}\label{fig:W_i}
\end{center}
\end{figure}

From now on we let $\delta$ be real, so that $\delta = [2]_q$ for some $q$, and we set $\alpha = [3]_q$ (cf. Lemma \ref{Lemma:delta^2=alpha+1}).
We define the tangle $\mathbf{1}_m$ to be the $m$-tangle with all strings vertical through strings. Then $\mathbf{1}_m$ is the identity of the algebra $\mathcal{V}^{A_2}_{m}$: $\mathbf{1}_m a = a = a \mathbf{1}_m$ for all $a \in \mathcal{V}^{A_2}_m$. We also define $W_i$ to be the $m$-tangle with all vertices along the top connected to the vertices along the bottom by vertical lines, except for the $i^{\textrm{th}}$ and $(i+1)^{\textrm{th}}$ vertices. The strings attached to the $i^{\textrm{th}}$ and $(i+1)^{\textrm{th}}$ vertices along the top are connected at an incoming trivalent vertex, with the third string coming from an outgoing trivalent vertex connected to the strings attached to the $i^{\textrm{th}}$ and $(i+1)^{\textrm{th}}$ vertices along the bottom. The tangle $W_i$ is illustrated in Figure \ref{fig:W_i}.

\setcounter{topnumber}{3}

\begin{figure}[tb]
\begin{center}
  \includegraphics[width=80mm]{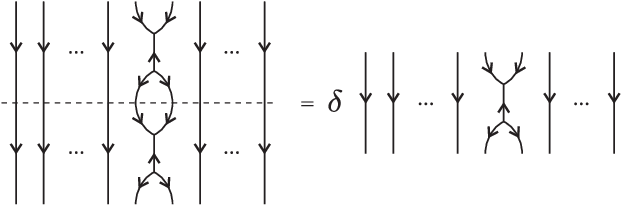}\\
 \caption{$w_i^2 = \delta w_i$}\label{fig:W_i_satisfies_Hecke1}
\end{center}
\end{figure}

\begin{figure}[tb]
\begin{center}
  \includegraphics[width=120mm]{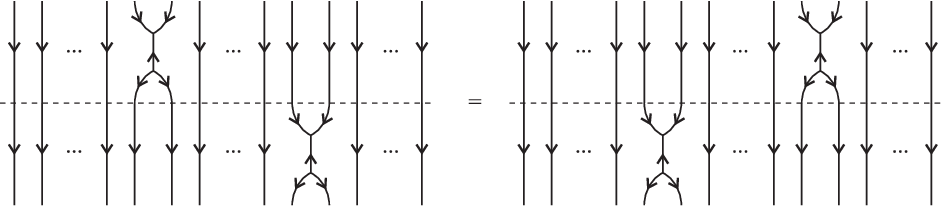}\\
 \caption{$w_i w_j = w_j w_i$ for $|i-j|>1$.}\label{fig:W_i_satisfies_Hecke2}
\end{center}
\end{figure}

\begin{figure}[tb]
\begin{center}
 \begin{minipage}[b]{16cm}
  \includegraphics[width=155mm]{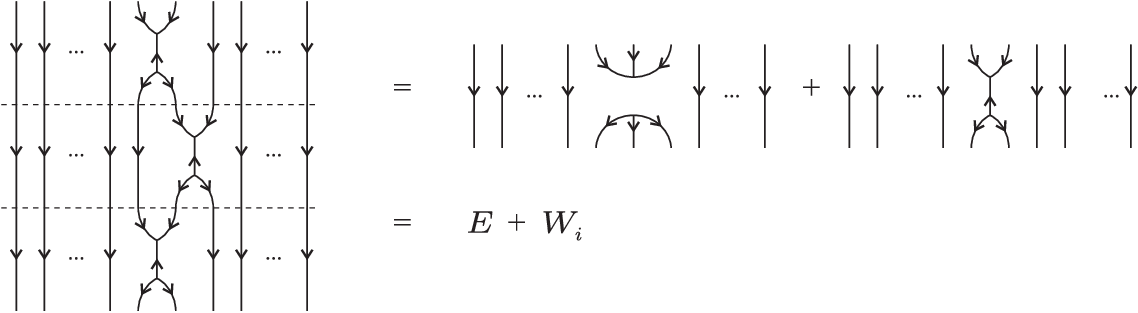}
 \end{minipage} \\
 \mbox{} \\
 \begin{minipage}[b]{16cm}
  \includegraphics[width=155mm]{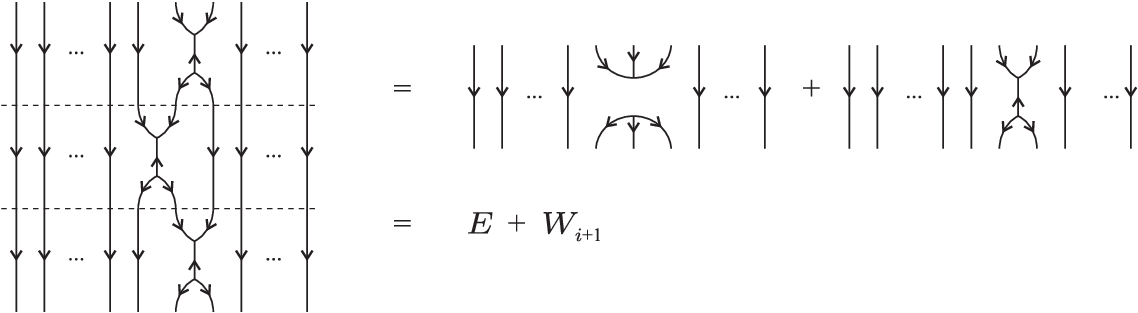}
 \end{minipage} \\
 \caption{$w_i w_{i+1} w_i - w_i = w_{i+1} w_i w_{i+1} - w_{i+1}$} \label{fig:W_i_satisfies_Hecke3}
\end{center}
\end{figure}

For $m \in \mathbb{N} \cup \{ 0 \}$ we define the algebra $A_2\textrm{-}TL_m$ to be $\mathrm{alg}(\mathbf{1}_m, w_i|i = 1, \ldots, m-1)$, where $w_i = W_i + I_m$. The $w_i$'s in $A_2\textrm{-}TL_m$ are clearly self-adjoint, and satisfy the relations H1-H3, as illustrated in Figures \ref{fig:W_i_satisfies_Hecke1}, \ref{fig:W_i_satisfies_Hecke2} and \ref{fig:W_i_satisfies_Hecke3}.

\begin{figure}[tb]
\begin{center}
  \includegraphics[width=40mm]{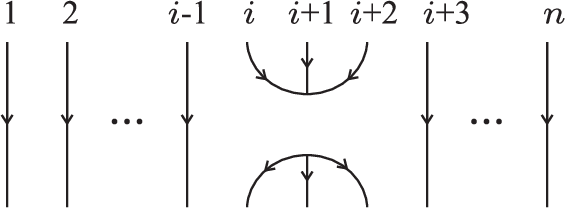}\\
 \caption{The $n$-tangle $F_i$, $i=1,\ldots, n-2$.}\label{fig:F_i}
\end{center}
\end{figure}

Let $F_i$ be the $m$-tangle illustrated in Figure \ref{fig:F_i}, and define $f_i = F_i + I_m$ so that $f_i = w_i w_{i+1} w_i - w_i = w_{i+1} w_i w_{i+1} - w_{i+1}$. By drawing pictures, it is easy to see that
$$f_i f_{i \pm 1} f_i = \delta^2 f_i, \qquad
f_i f_{i+2} f_i = \delta f_i w_{i+3}, \quad
\textrm{ and } \quad
f_i f_{i-2} f_i = \delta f_i w_{i-2}.$$
We also find that the $w_i$ satisfy the $SU(3)$ relation (\ref{eqn:SU(3)q_condition}):
$$(w_i - w_{i+2} w_{i+1} w_i + w_{i+1}) f_{i+1} = 0.$$

The following lemma is found in \cite[Lemma 3.3, p.385]{ohtsuki/yamada:1997}:

\begin{Lemma}\label{Lemma:ohtsuki/yamada_lemma}
Let $T$ be a basis $A_2$-$(m,n)$-tangle. Then $T$ must satisfy one of the following three conditions:
\begin{itemize}
\item[$(1)$] There are two consecutive vertices along the top which are connected by a cup or whose strings are joined at an (incoming) trivalent vertex,
\item[$(2)$] There are two consecutive vertices along the bottom which are connected by a cap or whose strings are joined at an (outgoing) trivalent vertex,
\item[$(3)$] $T$ is the identity tangle.
\end{itemize}
\end{Lemma}

Thus for any basis $A_2$-$m$-tangle which is not the identity tangle, there must be two (consecutive) vertices along the top or bottom whose strings are joined at an incoming or outgoing trivalent vertex respectively. In fact, by a Euler characteristic argument, this must be true for two vertices along both the top and bottom.

\begin{figure}[tb]
\begin{center}
  \includegraphics[width=45mm]{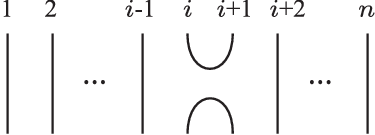}\\
 \caption{The $n$-diagram $E_i$, $i=1,\ldots, n-1$}\label{Fig:fig_Ei}
\end{center}
\end{figure}

Then we have the following lemma which says that the $A_2$-Temperley-Lieb algebra $A_2\textrm{-}TL$ is equal to the algebra $V^{A_2}$ of all $A_2$-tangles subject to the relations K1-K3. This is the $A_2$ analogue of the fact that the Temperley-Lieb algebra $TL_n = \mathrm{alg}(1,e_1,e_2,\ldots,e_{n-1})$ is isomorphic to Kauffman's diagram algebra \cite{kauffman:1987}, which is the algebra generated by the elements $E_1, E_2,\ldots, E_{n-1}$ on $n$ strings, illustrated in Figure \ref{Fig:fig_Ei}, along with the identity tangle $\mathbf{1}_n$ where every vertex along the top is connected to a vertex along the bottom by a vertical through string.
This lemma appeared in \cite{pugh:2008}, and also independently with an alternate proof in \cite[Theorem 2.2]{pylyavskyy:2007}.

\begin{figure}[htb]
\begin{center}
 \includegraphics[width=40mm]{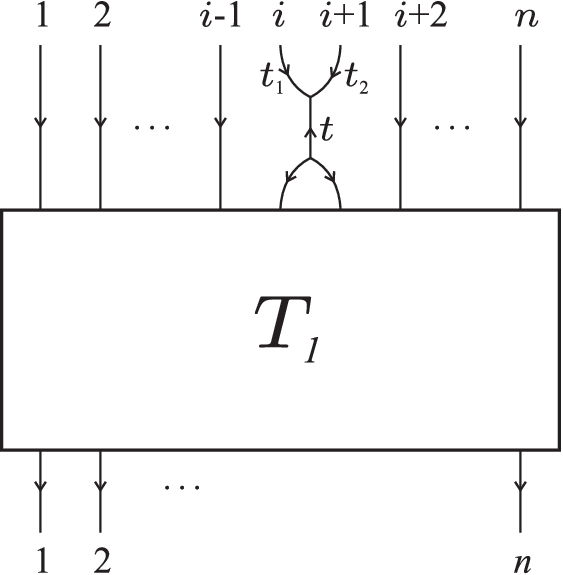}\\
 \caption{}\label{fig:fig8}
\end{center}
\end{figure}

\begin{Lemma}\label{Lemma:tangles_as_Wi's}
The algebra $V^{A_2}_m$ is generated by $\mathbf{1}_m$ and $W_i \in V^{A_2}_m$, $i=1, \ldots, m-1$. So $V^{A_2}_m \cong A_2\textrm{-}TL_m$.
\end{Lemma}
\emph{Proof:}
Let $T$ be a basis $m$-tangle which is not the identity. Then by Lemma \ref{Lemma:ohtsuki/yamada_lemma}, $T$ has (at least) one pair of vertices along the top whose strings are connected at an incoming trivalent vertex.
For an incoming trivalent which is only connected to two vertices along the top, the third strand of this trivalent vertex must be connected to an outgoing trivalent vertex, since it cannot be connected to another incoming trivalent vertex or a vertex along the bottom due to its orientation. Suppose these two vertices along the top are consecutive vertices. We isotope the strings so that we pull out this pair of trivalent vertices from the rest of the tangle as shown in Fig \ref{fig:fig8}, where $T_1$ is the resulting $m$-tangle contained inside the rectangle.
We repeat this procedure for all incoming trivalent vertices connected to exactly two vertices along the top, where these two vertices are consecutive. We also perform a similar procedure for all outgoing trivalent vertices connected to exactly two vertices along the bottom, where these two vertices are consecutive.

For any remaining trivalent vertices with only two of its strands connected to vertices along the top, these two vertices must not be consecutive. The region bounded by these two strands and the top of the tangle is a closed region which contains a non-zero number of vertices (in fact this number must necessarily be a multiple of three).
The braiding is a linear combination of the identity tangle and $W_i$'s. Thus by composing with the braiding we can move the pair of vertices along the top to the left side of the tangle so that these two vertices are consecutive. The strings may be isotoped in such a way so that once the braided part along the top has been removed to give a linear combination of the identity tangle and $W_i$'s, the resulting diagram does not contain any crossings. The third strand at this incoming trivalent vertex must again be connected to an outgoing trivalent vertex, and we pull out this pair of vertices as before, giving a factor of $W_1$. We repeat this procedure and the one described above for all the remaining incoming trivalent vertices connected to exactly two vertices along the top, and similarly for all the remaining outgoing trivalent vertices connected to exactly two vertices along the bottom.

If the resulting tangle is not the identity, then by Lemma \ref{Lemma:ohtsuki/yamada_lemma} there will again be a pair of vertices along the top whose strings are connected at an incoming trivalent vertex. Since all the incoming trivalent vertices which are connected to exactly two vertices along the top have been removed, this trivalent vertex must have all its strands connected to vertices along the top. By a similar argument there will also be an outgoing vertex which is connected to three vertices along the bottom. Then using the braiding we move this pair of trivalent vertices to the left of the diagram, which gives a factor $F_1$. Repeating this procedure we remove all the remaining trivalent vertices in the tangle, and we are done.
\hfill
$\Box$

\subsection{Trace on $\mathcal{V}^{A_2}_n$}

The following proposition is from \cite[Prop. 1.2, p.375]{ohtsuki/yamada:1997}:
\begin{Prop} \label{Lemma:ohtsuki/yamada_Prop-1.2}
The quotient $V^{A_2}_0 = A_2\textrm{-}TL_0$ of the free vector space of all planar 0-tangles by the Kuperberg relations K1-K3 is isomorphic to $\mathbb{C}$.
\end{Prop}

\begin{figure}[tb]
\begin{minipage}[t]{7cm}
 \begin{center}
  \includegraphics[width=30mm]{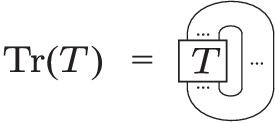}\\
  \caption{$\mathrm{Tr}(T)$} \label{fig:trace-tangle}
 \end{center}
\end{minipage}
\hfill
\begin{minipage}[t]{7cm}
 \begin{center}
  \includegraphics[width=35mm]{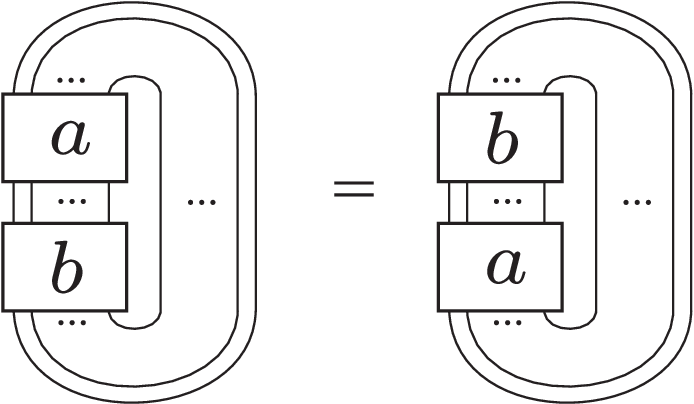}\\
 \caption{$\mathrm{Tr}(ab) = \mathrm{Tr}(ba)$}\label{fig:tr(ab)=tr(ba)}
 \end{center}
\end{minipage}
\end{figure}

\begin{figure}[b]
\begin{center}
  \includegraphics[width=120mm]{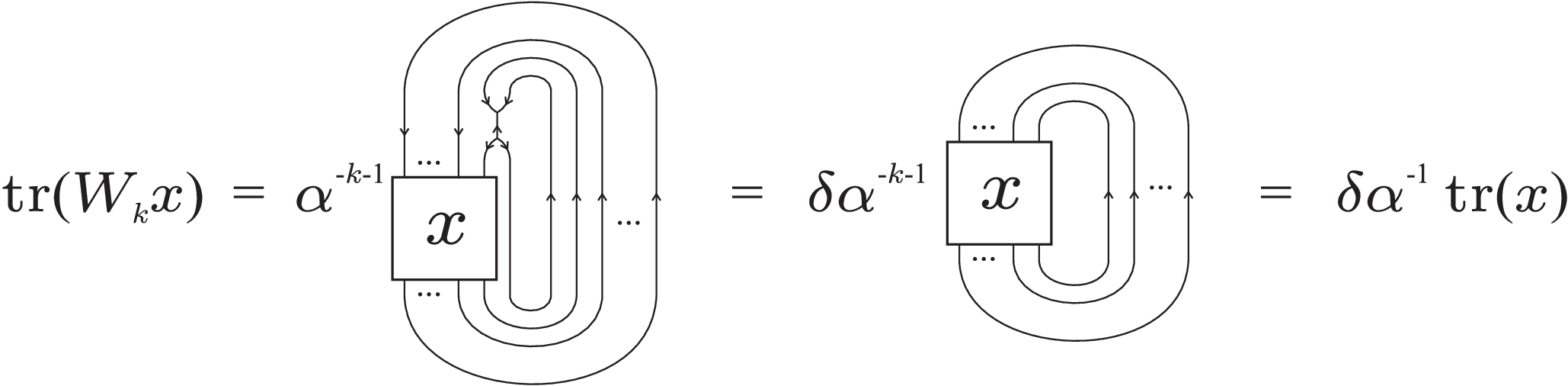}\\
 \caption{Markov trace on $V^{A_2}$}\label{fig:Markov_tr_for_Wi}
\end{center}
\end{figure}

We define a trace $\mathrm{Tr}$ on $\mathcal{V}^{A_2}_m$ as follows. For an $A_2$-$m$-tangle $T \in \mathcal{V}^{A_2}_m$, we form the 0-tangle $\mathrm{Tr}(T)$ as in Figure \ref{fig:trace-tangle} by joining the last vertex along the top of $T$ to the last vertex along the bottom by a string which passes round the tangle on the right hand side, and joining the other vertices along the top to those on the bottom similarly. Then $\mathrm{Tr}(T)$ gives a value in $\mathbb{C}$ by Proposition \ref{Lemma:ohtsuki/yamada_Prop-1.2}. We could define the above trace as a right trace, and define a left trace similarly where the strings pass round the tangle on the left hand side. However, by the comments after Proposition \ref{Prop:flat=spherical}, the right and left traces are equal.
The trace of a linear combination of tangles is given by linearity. Clearly $\mathrm{Tr}(ab) = \mathrm{Tr}(ba)$ for any $a,b \in \mathcal{V}^{A_2}_m$, as in Figure \ref{fig:tr(ab)=tr(ba)}. For any $x \in I_m$ we have $\mathrm{Tr}(x) = 0$, which follows trivially from the definition of Tr.
Then Tr is well defined on $V^{A_2}_m$. We define a normalized trace $\mathrm{tr}$ on $\mathcal{V}^{A_2}_m$ by $\mathrm{tr} = \alpha^{-m} \mathrm{Tr}$, so that $\mathrm{tr}(\mathbf{1}_m) = 1$. Then $\mathrm{tr}$ is a Markov trace on $V^{A_2}$ since for $x \in V^{A_2}_k$, $\mathrm{tr}(W_k x) = \delta \alpha^{-1} \mathrm{tr}(x)$, as illustrated in Figure \ref{fig:Markov_tr_for_Wi}, and in particular $\mathrm{tr} (W_i) = \delta \alpha^{-1}$. The Markov trace $\mathrm{tr}$ is positive by Lemma \ref{Lemma:Vn=Hn(q)} and \cite[Theorem 3.6(b)]{wenzl:1988}.

For each non-negative integer $m$ we define an inner-product on $\mathcal{V}^{A_2}_m$ by
\begin{equation} \label{inner-product_for_V}
\langle S, T \rangle = \mathrm{tr} (T^{\ast} S),
\end{equation}
which is well defined on $V^{A_2}_m$ since $\mathrm{tr}$ is.

For $\delta < 2$ (so $\delta = [2]_q = [2]$ where $q = e^{\pi i/n}$, $n \in \mathbb{N}$), we define $\widehat{V}^{A_2}_m$ to be the quotient of $V^{A_2}_m$ by the zero-length vectors in $V^{A_2}_m$ with respect to the inner-product defined in (\ref{inner-product_for_V}). Then the following lemma gives an identification between (a subalgebra of) the algebra of $A_2$-tangles and $\rho(H_{\infty}(q))$ where $\rho$ is one of Wenzl's Hecke representations for $SU(3)$ (see \cite{wenzl:1988}). This lemma will be used later in Section \ref{sect:planar_alg_for_A=PTL}.

\begin{Lemma} \label{Lemma:Vn=Hn(q)}
For $\delta \geq 2$, there is a $C^{\ast}$ representation $\rho$ of $H_{\infty}(q^2)$ such that $\rho(H_m(q^2)) \cong V^{A_2}_m$. The representation $\rho$ is equivalent to Wenzl's representation $\pi$ of the Hecke algebra, and consequently $V^{A_2}$ is isomorphic to the path algebra for $\mathcal{A}^{(\infty)}$. For $\delta = [2]_q$, $q = e^{\pi i/n}$, there is a $C^{\ast}$ representation $\rho$ of $H_{\infty}(q^2)$ such that $\rho(H_m(q^2)) \cong \widehat{V}^{A_2}_m$. In this case the representation $\rho$ is equivalent to Wenzl's representation $\pi^{(3,n)}$ of the Hecke algebra, and consequently $V^{A_2}$ is isomorphic to the path algebra for $\mathcal{A}^{(n)}$.
\end{Lemma}
\emph{Proof:}
Clearly $\delta^{-1} W_i$, $i=1,\ldots,m-1$, is a self-adjoint projection in $V^{A_2}_m$, and hence $\rho$ is a $C^{\ast}$-representation of $H_m(q^2)$ for any real $q \geq 1$ or $q = e^{\pi i/n}$. When $q = e^{x}$, $x \geq 0$, we have $\eta = (1-q^{2(-k+1)})/(1+q^2)(1-q^{-2k}) = \sinh((k-1)x)/2\cosh(x)\sinh(kx) = [k-1]_q/[2]_q[k]_q$, whilst for $q = e^{\pi i/n}$, $\eta = \sin((k-1)\pi/n)/2\cos(\pi/n)\sin(k\pi/n) = [k-1]/[2][k]$. Then for $k = 3$, $\eta = [3]_q^{-1}$ so that the Markov trace on $V^{A_2}_m$ satisfies the condition in \cite[Theorem 3.6]{wenzl:1988}.
\hfill
$\Box$

Then the algebra $V^{A_2}_m$ is finite-dimensional for all finite $m$ since the $m^{\textrm{th}}$ level of the path algebra for $\mathcal{A}^{(n)}$ is finite-dimensional.

\section{$A_2$-planar algebras} \label{sec:A_2-planar_algebras}
\subsection{General $A_2$-planar algebras} \label{sec:general_A_2-planar_algebras}

We will now define an $A_2$-version of Jones' planar algebra, using tangles generated by Kuperberg's $A_2$-webs. Under certain assumptions, these $A_2$-planar algebras will correspond to certain subfactors of $SU(3)$ $\mathcal{ADE}$ graphs which have flat connections.
The best way to describe planar algebras is in terms of operads (see \cite{jones:planar, may:1997}).

\begin{Def} \label{Def:operad}
An \textbf{operad} consists of a sequence $(\mathcal{C}(n))_{n \in \mathbb{N}}$ of sets. There is a unit element 1 in $\mathcal{C}(1)$, and a function $\mathcal{C}(n) \otimes \mathcal{C}(j_1) \otimes \cdots \otimes \mathcal{C}(j_n) \rightarrow \mathcal{C}(j_1 + \cdots + j_n)$ called composition, given by $(y \otimes x_1 \otimes \cdots \otimes x_n) \rightarrow y \circ (x_1 \otimes \cdots \otimes x_n)$, satisfying the following properties
\begin{itemize}
\item associativity: $y \circ (x_1 \circ (x_{1,1} \otimes \cdots \otimes x_{1,k_1}) \otimes \cdots \otimes x_n \circ (x_{n,1} \otimes \cdots \otimes x_{n,k_n}))$ \\
    ${} \qquad \qquad \qquad = (y \circ (x_1 \otimes \cdots \otimes x_n)) \circ (x_{1,1} \otimes \cdots \otimes x_{1,k_1} \otimes \cdots \otimes x_{n,1} \otimes \cdots \otimes x_{n,k_n})$,
\item identity: $y \circ (1 \otimes \cdots \otimes 1) = y = 1 \circ y$.
\end{itemize}
\end{Def}

\begin{figure}[htb]
\begin{center}
  \includegraphics[width=70mm]{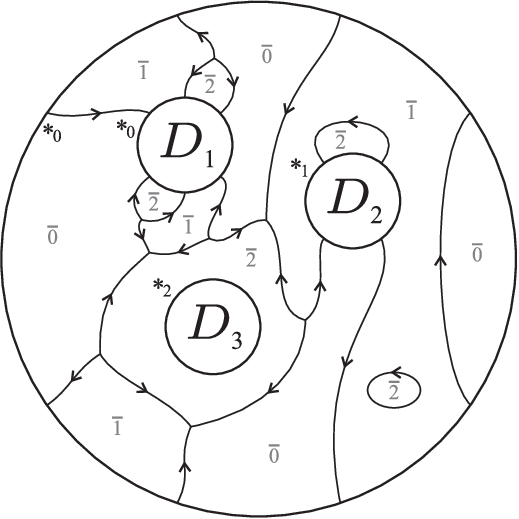}\\
 \caption{$A_2$-planar $\sigma$-tangle for $\sigma = -+-+-+-+$}\label{fig:A2_planar_tangle}
\end{center}
\end{figure}

Let $\sigma = \sigma_1 \cdots \sigma_m$ be a sign string, $\sigma_j \in \{ \pm \}$. An $A_2$-planar $\sigma$-tangle will be the unit disc $D=D_0$ in $\mathbb{C}$ together with a finite (possibly empty) set of disjoint sub-discs $D_1, D_2, \ldots, D_n$ in the interior of $D$. Each disc $D_k$, $k \geq 0$, will have $m_k \geq 0$ vertices on its boundary $\partial D_k$, whose orientations are determined by sign strings $\sigma^{(k)} = \sigma^{(k)}_1 \cdots \sigma^{(k)}_{m_k}$ where `$+$' denotes a sink and `$-$' a source, and such that the difference between the number of `$+$' and `$-$' is 0 mod 3. The disc $D_k$ will be said to have pattern $\sigma^{(k)}$.
Inside $D$ we have an $A_2$-tangle where the endpoint of any string is either a trivalent vertex (see Figure \ref{fig:A2-webs}) or one of the vertices on the boundary of a disc $D_k$, $k=0, \ldots, n$, or else the string forms a closed loop. Each vertex on the boundaries of the $D_k$ is the endpoint of exactly one string, which meets $\partial D_k$ transversally.
An example of an $A_2$-planar $\sigma$-tangle is illustrated in Figure \ref{fig:A2_planar_tangle} for $\sigma = -+-+-+-+$.

The regions inside $D$ have as boundaries segments of the $\partial D_k$ or the strings. These regions ar labelled $\overline{0}$, $\overline{1}$ or $\overline{2}$, called the colouring, such that if we pass from a region $R$ of colour $\overline{a}$ to an adjacent region $R'$ by passing to the right over a vertical string with downwards orientation, then $R'$ has colour $\overline{a+1}$ (mod 3). We mark the segment of each $\partial D_k$ between the last and first vertices with $\ast_{b_k}$, $b_k \in \{0,1,2\}$, so that the region inside $D$ which meets $\partial D_k$ at this segment is of colour $\overline{b_k}$, and the choice of these $\ast_{b_k}$ must give a consistent colouring of the regions.
For each $\sigma$ we have three types of tangle, depending on the colour $\overline{b}$ of the marked segment, or of the marked region near $\partial D$ for $\sigma = \varnothing$.

We define $\widetilde{\mathcal{P}}_{\sigma}(L)$ to be the free vector space generated by orientation-preserving diffeomorphism classes of $A_2$-planar $\sigma$-tangles with labelling sets $L$. The diffeomorphisms preserve the boundary of $D$, but may move the $D_k$'s, $k \geq 1$. Let $\mathcal{P}_{\sigma}(L)$ be the quotient of $\widetilde{\mathcal{P}}_{\sigma}(L)$ by the Kuperberg relations K1-K3. The \textbf{$A_2$-planar operad} $\mathcal{P}(L)$ is defined to be $\mathcal{P}(L) = \bigcup_{\sigma} \mathcal{P}_{\sigma}(L)$. We will usually simply write $\mathcal{P}$ for $\mathcal{P}(L)$.

We define composition in $\mathcal{P}$ as follows. Given an $A_2$-planar $\sigma$-tangle $T$ with an internal disc $D_l$ with pattern $\sigma_l = \sigma'$, and an $A_2$-planar $\sigma'$-tangle $S$ with external disc $D'$ and $\ast_{D'}=\ast_{D_l}$, we define the $\sigma$-tangle $T \circ_l S$ by isotoping $S$ so that its boundary and vertices coincide with those of $D_l$, joining the strings at $\partial D_l$ and smoothing if necessary. We then remove $\partial D_l$ to obtain the tangle $T \circ_l S$ whose diffeomorphism class clearly depends only on those of $T$ and $S$. This gives $\mathcal{P}$ the structure of a \textbf{coloured operad}, where each $D_k$, $k>0$, is assigned the colour $\sigma_k$, and composition is only allowed when the colouring of the regions match (which forces the orientations of the vertices to agree). The $D_k$'s, $k \geq 1$ are to be thought of as inputs, and $D=D_0$ is the output.

\begin{figure}[tb]
\begin{center}
  \includegraphics[width=130mm]{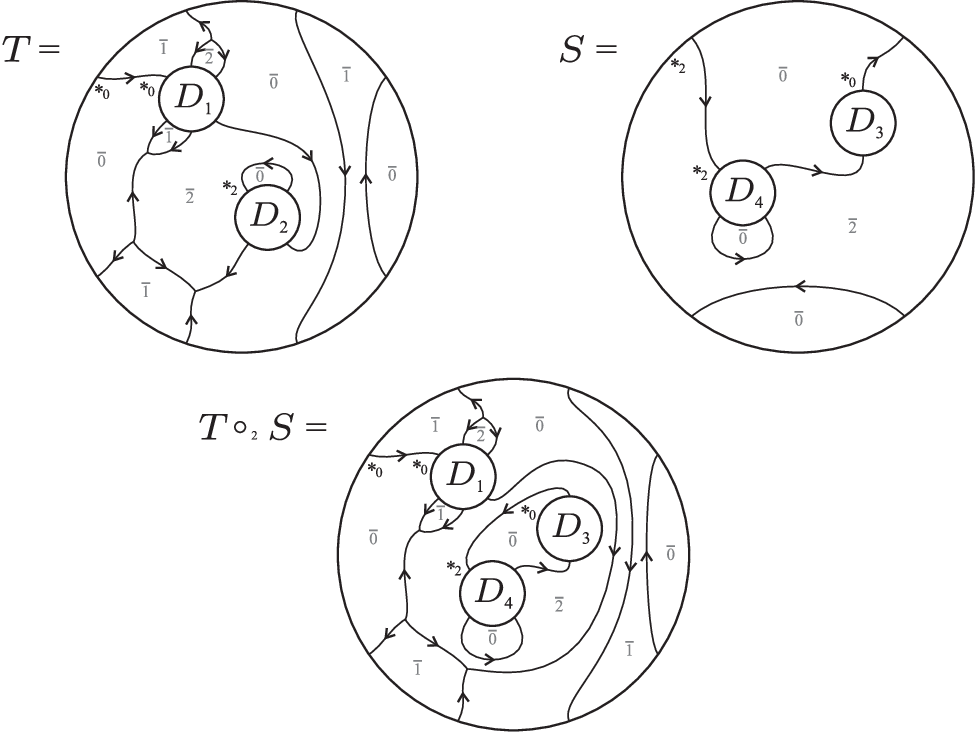}\\
 \caption{Composition of planar tangles}\label{fig:composition_of_planar_tangle}
\end{center}
\end{figure}

The most general notion of an $A_2$-planar algebra will be an algebra over the operad $\mathcal{P}$, i.e. a \textbf{general $A_2$-planar algebra} $P$ is a family
$$P = \left( P_{\sigma}^{\overline{a}}, \textrm{ for all sign strings } \sigma, \textrm{ and all } a \in \{0,1,2\} \right)$$
of vector spaces
with the following property: for every labelled $\sigma$-tangle $T \in \mathcal{P}_{\sigma}$ with internal discs $D_1, D_2, \ldots, D_n$, where $D_k$ has pattern $\sigma_k$ and outer disc marked by $\ast_{b_k}$, there is associated a linear map $Z(T): \otimes_{k=1}^n P_{\sigma_k}^{\overline{b_k}} \longrightarrow P_{\sigma}^{\overline{b}}$ which is compatible with the composition of tangles in the following way. If $S$ is a $\sigma_k$-tangle with internal discs $D_{n+1}, \ldots, D_{n+m}$, where $D_k$ has pattern $\sigma_k$, then the composite tangle $T \circ_l S$ is a $\sigma$-tangle with $n+m-1$ internal discs $D_k$, $k = 1,2, \ldots l-1, l+1, l+2, \ldots, n+m$. From the definition of an operad, associativity means that the following diagram commutes:
\begin{equation} \label{eqn:compatability_condition_for_Z(T)}
\xymatrix{
{\left( \bigotimes_{\stackrel{k=1}{\scriptscriptstyle{k \neq l}}}^n P_{\sigma_k}^{\overline{b_k}} \right) \otimes \left( \bigotimes_{k=n+1}^{n+m} P_{\sigma_k}^{\overline{b_k}} \right)} \ar[d]_{\mathrm{id} \otimes Z(S)} \ar[dr]^(.6){Z(T \circ_l S)} \\
{\bigotimes_{k=1}^n P_{\sigma_k}}^{\overline{b_k}} \ar[r]_{Z(T)} & P_{\sigma}^{\overline{b}} }
\end{equation}
so that $Z(T \circ_l S) = Z(T')$, where $T'$ is the tangle $T$ with $Z(S)$ used as the label for disc $D_l$. We also require $Z(T)$ to be independent of the ordering of the internal discs, that is, independent of the order in which we insert the labels into the discs.
If $\sigma = \varnothing$, we adopt the convention that the empty tensor product is the complex numbers $\mathbb{C}$. By using the tangle
\begin{center}
  \includegraphics[width=25mm]{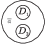}
\end{center}
we see that each $P_{\varnothing}^{\overline{a}}$ (sometimes denoted by $P_{0}^{\overline{a}}$) is a commutative associative algebra, $a \in \{0,1,2\}$. Each $P_{\sigma}^{\overline{a}}$ has a distinguished subset, given by the elements $Z(T)$ for all $\sigma$-tangles without internal discs, with outer disc marked by $\ast_{a}$. This is the unital operad (see \cite{may:1997}). Following Jones' terminology, we call the linear map $Z$ the \textbf{presenting map} for $P$.

Jones' planar algebra is contained in the $A_2$-planar algebra in the following way.
Let $(\pm,n)$ denote the alternating sign string of length $n$, where the first sign is $\pm$. If we consider the sub-operad $\mathcal{Q} = \bigcup \mathcal{Q}_n$ where $\mathcal{Q}_n$ is the subset of $\mathcal{P}_{(\pm,n)}$ generated by tangles with no trivalent vertices (and hence no crossings) and where each internal disc $D_k$ only has pattern $(\pm,n_k)$, then $\mathcal{Q}$ is the coloured planar operad of Jones in \cite{jones:planar}, where instead of the three colours $a = 0,1,2$ of the $A_2$-planar algebras, in $\mathcal{Q}$ there are now only two colours, usually called black and white. Jones' planar algebra is then $Q = Z(\mathcal{Q})$.

\subsection{Partial Braiding} \label{sect:braiding}

\begin{figure}[htb]
\begin{center}
  \includegraphics[width=51mm]{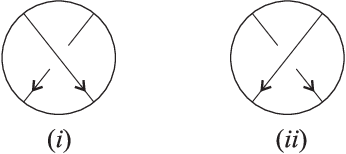}\\
 \caption{Disc $b$ for $(i)$ negative crossing, $(ii)$ positive crossing}\label{fig:disc_b}
\end{center}
\end{figure}

We now introduce the notion of a partial braiding in our $A_2$-planar operad. We will allow over and under crossings in our diagrams, which are interpreted as follows. For a tangle $T$ with $n$ crossings $c_1, \ldots, c_n$, choose one of the crossings $c_i$ and, isotoping any strings if necessary, we enclose $c_i$ in a disc $b$, as shown in Figure \ref{fig:disc_b} for $c_i$ a $(i)$ negative crossing and $(ii)$ positive crossing (up to some rotation of the disc).

\begin{figure}[htb]
\begin{center}
  \includegraphics[width=65mm]{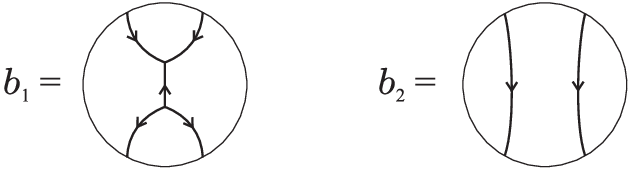}\\
 \caption{Discs $b_1$ and $b_2$}\label{fig:discs_b1_and_b2}
\end{center}
\end{figure}

Let $b_1$, $b_2$ be the discs illustrated in Figure \ref{fig:discs_b1_and_b2}. We form two new tangles $S_1^{(1)}$ and $T_1^{(1)}$ which are identical to $T$ except that we replace the disc $b$ by $b_1$ for $S_1^{(1)}$ and by $b_2$ for $T_1^{(1)}$. If $c_i$ is a negative crossing then $T$ is equal to the linear combination of tangles $q^{-2/3} S_1^{(1)} - q^{1/3} T_1^{(1)}$, and if $c_i$ is a positive crossing $T = q^{2/3} S_1^{(1)} - q^{-1/3} T_1^{(1)}$, where $q > 0$ satisfies $q + q^{-1} = \delta$ (cf. (\ref{braiding1}) and (\ref{braiding2})). Then for both $S_1^{(1)}$ and $T_1^{(1)}$ we consider another crossing $c_j$ and repeat the above process to obtain $S_1^{(1)} = r_1 S_1^{(2)} - r_1' T_1^{(2)}$, $T_1^{(1)} = r_2 S_2^{(2)} - r_2' T_2^{(2)}$, where $r_1, r_2 \in \{ q^{\pm 2} \}$ and $r_1', r_2' \in \{ q^{\pm 1} \}$ depending on whether $c_j$ is a positive or negative crossing. Since this expansion of the crossings is independent of the order in which the crossings are selected, repeating this procedure we obtain a linear combination $T = \sum_{i=1}^{2^{(n-1)}} (s_i S_i^{(n)} + s_i' T_i^{(n)})$, where the $s_i$, $s_i'$ are powers of $q^{\pm 1/3}$.

With this definition of a partial braiding, two tangles give identical elements of the planar algebra if one can be deformed into the other using relations (\ref{braiding_relations1-3}), (\ref{braiding_relation4}). It is not a braiding as we cannot in general pull strings over or under labelled inner discs $D_k$.

\begin{figure}[htb]
\begin{center}
  \includegraphics[width=25mm]{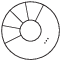}\\
 \caption{Tangle $I_{\sigma}$} \label{fig:I_sigma}
\end{center}
\end{figure}

The tangles $I_{\sigma} \in \mathcal{P}_{\sigma}$ illustrated in Figure \ref{fig:I_sigma} have pattern $\sigma$ on the inner and outer discs and all strings are through strings. For any $\sigma$-tangle $T$ these tangles satisfy $I_{\sigma} \circ T = T$, and also inserting $I_{\sigma_k}$ inside every inner disc $D_k$ with pattern $\sigma_k$ also gives the original tangle $T$. Then $I_{\sigma}$ is the unit element (see Definition \ref{Def:operad}). We let $I_{\sigma}(x)$ denote the tangle $I_{\sigma}$ with $x \in P_{\sigma}$ as the label for the inner disc.

The condition $\mathrm{dim}(P_0^{\overline{a}})=1$, $a=0,1,2$, implies that there is a unique way to identify each $P_0^{\overline{a}}$ with $\mathbb{C}$ as algebras, with $Z \left( \bigcirc_a \right) = 1$, $a=0,1,2$, where $\bigcirc_a$ is the empty tangle with no vertices or strings at all, with the interior coloured $a$. By Lemma \ref{Lemma:delta^2=alpha+1} there is thus also one scalar, or parameter, associated to a general $A_2$-planar algebra:
\setlength{\unitlength}{1mm}
\begin{equation}\label{fig:Z(O)=alpha}
Z(\begin{picture}(6,5)
\put(3,1){\circle{5}}
\put(3,1){\circle{2}}
\end{picture}) = \alpha,
\end{equation}
where the inner circle is a closed loop not an internal disc.

It follows from the compatability condition (\ref{eqn:compatability_condition_for_Z(T)}) that $Z$ is multiplicative on connected components, i.e. if a part of a tangle $Y$ can be surrounded by a disc so that $T = T' \circ_l S$ for a tangle $T'$ and 0-tangle $S$, then $Z(T)=Z(S)Z(T')$ where $Z(S)$ is a multilinear map from $\mathcal{P}_{0}^{\overline{a}}$ into the field $\mathbb{C}$, where the region which meets the outer boundary of $S$ is coloured $a$, $a \in \{ 0,1,2 \}$.

Every general $A_2$-planar algebra contains the $A_2$-planar subalgebra $PTL$, the \textbf{planar $A_2$-Temperley-Lieb algebra}, which is defined by $PTL_{\sigma} = \mathcal{P}_{\sigma}(\varnothing)$, i.e. there is no labelling set. We have $PTL_{0}^{\overline{a}} \cong \mathbb{C}$. The presenting map $Z$ is just the identity map. Note that the partial braiding defined above is a genuine braiding in $PTL$. The $A_2$-Temperley-Lieb algebra, introduced in section \ref{sec:A_2-tangles}, is a subalgebra of $PTL$, given by $A_2\textrm{-}TL_n = PTL_{-^n +^n}$, where $+^n$ denotes the sign string $++ \cdots +$ ($n$ copies), and $-^n = -- \cdots -$ ($n$ copies). The action of an $A_2$-planar $\sigma$-tangle $T$ on $PTL$ is given by filling the internal discs of $T$ with basis elements of $PTL$, where we ignore the colouring of the regions in $T$. The resulting tangle may then contain digons or embedded squares, which are removed using K2 and K3, and closed curves are removed using (\ref{fig:Z(O)=alpha}). The result is a linear combination of elements of $PTL$. In the $A_1$ case, the planar algebra for which there is no labelling set is the Temperley-Lieb algebra itself, $\mathcal{P}_n(\varnothing) = TL_n$.

Suppose $\sigma$ is a sign string. We define $\sigma^{\ast}$ to be the sign string obtained by reversing the string $\sigma$ and flipping all its signs.

We define multiplication tangles $M_{\sigma\sigma^{\ast}}: \mathcal{P}_{\sigma\sigma^{\ast}} \times \mathcal{P}_{\sigma\sigma^{\ast}} \rightarrow \mathcal{P}_{\sigma\sigma^{\ast}}$ by:
\begin{center}
  \includegraphics[width=47mm]{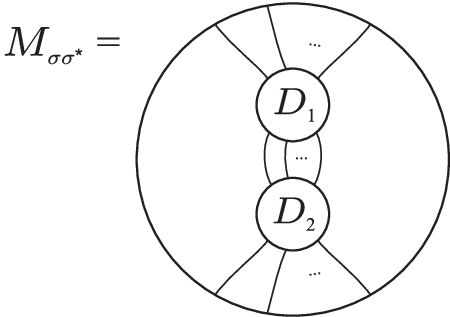}
\end{center}

Each $P_{\sigma\sigma^{\ast}}$ is then an associative algebra, with multiplication being defined by $x_1 x_2 = Z(M_{\sigma\sigma^{\ast}}(x_1, x_2))$, where $M_{\sigma\sigma^{\ast}}(x_1, x_2)$ has $x_k \in P_{\sigma\sigma^{\ast}}$ as the insertion in disc $D_k$, $k=1,2$. The multiplication is also clearly compatible with the inclusion tangles, as can be seen by drawing pictures.

\begin{figure}[tb]
\begin{minipage}[t]{7.5cm}
\begin{center}
  \includegraphics[width=43mm]{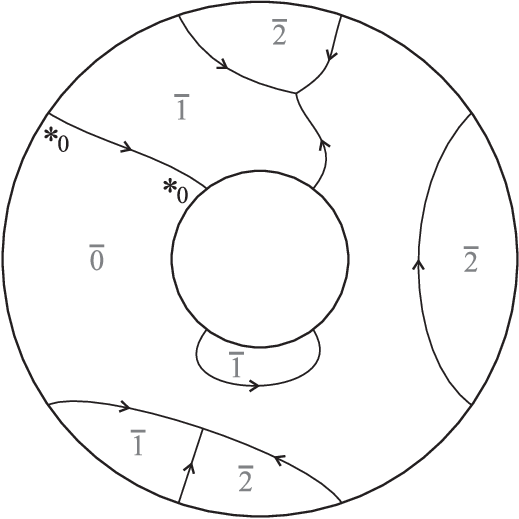}\\
 \caption{Annular tangle} \label{fig:annular_tangle}
\end{center}
\end{minipage}
\hfill
\begin{minipage}[t]{7.5cm}
\begin{center}
  \includegraphics[width=40mm]{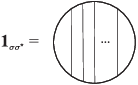}\\
 \caption{Identity Tangle $\mathbf{1}_{\sigma\sigma^{\ast}} \in \mathcal{P}_{\sigma\sigma^{\ast}}$} \label{fig:identity_tangle}
\end{center}
\end{minipage}
\end{figure}

An annular tangle with outer disc with pattern $\sigma$ and inner disc with pattern $\sigma'$ will be called an annular $(\sigma,\sigma')$-tangle. An example of an annular $(\sigma,\sigma')$-tangle is illustrated in Figure \ref{fig:annular_tangle}, where $\sigma = ---+-+++$, $\sigma' = -+-+$.

\begin{figure}[tb]
\begin{center}
  \includegraphics[width=100mm]{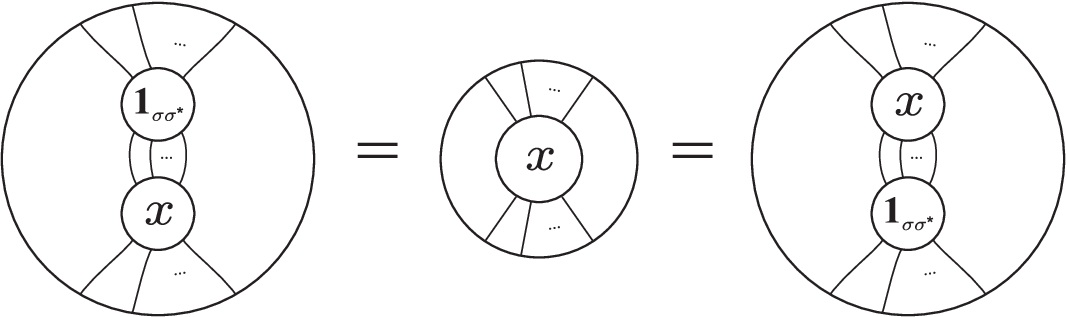}\\
 \caption{$Z(\mathbf{1}_{\sigma\sigma^{\ast}}) x = x = x Z(\mathbf{1}_{\sigma\sigma^{\ast}})$} \label{fig:1x=x=x1}
\end{center}
\end{figure}

The tangle $\mathbf{1}_{\sigma\sigma^{\ast}} \in \mathcal{P}_{\sigma\sigma^{\ast}}$ illustrated in Figure \ref{fig:identity_tangle} is called the identity tangle. By inserting $\mathbf{1}_{\sigma\sigma^{\ast}}$ and $x \in P_{\sigma\sigma^{\ast}}$ into the discs of the multiplication tangle $M_{\sigma\sigma^{\ast}}$ as in Figure \ref{fig:1x=x=x1} we see that $Z(\mathbf{1}_{\sigma\sigma^{\ast}}) x = x = x Z(\mathbf{1}_{\sigma\sigma^{\ast}})$, hence $Z(\mathbf{1}_{\sigma\sigma^{\ast}})$ is the left and right identity for $P_{\sigma\sigma^{\ast}}$.

The following proposition shows that the $A_2$-planar operad $\mathcal{P}$ is generated by the algebra $PTL$, multiplication tangles $M$, and annular tangles, which are tangles with only one internal disc. We note that this result is only one possible choice for the generators of the $A_2$-planar operad and that there is much freedom in the choice of such generators.

\begin{figure}[h]
\begin{center}
  \includegraphics[width=30mm]{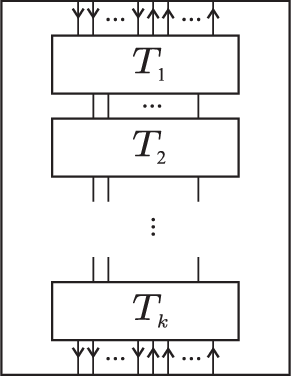}\\
 \caption{An arbitrary tangle $T \in \mathcal{P}_{\sigma\sigma^{\ast}}$, for $\sigma = -^k +^{k'}$} \label{fig:Prop2.4(1)}
\end{center}
\end{figure}

\begin{Prop}
The $A_2$-planar operad $\mathcal{P}$ is generated by the algebra $PTL$, multiplication tangles $M$, and annular tangles.
\end{Prop}
\emph{Proof:}
Consider first an arbitrary tangle $T \in \mathcal{P}_{\sigma\sigma^{\ast}}$ which has $k$ inner discs $D_l$ with labels $x_l$, $l=1,\ldots,k$, and where the sign string $\sigma$ is of the form $-^k +^{k'}$ (we can always insert the tangle $T$ inside an annular tangle which uses the braiding to permute the vertices if $\sigma$ is not of this form). We isotope the tangle to move all the inner discs so that the tangle can be divided into horizontal strips in such a way that in any horizontal strip there is only one disc. Then we may draw $T$ as in Figure \ref{fig:Prop2.4(1)}, where the $T_l$ are all tangles with one inner disc labelled by $x_l$, $l=1,\ldots,k$, and where we draw the tangles inside rectangles rather than discs.

Consider first the tangle $T_1$, which has pattern $\sigma$ along the top edge. Using the braiding we may permute all the strings along the bottom of $T_1$ so that they are of the form $-^{k_1} +^{k_1'}$ (reading from left to right), i.e. all the strings with downwards orientation are moved to the left.
Now $k+k_1' \equiv k'+k_1 \textrm{ mod } 3$, so we have $k-k' = k_1-k_1'+3p$, for some $p \in \mathbb{Z}$.
Suppose $p>0$.
Then we add $p$ double loops $\;$ \includegraphics[width=5mm]{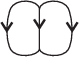} $\;$ at the bottom of $T_1$ to the left of the leftmost string (and multiply the tangle $T$ by a scalar factor $\alpha^{-p} \delta^{-p}$):
\begin{center}
  \includegraphics[width=77mm]{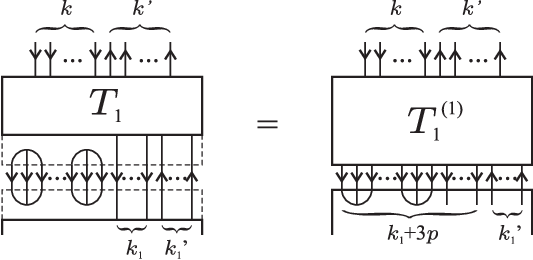}
\end{center}
If $p<0$ we instead add $p$ double loops at the top of $T_1$ to the right of the rightmost string, and similarly at the bottom of $T_k$ (and multiply $T$ by a scalar factor $\alpha^{-2p}\delta^{-2p}$).
We now have $k-k' = k_1-k_1'$, and the number of vertices along the top and bottom of $T_1^{(1)}$ differs by an even integer, i.e. $k+k' = k_1+k_1'+2p'$, for some $p' \in \mathbb{Z}$.
Suppose $p'>0$.
Then we add $p'$ concentric closed loops (with anti-clockwise orientation) beneath $T_1^{(1)}$, between the rightmost string with downwards orientation and the leftmost string with upwards orientation (and multiply the tangle $T$ by a scalar factor $\alpha^{-p'}$):
\begin{center}
  \includegraphics[width=70mm]{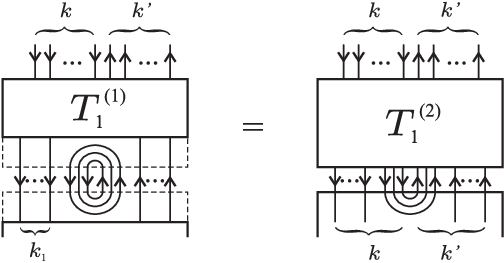}
\end{center}
If $p'<0$ we instead add $p'$ concentric closed loops (with clockwise orientation) above $T_1^{(1)}$, between the rightmost string with downwards orientation and the leftmost string with upwards orientation, and similarly at the bottom of $T_k$ (and multiply by a scalar factor $\alpha^{-2p'}$.
Then we have a multiplication tangle $M_{\widetilde{\sigma}\widetilde{\sigma}^{\ast}}$ surrounded by an annular $(\sigma\sigma^{\ast},\widetilde{\sigma}\widetilde{\sigma}^{\ast})$-tangle (where $\widetilde{\sigma}$ is possibly equal to $\sigma$), with $T_1^{(2)}$ as the insertion for the first disc of $M_{\widetilde{\sigma}\widetilde{\sigma}^{\ast}}$, and the rest of the tangle, which we will call $T'$, as the insertion for the second disc. So $T'$ is an $\widetilde{\sigma}\widetilde{\sigma}^{\ast}$-tangle with $k-1$ inner discs, and by the above procedure we can write $T'$ as a multiplication tangle (possibly surrounded by an annular tangle), where the insertion for the second disc now only has $k-2$ inner discs. Continuing in this way we see inductively that $T$ is generated by multiplication tangles and annular tangles.
Suppose now that $T \in \mathcal{P}_{\sigma}$, where $\sigma$ is not of the form $\widetilde{\sigma}\widetilde{\sigma}^{\ast}$ for some sign string $\widetilde{\sigma}$. By using a similar procedure to that given above we can write the tangle $T$ as $T' \in \mathcal{P}_{\widetilde{\sigma}\widetilde{\sigma}^{\ast}}$ surrounded by an annular $(\sigma,\widetilde{\sigma}\widetilde{\sigma}^{\ast})$-tangle.
Finally, tangles with no inner discs are elements of $PTL$.
\hfill
$\Box$

\begin{Def}
A general $A_2$-planar algebra $P$ will be called \textbf{finite-dimensional} if $\mathrm{dim}P_{\sigma} < \infty$ for all $\sigma$.
\end{Def}

\emph{Remark.} The algebras $A_2\textrm{-}TL_n$ are finite dimensional, since from section \ref{sec:A_2-tangles} we know that they are isomorphic to the path algebra for the $SU(3)$ graph $\mathcal{A}^{(\infty)}$. By Theorem 6.3 in \cite{kuperberg:1996} the dimensions of $PTL_{\sigma}$ and $PTL_{\sigma'}$ are the same for $\sigma'$ any permutation of $\sigma$. Thus $PTL_{\sigma}$ is finite dimensional for any $\sigma$ which is a permutation of $+^n -^n$. It follows from Corollary \ref{Cor:P-f.d.} at the end of Section \ref{sect:involution_on_P} that $PTL_{\sigma}$ is thus finite dimensional for all sign strings $\sigma$.

\subsection{$A_2$-Planar Algebras} \label{sect:A2-planar_algebras+flatness}

We now define an $A_2$-planar algebra $P$, where unlike for general $A_2$-planar algebras, there are restrictions on the dimensions of the lowest graded parts. The $A_2$-planar algebra $P$ comes with two traces. We will also define notions of non-degeneracy and sphericity in the same way as Jones \cite[Definition 1.27]{jones:planar}, and the notion of flatness.

\begin{Def} \label{Def:A_2-planar_algebra}
\begin{itemize}
\item[(a)]
An \textbf{$A_2$-planar algebra} will be a general $A_2$-planar algebra $P$ which has $\mathrm{dim} \left( P_0^{\overline{0}} \right) = \mathrm{dim} \left( P_0^{\overline{1}} \right) = \mathrm{dim} \left( P_0^{\overline{2}} \right) = 1$, and $Z(\begin{picture}(6,5)
\put(3,1){\circle{5}}
\put(3,1){\circle{2}}
\end{picture}) = \alpha$ non-zero.
\item[(b)]
We call the presenting map $Z$ the \textbf{partition function} when it is applied to a closed $0$-tangle $T$ with internal discs $D_k$ of pattern $\sigma_k$. We identify $P_0^{\overline{a}}$ with $\mathbb{C}$, so that $Z(T): \otimes_k P_{\sigma_k} \longrightarrow \mathbb{C}$.
\item[(c)]
Let $A_{\sigma}$ be the set of all $0$-tangles with only one internal disc, where the internal disc has pattern $\sigma$. An $A_2$-planar algebra will be called \textbf{non-degenerate} if, for $x \in P_{\sigma}$, $x=0$ if and only if $Z(T(x)) = 0$ for all $T \in A_{\sigma}$.
An $A_2$-planar algebra will be called \textbf{spherical} if its partition function is an invariant of tangles on the two-sphere $S^2$ (obtained from $\mathbb{R}^2$ by adding a point at infinity).
\item[(d)]
Let $P$ be an $A_2$-planar algebra, and $\sigma$ a sign string. Define two traces ${}_L \mathrm{Tr}_{\sigma\sigma^{\ast}}$ and ${}_R \mathrm{Tr}_{\sigma\sigma^{\ast}}$ on $P_{\sigma\sigma^{\ast}}$ by
\begin{center}
  \includegraphics[width=100mm]{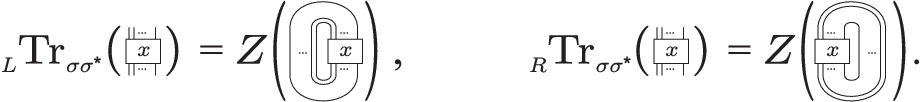}
\end{center}
\end{itemize}
\end{Def}

For a spherical $A_2$-planar algebra ${}_L \mathrm{Tr}_{\sigma\sigma^{\ast}} = {}_R \mathrm{Tr}_{\sigma\sigma^{\ast}} =: \mathrm{Tr}_{\sigma\sigma^{\ast}}$. The converse is also true- that is, if ${}_L \mathrm{Tr}_{\sigma\sigma^{\ast}} = {}_R \mathrm{Tr}_{\sigma\sigma^{\ast}}$ on $P_{\sigma\sigma^{\ast}}$ for all sign strings $\sigma$ then $P$ is spherical.

The proof of the following proposition given in \cite{jones:planar} in the setting of his $A_1$-planar algebras yields:

\begin{Prop}
A spherical $A_2$-planar algebra $P$ is non-degenerate if and only if $\mathrm{Tr}_{\sigma\sigma^{\ast}}$ defines a non-degenerate bilinear form on $P_{\sigma\sigma^{\ast}}$ for each sign string $\sigma$.
\end{Prop}

\begin{figure}[htb]
\begin{center}
  \includegraphics[width=70mm]{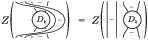}\\
 \caption{Flatness} \label{fig:flatness}
\end{center}
\end{figure}

\begin{Def} \label{def:flatness}
Let $T$ be any tangle with internal discs $D_k$, $k = 1, \ldots, n$. We call an $A_2$-planar algebra \textbf{flat} if $Z(T) = Z(T')$ where $T'$ is any tangle obtained from $T$ by pulling strings over an internal disc $D_k$, for any $k = 1, \ldots, n$. This is illustrated in Figure \ref{fig:flatness}, where we only show a local part of the tangle.
\end{Def}

We could alternatively have defined a flat $A_2$-planar algebra to be one where strings can be pulled under internal discs instead of over. Such an $A_2$-planar algebra is isomorphic to the one defined above, with the isomorphism given by replacing $q$ by $q^{-1}$, equivalent to reversing all crossings in any tangle.
Note that our definition of flatness does not imply that we can also pull strings under internal discs, which in general will not be the case -- c.f. the relative braiding notion in the theory of $\alpha$-induction as explained in \cite[Section 3.3]{bockenhauer/evans:1999ii} and \cite[Section 2]{bockenhauer/evans/kawahigashi:2000}.

\begin{figure}[tb]
\begin{center}
  \includegraphics[width=105mm]{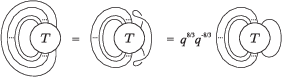}\\
 \caption{Flatness gives sphericity} \label{fig:flatness_gives_sphericity}
\end{center}
\end{figure}

\begin{Prop} \label{Prop:flat=spherical}
A flat $A_2$-planar algebra is spherical.
\end{Prop}
\emph{Proof:}
Given a 0-tangle, we isotope the strings so that we have a $\sigma\sigma^{\ast}$-tangle $T$, where $|\sigma| = n$ for $n \in \mathbb{N}$, with the $n$ vertices along the top and bottom of $T$ connected by closed strings which pass to the left of $T$. Then the string from the $n^{\textrm{th}}$ vertex along the top and bottom of $T$ can be pulled over all the other strings and all internal discs of $T$, introducing two opposite kinks, which contribute a scalar factor $q^{8/3} q^{-8/3} = 1$ (see Figure \ref{fig:flatness_gives_sphericity}). We may similarly pull the other strings which pass to the left of $T$ over $T$.
\hfill
$\Box$

The $A_2$-planar algebra $PTL$ is clearly flat, since the labelling set $L_{\pm} = \varnothing$. Then by Proposition \ref{Prop:flat=spherical} we see that there is only one trace on the algebra $\mathcal{V}^{A_2}_m$ in Section \ref{sec:A_2-tangles}.

\subsection{The involution on $P$} \label{sect:involution_on_P}

We can define the adjoint $T^{\ast} \in \mathcal{P}_{\sigma^{\ast}}(L)$ of a tangle $T \in \mathcal{P}_{\sigma}(L)$ , where $L$ has a $\ast$ operation defined on it, by reflecting the whole tangle about the horizontal line that passes through its centre and reversing all orientations. The labels $x_k \in L$ of $T$ are replaced by labels $x_k^{\ast}$ in $T^{\ast}$. If $\varphi$ is the map which sends $T \rightarrow T^{\ast}$, then every region $\varphi (R)$ of $T^{\ast}$ has the same colour as the region $R$ of $T$. For any linear combination of tangles in $\mathcal{P}_{\sigma}(L)$ we extend $\ast$ by conjugate linearity.
Then $P$ is an $A_2$-planar $\ast$-algebra if each $P_{\sigma}$ is a $\ast$-algebra, and for a $\sigma$-tangle $T$ with internal discs $D_k$ with patterns $\sigma_k$, labelled by $x_k \in P_{\sigma_k}$, we have
$$Z(T)^{\ast} = Z(T^{\ast}),$$
where the labels of the discs in $T^{\ast}$ are $x_k^{\ast}$, and where the definition of $Z(T)^{\ast}$ is extended to linear combinations of $\sigma$-tangles by conjugate linearity.
For $x_j \in P_{\sigma_j}$, $j=1,2$, we define the tangle $m(x_1,x_2) \in \mathcal{P}_{\sigma_1\sigma_2}$ by:
\begin{center}
  \includegraphics[width=25mm]{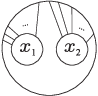}.
\end{center}

\begin{figure}[tb]
\begin{center}
  \includegraphics[width=50mm]{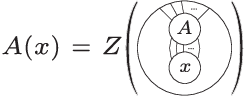}\\
 \caption{Action of $A \in P_{\sigma\sigma^{\ast}}$ on $x \in P_{\sigma}$} \label{fig:C^star-module_action}
\end{center}
\end{figure}

\begin{Prop}
Let $P$ be an $A_2$-planar $\ast$-algebra. Then $\mathrm{dim}(P_{\sigma}) \leq \mathrm{dim}(P_{\sigma\sigma^{\ast}})$ for any sign string $\sigma$.
\end{Prop}
\emph{Proof:}
Fix an element $y \in PTL_{\sigma^{\ast}}$ (i.e. the tangle $y$ does not contain any internal discs) such that $0 \neq c_y = {}_R \mathrm{Tr}(m(y,y^{\ast})) \in \mathbb{C}$. We have an embedding $\iota_y: P_{\sigma} \hookrightarrow P_{\sigma\sigma^{\ast}}$ given by $\iota_y(\; \cdot \;) = Z(m(\; \cdot \;,y))$. Let $\iota'_y: P_{\sigma\sigma^{\ast}} \rightarrow P_{\sigma}$ be the map defined by $\iota'_y (A) = c_y^{-1} A(y^{\ast})$, and the action of $P_{\sigma\sigma^{\ast}}$ on $P_{\sigma}$ is given in Figure \ref{fig:C^star-module_action}, for $A \in P_{\sigma\sigma^{\ast}}$, $x \in P_{\sigma}$.
Then $\iota'_y \circ \iota_y = \mathrm{id}$ on $P_{\sigma}$, and thus $\mathrm{dim}(P_{\sigma}) \leq \mathrm{dim}(P_{\sigma\sigma^{\ast}})$.

\begin{Cor} \label{Cor:P-f.d.}
An $A_2$-planar $\ast$-algebra $P$ is finite dimensional if and only if $\mathrm{dim}(P_{\sigma\sigma^{\ast}}) < \infty$ for any sign string $\sigma$.
\end{Cor}

The partition function $Z:P_{\sigma} \rightarrow \mathbb{C}$ on an $A_2$-planar algebra will be called \textbf{positive} if ${}_R \mathrm{Tr}_{\sigma^{\ast}\sigma} (m(x^{\ast},x)) \geq 0$, for all $x \in P_{\sigma}$, and \textbf{positive definite} if ${}_R \mathrm{Tr}_{\sigma^{\ast}\sigma} (m(x^{\ast},x)) > 0$, for all non-zero $x \in P_{\sigma}$.
The proof of \cite[Prop. 1.33]{jones:planar} in the $A_1$-case carries over to $A_2$-planar algebras where the only modification is that we allow possibly an odd number of vertices on discs, and different orientations on the strings.

\begin{Prop}\label{Prop:non-degenerate=tr_pos_def}
Let $P$ be an $A_2$-planar $\ast$-algebra with positive partition function $Z$. The following three conditions are equivalent:
(i) $P$ is non-degenerate,
(ii) ${}_R \mathrm{Tr}_{\sigma\sigma^{\ast}}$ is positive definite,
(iii) ${}_L \mathrm{Tr}_{\sigma\sigma^{\ast}}$ is positive definite.
\end{Prop}

Then we have the following Corollary, c.f. \cite[Cor. 1.36]{jones:planar}:
\begin{Cor}
If $P$ is a non-degenerate finite-dimensional $A_2$-planar $\ast$-algebra with positive partition function then $P_{\sigma\sigma^{\ast}}$ is semisimple for all sign strings $\sigma$, so there is a unique norm on $P_{\sigma\sigma^{\ast}}$ making it into a $C^{\ast}$-algebra. Each $P^{\sigma}$ is a Hilbert $C^{\ast}$-module over $P_{\sigma\sigma^{\ast}}$, for the action of $P_{\sigma\sigma^{\ast}}$ on $P_{\sigma}$ given above.
\end{Cor}

\begin{Def}
We call an $A_2$-planar algebra over $\mathbb{R}$ or $\mathbb{C}$ an \textbf{$A_2$-$C^{\ast}$-planar algebra} if it is a non-degenerate finite-dimensional $A_2$-planar $\ast$-algebra with positive definite partition function.
\end{Def}

If $P$ is a spherical $A_2$-$C^{\ast}$-planar algebra we can define an inner-product on $P_{\sigma}$, for $\sigma$ a sign string of length $n$, by $\langle x,y \rangle = \alpha^{-n/2} \mathrm{Tr}_{\sigma^{\ast}\sigma}(m(x^{\ast},y))$ for $x,y \in P_{\sigma}$. This inner product is normalized in the sense that $\langle \mathbf{1}_{\sigma\sigma^{\ast}}, \mathbf{1}_{\sigma\sigma^{\ast}} \rangle = 1$ for any sign string $\sigma$.

\section{$A_2$-planar $i,j$-tangles} \label{sect:i,j-tangles}

We will be particularly interested in the vector spaces $P_{\sigma}$ for sign strings $\sigma$ with a particular form, since these will correspond exactly to the vector spaces in the double complex associated to the $SU(3)$-subfactors. We describe these vector spaces in the next sections, and introduce certain basic tangles which will play an important role later.

An $A_2$-planar $i,j$-tangle will be an $A_2$-planar $\sigma$-tangle with external disc $D=D_0$ and internal discs $D_1, \ldots, D_n$, where each disc $D_k$, $k \geq 0$, has pattern $\sigma^{(k)} = -^{j_k} \cdot \widetilde{\sigma}^{(k)} \cdot +^{j_k}$, where $\widetilde{\sigma}^{(k)}$ is the alternating string of length $2i_k$ which begins with `$-$'.
We will position the vertices so that the first $i_k + j_k$ are along the boundary for the upper half of the disc, which we will call the top edge, and the next $i_k + j_k$ vertices are along the boundary for the bottom half of the disc, which we will call the bottom edge. We will use the convention of numbering the vertices along the bottom edge in reverse order, so that the $2(i_k+j_k)$-th vertex is called the first vertex along the bottom edge. The total number of source vertices along the top edge is $\lfloor j_k + (i_k+1)/2 \rfloor$, and the number of sink vertices is $\lfloor i_k/2 \rfloor$.
For the outer boundary $\partial D$ we impose the restriction $b_0 =0$.

It is important to note that what we here call an $i,j$-tangle is different from the $(i,j)$-tangles of Section \ref{sect:A2tangles}. In both cases the integers $i,j$ refer to the number of vertices along the top (and bottom) edge of the disc, however in an $(i,j)$-tangle the first $i$ vertices are all sources, and the next $j$ vertices are all sinks.

In the figures which follow we omit the orientation on the strings from the last $i$ vertices along the top and bottom of an $i,j$-tangle -- these will be alternating.

\subsection{Some basic $A_2$-planar $i,j$-tangles} \label{Sect:Basic_planar_tangles}

The following basic tangles will be of importance to us: \\
$\bullet \quad$ \emph{Inclusion tangles} $IR^{i,j}_{i+1,j}$, $IR^{i,j}_{i,j+1}$ and $\widetilde{IR}^{i,j}_{i,j+1}$:\\
\begin{center}
  \includegraphics[width=65mm]{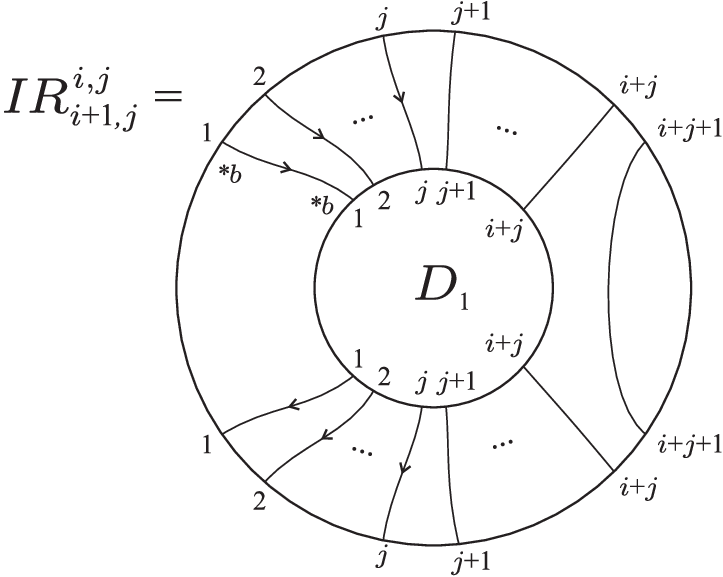} \\
  \mbox{} \\
  \includegraphics[width=140mm]{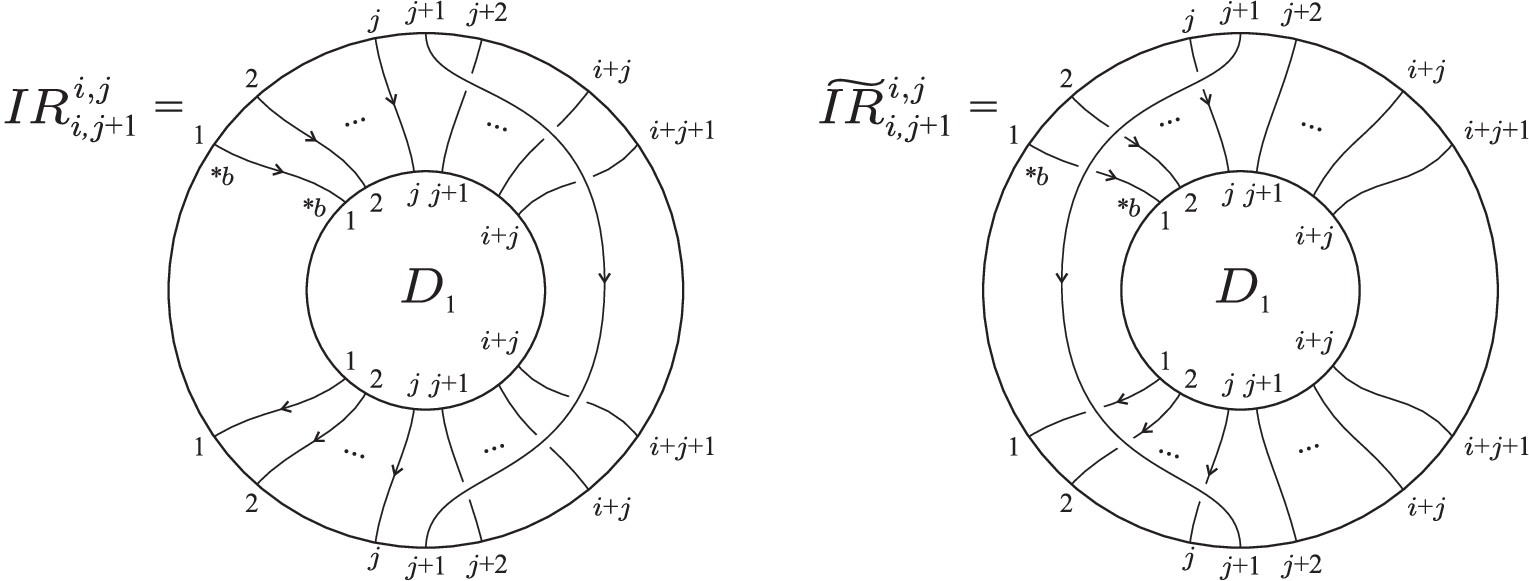}
\end{center}
where the orientation of the rightmost string in $IR^{i,j}_{i+1,j}$ is downwards for $i$ even and upwards for $i$ odd. Both $IR^{i,j}_{i,j+1}$ and $\widetilde{IR}^{i,j}_{i,j+1}$ add a new source vertex along the top immediately to the right of the first $j$ source vertices, and a sink vertex along the bottom immediately to the right of the first $j$ sink vertices along the bottom. These new vertices are regarded as being among the downwards oriented vertices rather than the alternating vertices. They are connected by a through string, and $IR^{i,j}_{i,j+1}$, $\widetilde{IR}^{i,j}_{i,j+1}$ differ only in that the through string passes to the right of the inner disc in $IR^{i,j}_{i,j+1}$ and to the left in $\widetilde{IR}^{i,j}_{i,j+1}$.
We have $Z(IR^{i,j}_{i+1,j}): P_{i,j} \rightarrow P_{i+1,j}$, and $Z(IR^{i,j}_{i,j+1}), Z(\widetilde{IR}^{i,j}_{i,j+1}), : P_{i,j} \rightarrow P_{i,j+1}$.

For a flat $A_2$-planar algebra, the two right inclusion tangles $IR^{i,j}_{i,j+1}$ and $\widetilde{IR}^{i,j}_{i,j+1}$ are equal, and we will simply write $IR^{i,j}_{i,j+1}$.
For a spherical $A_2$-planar algebra $P$, we define $\mathrm{tr}(x) = \alpha^{-i-j} \, \mathrm{Tr}_{i,j}(x)$ for $x \in P_{i,j}$. Then $\mathrm{tr}$ is compatible with the inclusions $P_{i,j} \subset P_{i,j+1}$ and $P_{i,j} \subset P_{i+1,j}$, given by $IR^{i,j}_{i,j+1}$, $IR^{i,j}_{i+1,j}$ respectively, and $\mathrm{tr}(1) = 1$, and so defines a trace on $P$ itself.
If $P$ is a spherical $A_2$-$C^{\ast}$-planar the inner-product defined at the end of Section \ref{sect:involution_on_P} is given on $P_{i,j}$ by $\langle x,y \rangle = \mathrm{tr}(x^{\ast}y)$ for $x,y \in P_{i,j}$, and is consistent with the inclusions $P_{i,j} \subset P_{i,j+1}$ and $P_{i,j} \subset P_{i+1,j}$ given above, since $\mathrm{tr}$ is.
\\ \\
$\bullet \quad$ \emph{Conditional expectation tangles} $ER^{i+1,j}_{i,j}$ and $ER^{i,j+1}_{i,j}$:
\vspace{-1cm}\\
\begin{center}
 \begin{minipage}[b]{16cm}
  \begin{minipage}[t]{14cm}
  \mbox{} \\
   \includegraphics[width=140mm]{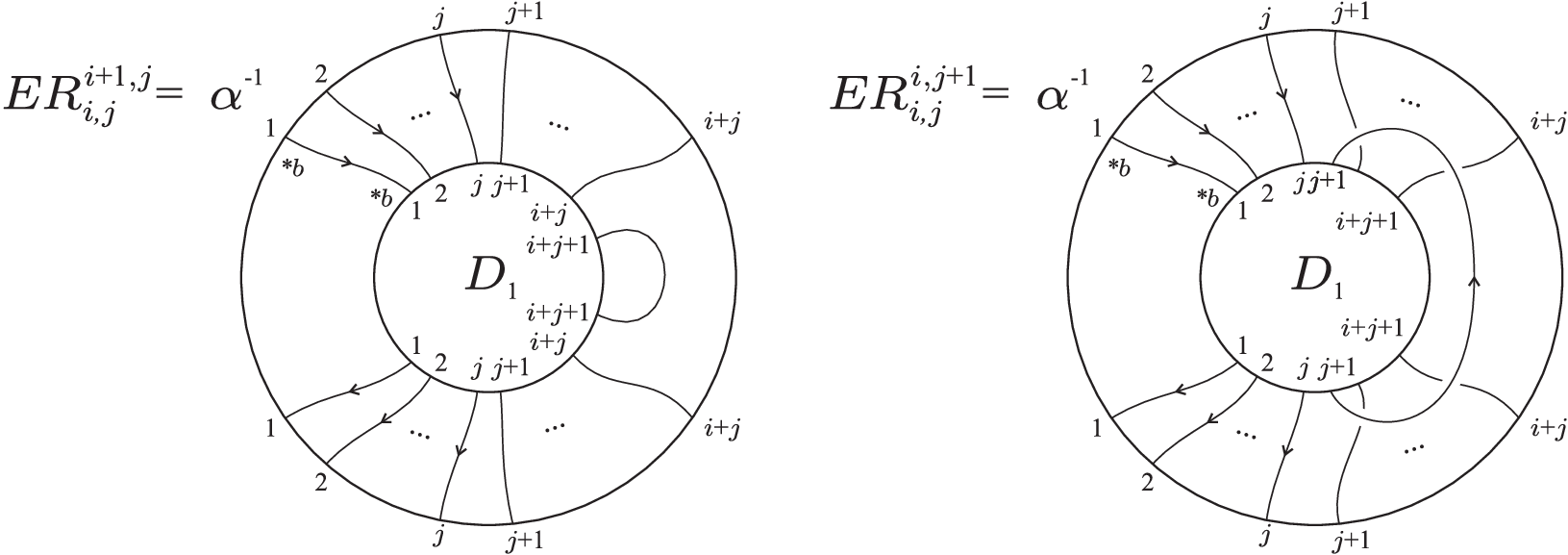}
  \end{minipage}
  \begin{minipage}[t]{1.5cm}
   \mbox{} \\
   \hfill
   \parbox[t]{1.5cm}{ \begin{eqnarray}\label{fig:conditional_expectation_tangles}\end{eqnarray}}
  \end{minipage}
 \end{minipage}
\end{center}
The orientation of the string from vertex $i+j+1$ on the inner disc of $ER^{i+1,j}_{i,j}$ is clockwise for $i$ odd and anticlockwise for $i$ even. We have $Z(ER^{i+1,j}_{i,j}): P_{i+1,j} \rightarrow P_{i,j}$ and $Z(ER^{i,j+1}_{i,j}): P_{i,j+1} \rightarrow P_{i,j}$.

Let $\mathcal{P}^{(1)}_{i,j}$ denote the subset of $\mathcal{P}_{i,j}$ spanned by all tangles where vertices $j+1$ along the top and bottom are connected by a through string which passes over every string it crosses and such that there are no internal discs in the region between this string and the outer boundary of the tangle to the left of it. If $P$ is a general $A_2$-planar algebra with presenting map $Z$, we define $P^{(1)}_{i,j} = Z(\mathcal{P}^{(1)}_{i,j}) \subset P_{i,j}$, and denote by $P^{(1)} \subset P$ the subspace $P^{(1)} = \bigcup_{i,j} P^{(1)}_{i,j}$.
We also have left conditional expectation tangles $EL^{i+1,j}_{i+1,j}$ and $EL^{i,j+1}_{i,j+1}$:
\begin{center}
  \includegraphics[width=150mm]{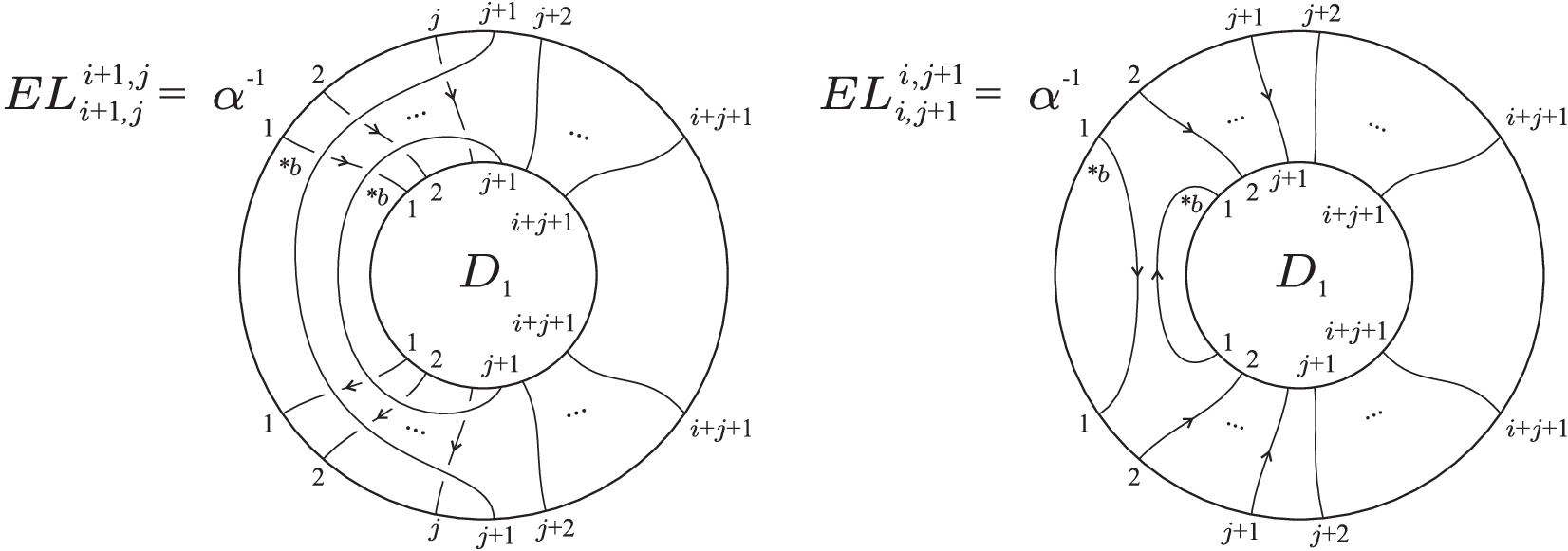}
\end{center}
where $Z(EL^{i+1,j}_{i+1,j}): P_{i+1,j} \rightarrow P^{(1)}_{i+1,j}$. \\

\begin{figure}[b]
\begin{center}
  \includegraphics[width=48mm]{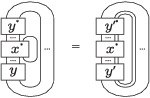}\\
 \caption{} \label{fig:fig34}
\end{center}
\end{figure}

The justification for calling the tangles in (\ref{fig:conditional_expectation_tangles}) conditional expectation tangles is seen in the following Lemma:

\begin{Lemma}
Let $P$ be an $A_2$-$C^{\ast}$-planar algebra. For the tangles $ER^{i+1,j}_{i,j}$ and $ER^{i,j+1}_{i,j}$ defined in (\ref{fig:conditional_expectation_tangles}), $E_1(x) = Z(ER^{i+1,j}_{i,j}(x))$ is the conditional expectation of $x \in P_{i+1,j}$ onto $P_{i,j}$ with respect to the trace, and $E_2(y) = Z(ER^{i,j+1}_{i,j}(y))$ is the conditional expectation of $y \in P_{i,j+1}$ onto $P_{i,j}$ with respect to the trace.
\end{Lemma}
\emph{Proof:}
We first check positivity of $E_1(x)$ for positive $x \in P_{i+1,j}$. As $P$ is an $A_2$-$C^{\ast}$-planar algebra, the inner-product defined above is positive definite. We need to show that $\langle E_1(x) y, y \rangle \geq 0$ for all $y \in P_{i,j}$. From Figure \ref{fig:fig34} we see that $\mathrm{tr}(y^{\ast} ER^{i+1,j}_{i,j}(x)^{\ast} y) = \mathrm{tr}(y'^{\ast} x^{\ast} y') = \langle xy', y' \rangle \geq 0$ for all $y \in P_{i,j}$, where $y' =$ \includegraphics[width=7mm]{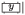} $\in P_{i+1,j}$.
From
\begin{center}
  \includegraphics[width=55mm]{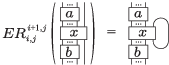}
\end{center}
we see that $E_1(axb) = aE_1(x)b$, for $x \in P_{i+1,j}$, $a,b \in P_{i,j}$.
Since also $\langle E_1(x), y \rangle = \langle x, y' \rangle$, $E_1$ is the trace-preserving conditional expectation from $P_{i+1,j}$ onto $P_{i,j}$. The proof for $E_2$ is similar.
\hfill
$\Box$

Similarly, $Z(EL^{i+1,j}_{i+1,j}(x))$ is the conditional expectation of $x \in P_{i+1,j}$ onto $P^{(1)}_{i+1,j}$.

\subsection{Dimensions in $A_2$-planar algebras and $PTL$.} \label{sect:dimP}

\begin{figure}[tb]
\begin{center}
  \includegraphics[width=130mm]{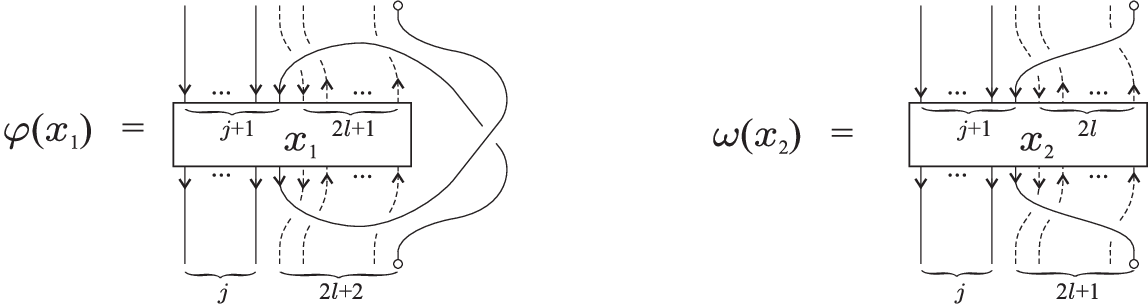}
  \caption{Maps $\varphi: \mathcal{P}_{2l+1,j+1}(L) \rightarrow \mathcal{P}_{2l+2,j}(L)$, $\omega: \mathcal{P}_{2l,j+1}(L) \rightarrow \mathcal{P}_{2l+1,j}(L)$} \label{fig:varphi-omega}
\end{center}
\end{figure}

\begin{figure}[tb]
\begin{center}
  \includegraphics[width=135mm]{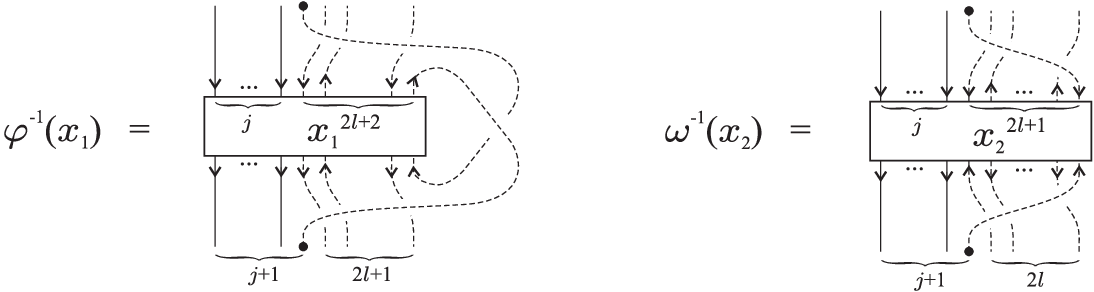}
  \caption{Maps $\varphi^{-1}: \mathcal{P}_{2l+2,j}(L) \rightarrow \mathcal{P}_{2l+1,j+1}(L)$, $\omega^{-1}: \mathcal{P}_{2l+1,j}(L) \rightarrow \mathcal{P}_{2l,j+1}(L)$} \label{fig:inv_varphi-omega}
\end{center}
\end{figure}

We now present some results regarding the dimensions of the different graded parts of $A_2$-planar algebras. These will be needed later in Section \ref{sect:Planar_algebras_give_subfactors}.
We define maps $\varphi: \mathcal{P}_{2l+1,j+1}(L) \rightarrow \mathcal{P}_{2l+2,j}(L)$, $\omega: \mathcal{P}_{2l,j+1}(L) \rightarrow \mathcal{P}_{2l+1,j}(L)$ as in Figure \ref{fig:varphi-omega} for $x_1 \in \mathcal{P}_{2l+1,j+1}(L)$, $x_2 \in \mathcal{P}_{2l,j+1}(L)$, where the white circle at the end of a string indicates that this vertex is now regarded as one of the $i$ vertices of $\mathcal{P}_{i,j}$ with alternating orientation ($i=2l+2,2l+1$ for $\varphi$, $\omega$ respectively). The maps $\varphi$, $\omega$ are invertible, with $\varphi^{-1}$, $\omega^{-1}$ as in Figure \ref{fig:inv_varphi-omega} for $x_1 \in \mathcal{P}_{2l+2,j}(L)$, $x_2 \in \mathcal{P}_{2l+1,j}(L)$, where the solid black circle at the end of a string indicates that this vertex is now regarded as one of the $j+1$ vertices of $\mathcal{P}_{i,j+1}$ with alternating orientation ($i=2l+1,2l$ for $\varphi^{-1}$, $\omega^{-1}$ respectively).
Clearly $\varphi(\mathcal{P}_{2l+1,j+1}(L)) \subset \mathcal{P}_{2l+2,j}(L)$. Since $\mathcal{P}_{2l+1,j+1}(L) \supset \varphi^{-1}(\mathcal{P}_{2l+2,j}(L))$ then $\varphi(\mathcal{P}_{2l+1,j+1}(L)) \supset \mathcal{P}_{2l+2,j}(L)$. So $\varphi(\mathcal{P}_{2l+1,j+1}(L)) = \mathcal{P}_{2l+2,j}(L)$ and $\varphi$ is a bijection. Similarly the map $\omega$ is a bijection and $\omega(\mathcal{P}_{2l,j+1}(L)) = \mathcal{P}_{2l+1,j}(L)$. Let $Z:\mathcal{P}_{i,j}(L) \rightarrow P_{i,j}$ be the presenting map for an $A_2$-$C^{\ast}$-planar algebra $P$. We define bijections $\widetilde{\varphi}: P_{2l+1,j+1}(L) \rightarrow P_{2l+2,j}(L)$, $\widetilde{\omega}: P_{2l,j+1}(L) \rightarrow P_{2l+1,j}(L)$ by $\widetilde{\varphi}(x_1) = Z(\varphi(x_1))$ and $\widetilde{\omega}(x_2) = Z(\omega(x_2))$.
Then
\begin{equation} \label{eqn:dimP=dimP-1}
\mathrm{dim}(\mathcal{P}_{i,j}(L)) = \mathrm{dim}(\mathcal{P}_{i+k,j-k}(L)),
\end{equation}
\begin{equation} \label{eqn:dimP=dimP-2}
\mathrm{dim}(P_{i,j}) = \mathrm{dim}(P_{i+k,j-k}),
\end{equation}
for all integers $k$ such that $-i \leq k \leq j$.
Note, (\ref{eqn:dimP=dimP-1}) follows immediately from \cite[Theorem 6.3]{kuperberg:1996}.

For $L=\varnothing$, we define $\mathcal{PTL}_{i,j}$ to be the quotient of $PTL_{i,j} = \mathcal{P}_{i,j}(\varnothing)$ by the subspace of zero-length vectors with respect to the inner-product on $PTL_{i,j}$ defined by $\langle x,y \rangle = \widehat{x^{\ast}y}$, for $x,y \in PTL_{i,j}$, where $\widehat{T}$ is the tangle defined as in Figure \ref{fig:trace-tangle}.
The element $\varphi(x)$ is a zero-length vector in $PTL_{2l+2,j}$ if and only if $x$ is a zero-length vector in $PTL_{2l+1,j+1}$. Similarly, $\omega(x)$ is zero-length vector in $PTL_{2l+1,j}$ if and only if $x$ is a zero-length vector in $PTL_{2l,j+1}$.
Thus for all integers $k$ with $-i \leq k \leq j$, $\mathrm{dim}(\mathcal{PTL}_{i,j}) = \mathrm{dim}(\mathcal{PTL}_{i+k,j-k})$.

\section{$A_2$-Planar algebra description of subfactors} \label{sect:Planar_algebras_give_subfactors}

We are now going to associate flat $A_2$-planar $C^{\ast}$-algebras to the double sequences of subfactors associated to $\mathcal{ADE}$ graphs with flat connections. These double sequences, introduced in \cite{evans/kawahigashi:1994}, have a periodicity three coming from the $A_2$-Temperley-Lieb algebra in the horizontal direction, and a periodicity two coming from the subfactor basic construction in the vertical direction.
In Section \ref{sect:planar_alg_for_A=PTL} we give a diagrammatic form for the double sequences for the Wenzl subfactors.

\begin{figure}[tb]
\begin{center}
\includegraphics[width=100mm]{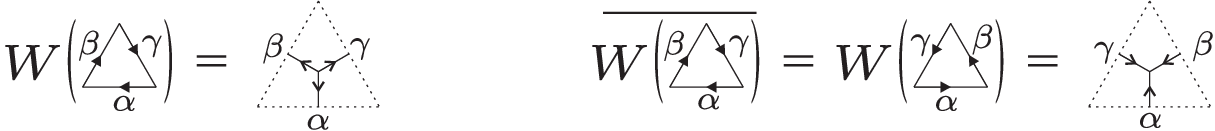}\\
 \caption{Cells associated to trivalent vertices} \label{fig:Oc-Kup}
\end{center}
\end{figure}

Let $\mathcal{G}$ be any finite $SU(3)$ $\mathcal{ADE}$ graph with Coxeter number $n$. Let $\alpha = [3]_q$, $q=e^{i \pi / n}$, be the Perron-Frobenius eigenvalue of $\mathcal{G}$ and let $(\phi_v)$ be the corresponding eigenvector.
Ocneanu \cite{ocneanu:2000ii} defined a cell system $W$ on $\mathcal{G}$ by associating a complex number $W \left( \triangle^{(\alpha \beta \gamma)} \right)$, called an Ocneanu cell, to each closed loop of length three $\triangle^{(\alpha \beta \gamma)}$ in $\mathcal{G}$ as in Figure \ref{fig:Oc-Kup}, where $\alpha$, $\beta$, $\gamma$ are edges on $\mathcal{G}$. These cells satisfy two properties, called Ocneanu's type I, II equations respectively, which are obtained by evaluating the Kuperberg relations K2, K3 respectively, using the identification in Figure \ref{fig:Oc-Kup}: \\
$(i)$ for any type I frame \includegraphics[width=16mm]{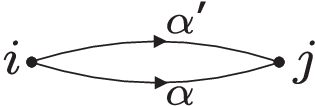} in $\mathcal{G}$ we have \\
\begin{flushright}
\begin{minipage}[b]{14cm}
 \begin{minipage}[t]{9cm}
  \mbox{} \\
 \includegraphics[width=80mm]{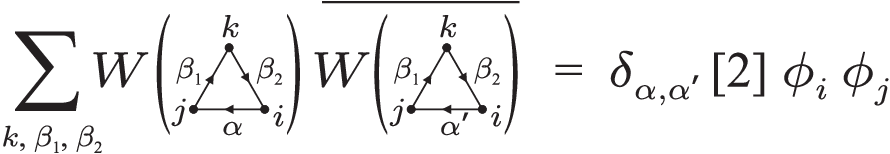}
 \end{minipage}
 \hfill
 \begin{minipage}[t]{1.5cm}
  \mbox{} \\
  \hfill
  \parbox[t]{7mm}{\begin{eqnarray}\label{eqn:typeI_frame}\end{eqnarray}}
 \end{minipage}
\end{minipage}
\end{flushright}
$(ii)$ for any type II frame \includegraphics[width=30mm]{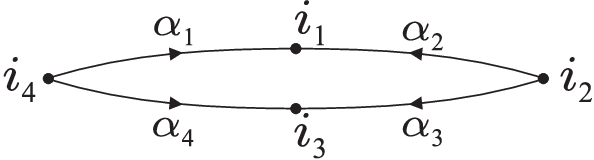} in $\mathcal{G}$ we have \\
\begin{flushright}
\begin{minipage}[b]{14cm}
 \begin{minipage}[b]{9cm}
  \mbox{} \\
 \includegraphics[width=90mm]{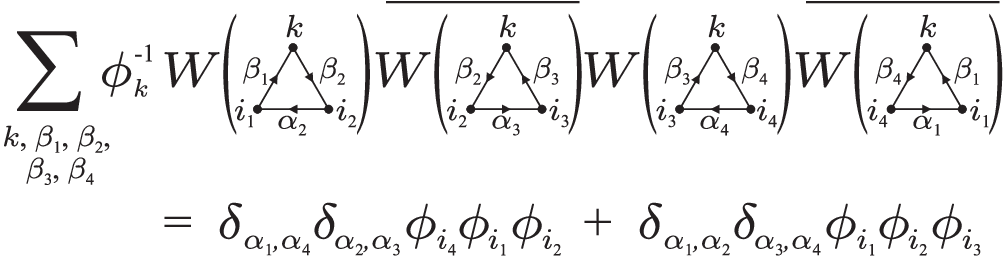}
 \end{minipage}
 \hfill
 \begin{minipage}[b]{1.5cm}
  \mbox{} \\
  \hfill
  \parbox[b]{7mm}{\begin{eqnarray}\label{eqn:typeII_frame}\end{eqnarray}}
 \end{minipage}
\end{minipage}
\end{flushright}
The existence of these cells for the finite $\mathcal{ADE}$ graphs was shown in \cite{evans/pugh:2009i} with the exception of the graph $\mathcal{E}_4^{(12)}$. Using these cells, we define a representation $\mathcal{U}^{\rho_1,\rho_2}_{\rho_3,\rho_4}$ of the Hecke algebra by
\begin{equation} \label{eqn:HeckeRep}
\mathcal{U}^{\rho_1,\rho_2}_{\rho_3,\rho_4} = \sum_{\lambda} \phi_{s(\rho_1)}^{-1} \phi_{r(\rho_2)}^{-1} W(\triangle^{(\lambda, \rho_3, \rho_4)}) \overline{W(\triangle^{(\lambda, \rho_1, \rho_2)})},
\end{equation}
for edges $\rho_1$, $\rho_2$, $\rho_3$, $\rho_4$, $\lambda$ of $\mathcal{G}$.

As in \cite{evans/kawahigashi:1994}, with any choice of distinguished vertex $\ast$, we define the double sequence $(B_{i,j})$ of finite dimensional algebras by:
$$
\begin{array}{cccccccc}
B_{0,0} & \subset & B_{0,1} & \subset & B_{0,2} & \subset & \cdots \qquad \longrightarrow & B_{0,\infty} \\
\cap & & \cap & & \cap & & & \cap \\
B_{1,0} & \subset & B_{1,1} & \subset & B_{1,2} & \subset & \cdots \qquad \longrightarrow & B_{1,\infty} \\
\cap & & \cap & & \cap & & & \cap \\
B_{2,0} & \subset & B_{2,1} & \subset & B_{2,2} & \subset & \cdots \qquad \longrightarrow & B_{2,\infty} \\
\cap & & \cap & & \cap & & & \cap \\
\vdots & & \vdots & & \vdots & & & \vdots
\end{array}
$$
The Bratteli diagrams for horizontal inclusions $B_{i,j} \subset B_{i,j+1}$ are given by $\mathcal{G}$. If $\mathcal{G}$ is three-colourable, the vertical inclusions $B_{i,j} \subset B_{i+1,j}$ are given by its $\overline{j}, \overline{j+1}$-part $\mathcal{G}_{\overline{j}, \overline{j+1}}$, where $\overline{p} = \tau(p)$ is the colour of $p$ for $p = j,j+1$. We identify $B_{0,0} = \mathbb{C}$ with the distinguished vertex $\ast$ of $\mathcal{G}$.

Then for the inclusions
\begin{equation} \label{commuting_square_of_Bij's}
\begin{array}{ccc}
B_{i,j} & \subset & B_{i,j+1} \\
\cap & & \cap \\
B_{i+1,j} & \subset & B_{i+1,j+1}
\end{array}
\end{equation}
with $i$ even, we define a connection by
\begin{equation} \label{eqn:connection}
X^{\rho_1,\rho_2}_{\rho_3,\rho_4} = \begin{array}{c}
\stackrel{\rho_1}{\longrightarrow}\\
\scriptstyle \rho_3 \textstyle \big\downarrow \qquad \big\downarrow \scriptstyle \rho_2 \\
\stackrel{\textstyle\longrightarrow}{\scriptstyle{\rho_4}}
\end{array}
= q^{2/3} \delta_{\rho_1, \rho_3} \delta_{\rho_2, \rho_4} - q^{-1/3} \, \mathcal{U}^{\rho_1,\rho_2}_{\rho_3,\rho_4} ,
\end{equation}
We denote by $\widetilde{\mathcal{G}}$ the reverse graph of $\mathcal{G}$, which is the graph obtained by reversing the direction of every edge of $\mathcal{G}$. For the inclusions (\ref{commuting_square_of_Bij's}) with $i$ odd, let $\rho_1$, $\rho_4$ be edges on $\mathcal{G}$ and let $\widetilde{\rho}_2$, $\widetilde{\rho}_3$ be edges on the reverse graph $\widetilde{\mathcal{G}}$ (so that $\rho_2$, $\rho_3$ are edges on $\mathcal{G}$). We define the connection by
\begin{equation} \label{eqn:inverse_graph_connection}
X^{\rho_1,\widetilde{\rho}_2}_{\widetilde{\rho}_3,\rho_4} = \begin{array}{c}
\stackrel{\rho_1}{\longrightarrow}\\
\scriptstyle \widetilde{\rho}_3 \textstyle \big\downarrow \qquad \big\downarrow \scriptstyle \widetilde{\rho}_2 \\
\stackrel{\textstyle\longrightarrow}{\scriptstyle{\rho_4}}
\end{array}
= \sqrt{\frac{\phi_{s(\rho_3)} \phi_{r(\rho_2)}}{\phi_{r(\rho_3)} \phi_{s(\rho_2)}}}
\overline{\begin{array}{c}
\stackrel{\rho_4}{\longrightarrow}\\
\scriptstyle \rho_3 \textstyle \big\downarrow \qquad \big\downarrow \scriptstyle \rho_2 \\
\stackrel{\textstyle\longrightarrow}{\scriptstyle{\rho_1}}
\end{array}}.
\end{equation}
It was shown in \cite{evans/pugh:2009i} that these connections satisfy the unitarity axiom
\begin{equation} \label{eqn:unitarity_property_of_connections}
\sum_{\rho_3,\rho_4} X^{\rho_1,\rho_2}_{\rho_3,\rho_4} \; \overline{X^{\rho_1',\rho_2'}_{\rho_3,\rho_4}} = \delta_{\rho_1, \rho_1'} \delta_{\rho_2, \rho_2'}.
\end{equation}

Then for the inclusions (\ref{commuting_square_of_Bij's}) an element indexed by paths in the basis
\setlength{\unitlength}{1mm}
\begin{picture}(6,3.5)
\put(1,3.5){\vector(0,-1){3.5}}
\put(1,0){\vector(1,0){4}}
\end{picture}
can be transformed to an element indexed by paths in the basis
\begin{picture}(5,3.5)
\put(0,3.5){\vector(1,0){4}}
\put(4,3.5){\vector(0,-1){3.5}}
\end{picture}
using the above connections: Let $(\sigma \cdot \sigma' \cdot \alpha_1 \cdot \alpha_2, \sigma \cdot \sigma' \cdot \alpha_1' \cdot \alpha_2')$ be an element in $B_{i+1,j+1}$ in the basis
\begin{picture}(6,3.5)
\put(1,3.5){\vector(0,-1){3.5}}
\put(1,0){\vector(1,0){4}}
\end{picture}
, where $\sigma$ is a horizontal path of length $j$, $\sigma'$ is a vertical path of length $i$, $\alpha_1$, $\alpha_1'$ are vertical paths of length 1, $\alpha_2$, $\alpha_2'$ are horizontal paths of length 1, and $r(\alpha_2) = r(\alpha_2')$. We transform this to an element in the basis
\begin{picture}(5,3.5)
\put(0,3.5){\vector(1,0){4}}
\put(4,3.5){\vector(0,-1){3.5}}
\end{picture}
by
$$(\sigma \cdot \sigma' \cdot \alpha_1 \cdot \alpha_2, \sigma \cdot \sigma' \cdot \alpha_1' \cdot \alpha_2') = \sum_{\beta_i, \beta_i'} \begin{array}{c}
\stackrel{\beta_1}{\longrightarrow}\\
\scriptstyle \alpha_1 \textstyle \big\downarrow \qquad \big\downarrow \scriptstyle \beta_2 \\
\stackrel{\textstyle\longrightarrow}{\scriptstyle{\alpha_2}}
\end{array}
\overline{\begin{array}{c}
\stackrel{\beta_1'}{\longrightarrow}\\
\scriptstyle \alpha_1' \textstyle \big\downarrow \qquad \big\downarrow \scriptstyle \beta_2' \\
\stackrel{\textstyle\longrightarrow}{\scriptstyle{\alpha_2'}}
\end{array}} (\sigma \cdot \sigma' \cdot \beta_1 \cdot \beta_2, \sigma \cdot \sigma' \cdot \beta_1' \cdot \beta_2'),$$
where the summation is over all horizontal paths $\beta_1$, $\beta_1'$ of length 1, and vertical paths $\beta_2$, $\beta_2'$ of length 1.

The Markov trace on $B_{i,j}$ is defined as in \cite{evans/kawahigashi:1994} by
\begin{equation} \label{eqn:trace}
\mathrm{tr}((\sigma_1,\sigma_2)) = \delta_{\sigma_1, \sigma_2} [3]^{-k} \phi_{r(\sigma_1)},
\end{equation}
for $(\sigma_1,\sigma_2) \in B_{i,j}$, where $k=i+j$.
We define $B_{i,\infty}$ to be the GNS-completion of $\bigcup_{k\geq0} B_{i,k}$ with respect to the trace. As in \cite{evans/kawahigashi:1994}, the braid elements
\begin{center}
  \includegraphics[width=60mm]{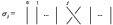}
\end{center}
appear as the connection.

If $\mathcal{G}$ is three-colourable then its adjacency matrix $\Delta_{\mathcal{G}}$ which may be written in the form
$$\Delta_{\mathcal{G}} = \left( { \begin{array}{ccc}
                 0 & \Delta_{01} & 0 \\
                 0 & 0 & \Delta_{12} \\
                 \Delta_{20} & 0 & 0
               \end{array} } \right),$$
where $\Delta_{01}$, $\Delta_{12}$ and $\Delta_{20}$ are matrices which give the number of edges between each 0,1,2-coloured vertex respectively of $\mathcal{G}$ to each 1,2,0-coloured vertex respectively. By a suitable ordering of the vertices the matrix $\Delta_{12}$ may be chosen to be symmetric. These matrices satisfy the conditions
\begin{equation} \label{eqn:DeltaDelta^T=Delta^TDelta}
\Delta_{01}^T \Delta_{01} = \Delta_{20} \Delta_{20}^T = \Delta_{12}^2,
\qquad \qquad
\Delta_{01} \Delta_{01}^T = \Delta_{20}^T \Delta_{20},
\end{equation}
which follow from the fact that $\Delta_{\mathcal{G}}$ is normal \cite{evans/pugh:2009ii}.

\begin{Lemma} \label{Lemma:dimB=dimB}
For the double sequence $(B_{i,j})$ defined above, $\textrm{dim}(B_{i,j}) = \textrm{dim}(B_{i+k,j-k})$ for all integers $k$ such that $-i \leq k \leq j$.
\end{Lemma}
\emph{Proof:}
If $\mathcal{G}$ is not three-colourable, then $B_{i,j}$ is the space of all pairs of paths of length $i+j$ on $\mathcal{G}$, hence the result is trivial. For the three-colourable graphs, let $\Lambda_{i,j}^1$ be the product of $j$ matrices $\Lambda_{i,j}^1 = \Delta_{01} \Delta_{12} \Delta_{20} \Delta_{01} \cdots \Delta_{\overline{j-1},\overline{j}}$, and $\Lambda_{i,j}^2$ the product of $i$ matrices $\Lambda_{i,j}^2 = \Delta_{\overline{j},\overline{j+1}} \Delta_{\overline{j},\overline{j+1}}^T \Delta_{\overline{j},\overline{j+1}} \Delta_{\overline{j},\overline{j+1}}^T \cdots \Delta'$, where $\Delta'$ is $\Delta_{\overline{j},\overline{j+1}}$ if $i$ is odd, $\Delta_{\overline{j},\overline{j+1}}^T$ if $i$ is even, and $\overline{p}$ is the colour of $p$. Then if $\Lambda_{i,j} = \Lambda_{i,j}^1 \Lambda_{i,j}^2$, the dimension of $B_{i,j}$ is given by $\left( \Lambda_{i,j} \Lambda_{i,j}^T \right)_{0,0}$. Using (\ref{eqn:DeltaDelta^T=Delta^TDelta}) it is easy to show by induction that $\Lambda_{i,j} \Lambda_{i,j}^T = (\Delta_{01} \Delta_{01}^T)^{i+j}$. So $\textrm{dim}(B_{i+k,j-k}) = \left( \Lambda_{i+k,j-k} \Lambda_{i+k,j-k}^T \right)_{0,0} = \left( (\Delta_{01} \Delta_{01}^T)^{i+j} \right)_{0,0} = \textrm{dim}(B_{i,j})$.
\hfill
$\Box$

For all $i,j \geq 0$ we define operators $U_{-k} \in B_{i,j}$, $k=0,1,\ldots, j-1$, which satisfy the Hecke relations H1-H3, by
\begin{eqnarray*}
U_{-k} & = & \sum_{\stackrel{|\zeta_1|=j-2-k, |\zeta'|= i}{\scriptscriptstyle{|\gamma_i|=|\eta_i|=1, |\zeta_2|=k}}} \mathcal{U}^{\gamma_2,\eta_2}_{\gamma_1,\eta_1}  \; (\zeta_1 \cdot \gamma_1 \cdot \eta_1 \cdot \zeta_2 \cdot \zeta', \zeta_1 \cdot \gamma_2 \cdot \eta_2 \cdot \zeta_2 \cdot \zeta'), \qquad 0 \leq k \leq j-2, \\
U_{-j+1} & = & \sum_{\stackrel{|\zeta|=j-1, |\zeta'| = i-1}{\scriptscriptstyle{|\gamma_i|=|\eta_i'|=1}}} \mathcal{U}^{\gamma_2,\eta_2}_{\gamma_1,\eta_1}  \; (\zeta \cdot \gamma_1 \cdot \eta_1' \cdot \zeta', \zeta \cdot \gamma_2 \cdot \eta_2' \cdot \zeta'),
\end{eqnarray*}
where $\xi, \xi'$ are horizontal, vertical paths respectively, and $\mathcal{U}^{\gamma_2,\eta_2}_{\gamma_1,\eta_1} $ are the Boltzmann weights for $\mathcal{A}^{(n)}$. The embedding of $U_{-k} \in B_{i,j}$ into $B_{i+1,j}$ is $U_{-k}$, whilst the embedding of $U_{-k} \in B_{i,j}$ into $B_{i,j+1}$ is $U_{-k-1}$.
We have $B_{i,j} \supset \textrm{alg} (U_{-j+1}, U_{-j+2}, \ldots, U_{-1}, U_0)$. When $\mathcal{G} = \mathcal{A}^{(n)}$, the algebra $B_{l,j} = \textrm{alg} (U_{-j+1}, U_{-j+2}, \ldots, U_{-l})$ for $l=0,1$ \cite{evans/kawahigashi:1994}.

\begin{Lemma} \label{Lemma:commuting_square_of_Bij's}
The square (\ref{commuting_square_of_Bij's}) is a commuting square.
\end{Lemma}
\emph{Proof:}
Note that for the $\mathcal{A}$ graphs, the result follows by \cite[Prop. 3.2]{wenzl:1988}. However, we prove the case for a general $SU(3)$ $\mathcal{ADE}$ graph $\mathcal{G}$.
By \cite[Theorem 11.2]{evans/kawahigashi:1998}, the square (\ref{commuting_square_of_Bij's}) is a commuting square if and only if the corresponding connection satisfies
\begin{equation} \label{eqn:condition_for_commuting_square}
\sum_{\sigma_2, \sigma_4} \frac{\phi_{r(\sigma_2)} \sqrt{\phi_{s(\sigma_3)} \phi_{s(\sigma_3')}}}{\phi_{s(\sigma_2)} \phi_{s(\sigma_4)}}
\begin{array}{c}
\stackrel{\sigma_1}{\longrightarrow}\\
\scriptstyle \sigma_3 \textstyle \big\downarrow \qquad \big\downarrow \scriptstyle \sigma_2 \\
\stackrel{\textstyle\longrightarrow}{\scriptstyle{\sigma_4}}
\end{array}
\overline{\begin{array}{c}
\stackrel{\sigma_1'}{\longrightarrow}\\
\scriptstyle \sigma_3' \textstyle \big\downarrow \qquad \big\downarrow \scriptstyle \sigma_2 \\
\stackrel{\textstyle\longrightarrow}{\scriptstyle{\sigma_4}}
\end{array}}
= \delta_{\sigma_1, \sigma_1'} \delta_{\sigma_3, \sigma_3'},
\end{equation}
where $\sigma_1$, $\sigma_1'$ are any edges on the graph of the Bratteli diagram for $B_{i,j} \subset B_{i,j+1}$, $\sigma_3$, $\sigma_3'$ are any edges on the graph of the Bratteli diagram for $B_{i,j} \subset B_{i+1,j}$, $\sigma_2$ is any edge on the graph of the Bratteli diagram for $B_{i,j+1} \subset B_{i+1,j+1}$, and $\sigma_4$ is any edge on the graph of the Bratteli diagram for $B_{i+1,j} \subset B_{i+1,j+1}$, such that $s(\sigma_2) = r(\sigma_1) = r(\sigma_1')$ and $s(\sigma_4) = r(\sigma_3) = r(\sigma_3')$.
Equation (\ref{eqn:condition_for_commuting_square}) is easily verified for both connections (\ref{eqn:connection}), (\ref{eqn:inverse_graph_connection}) using equations (\ref{eqn:typeI_frame}) and (\ref{eqn:typeII_frame}) and the fact that $[3]$ is the Perron-Frobenius eigenvalue for $\mathcal{G}$. This computation is essentially the algebraic verification of the first diagrammatic relation given in (\ref{braiding_relations1-3}).
\hfill
$\Box$

Then as in \cite{evans/kawahigashi:1994}, we define the Jones projections in $B_{i,j}$, for $i=1,2,\ldots$, by:
\begin{equation}\label{def:e_i}
e_{i-1} = \sum_{\stackrel{|\zeta|=j, |\zeta'| = i-2}{\scriptscriptstyle{|\gamma'|=|\eta'|=1}}} \frac{1}{[3]} \frac{\sqrt{\phi_{r(\gamma')}\phi_{r(\eta')}}}{\phi_{r(\zeta')}} \; (\zeta \cdot \zeta' \cdot \gamma' \cdot \widetilde{\gamma'}, \zeta \cdot \zeta' \cdot \eta' \cdot \widetilde{\eta'})
\end{equation}
where $\widetilde{\xi}$ denotes the reverse edge of $\xi$.
Let $E_{M_{i-1}}$ be the conditional expectation from $B_{i+1,\infty}$ onto $B_{i,\infty}$ with respect to the trace. For $x \in B_{i+1,j}$, $E_{M_{i-1}}(x)$ is given by the conditional expectation of $x$ onto $B_{i,j}$, because of Lemma \ref{Lemma:commuting_square_of_Bij's}. Clearly $e_l x = x e_l$, for $x \in B_{l-1, \infty}$, since $x$ and $e_l$ live on distinct parts of the Bratteli diagram.
It can be shown that $e_l x e_l = E_{M_{l-1}}(x) e_l$ for all $x \in B_{l,\infty}$, and that $B_{{l+1}, \infty}$ is generated by $B_{l, \infty}$ and $e_l$. Then $e_l$ is the Jones projection for the basic construction $B_{l-1,\infty} \subset B_{l,\infty} \subset B_{l+1,\infty}$, $l=1,2,\ldots \;$.
By \cite[Prop. 1.2]{pimsner/popa:1988} if we set $N = B_{0,\infty}$ and $M = B_{1,\infty}$, the sequence $B_{0,\infty} \subset B_{1,\infty} \subset B_{2,\infty} \subset B_{3,\infty} \subset \cdots$ can be identified with the Jones tower $N \subset M \subset M_1 \subset M_2 \subset \cdots$. It was shown in \cite{evans/kawahigashi:1994} that for $\mathcal{G} = \mathcal{A}^{(n)}$, $n < \infty$, if $\ast$ is now the apex vertex $(0,0)$ of $\mathcal{A}^{(n)}$, then this subfactor is the same as Wenzl's subfactor in \cite{wenzl:1988} for $SU(3)$, and we have the following theorem from \cite{evans/kawahigashi:1994} (Theorems 3.3, 5.8 and Corollary 3.4):

\begin{Thm} \label{thm:Wenzl_subfactor_has_pg_the_01part_of_A}
In the double sequence $(B_{i,j})$ above for $\mathcal{G} = \mathcal{A}^{(n)}$ or $\mathcal{D}^{(n)}$, $n < \infty$, with $\ast$ the vertex with lowest Perron-Frobenius weight, we have $B_{0,\infty}' \cap B_{i,\infty} = B_{i,0}$, i.e. $N' \cap M_{i-1} = B_{i,0}$.
The principal graph for the above subfactors is given by the 01-part $\mathcal{G}_{01}$ of $\mathcal{G}$.
\end{Thm}

The connection will be called flat \cite{ocneanu:1988, ocneanu:1991} if any two elements $x \in B_{k,0}$ and $y \in B_{0,l}$ commute. This is equivalent to the relation
\begin{center}
\begin{minipage}[b]{16cm}
 \begin{minipage}[t]{6.5cm}
  \mbox{} \\
  \parbox[t]{1cm}{}
 \end{minipage}
 \begin{minipage}[t]{4cm}
  \mbox{} \\
   \includegraphics[width=40mm]{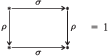}
 \end{minipage}
 \begin{minipage}[t]{5cm}
  \hfill
  \parbox[t]{2cm}{\begin{eqnarray}\label{flat_connection}\end{eqnarray}}
 \end{minipage}
\end{minipage}
\end{center}
for any paths $\sigma$, $\rho$ on the graphs $\mathcal{G}$ and $\widetilde{\mathcal{G}}$.

Then for graphs where the connection (\ref{eqn:connection}) is flat, the higher relative commutants are given by the $B_{k,0}$, that is, $B_{0,\infty}' \cap B_{k,\infty} = B_{k,0}$, by Ocneanu's compactness argument \cite{ocneanu:1991} in the setting of our $SU(3)$ subfactors. If $\mathcal{G}$ is a graph with flat connection, then the principal graph of the subfactor $B_{0,\infty} \subset B_{1,\infty}$ will be the 01-part $\mathcal{G}_{01}$ of $\mathcal{G}$.

Flatness of the connection for the $\mathcal{A}$, $\mathcal{D}$ graphs was shown in Theorem \ref{thm:Wenzl_subfactor_has_pg_the_01part_of_A}, where the distinguished vertex $\ast$ was chosen to be the vertex with lowest Perron-Frobenius weight. The flatness of the connection for the exceptional $\mathcal{E}$ graphs in not decided here. The determination of whether the connection is flat in these cases is a finite problem, involving checking the identity (\ref{flat_connection}) for diagrams of size $2d_{\mathcal{G}_{01}} \times 2(d_{\mathcal{G}}+3)$, where $d_{\mathcal{G}}$ is the depth of $\mathcal{G}$ and $d_{\mathcal{G}_{01}}$ is the depth of its 01-part $\mathcal{G}_{01}$. This is because for the vertical paths, the algebras $B_{l+1,j}$ are generated by $B_{l,j}$ and the Jones projection $e_l$ for all $l \geq d_{\mathcal{G}_{01}}$, and $e_l$ does not change its form under the change of basis using the connection. For the horizontal paths, by \cite[Lemma 4.7]{evans/pugh:2009ii} we see that the algebras $B_{i,l+1}$ are generated by $B_{i,l}$ and $U_{-l}$ for $l \geq d_{\mathcal{G}} + 3$, and the Hecke operators $U_{-l}$ do not change their form under the change of basis, as is shown in the proof of Theorem \ref{thm:type_I_subfactors_give_flat_planar_algebras} below.

We have not yet been able to determine whether or not the connection defined by (\ref{eqn:connection}), (\ref{eqn:inverse_graph_connection}) is flat for the $\mathcal{E}$ cases, where the vertex $\ast$ is chosen to be the vertex with lowest Perron-Frobenius weight, since the number of computations involved, though finite, is extremely large. We expect that this connection will be flat for the exceptional graphs $\mathcal{E}^{(8)}$, $\mathcal{E}_1^{(12)}$ and $\mathcal{E}^{(24)}$, since these graphs appear as the $M$-$N$ graphs for type I inclusions $N \subset M$. We expect that this connection will not be flat for the remaining exceptional graphs $\mathcal{E}_2^{(12)}$, $\mathcal{E}_4^{(12)}$ and $\mathcal{E}_5^{(12)}$ for any choice of distinguished vertex $\ast$. We also expect that the connection will not be flat for the $\mathcal{A}^{\ast}$, $\mathcal{D}^{\ast}$ graphs, for any choice of distinguished vertex $\ast$. The principal graph for the graphs with a non-flat connection is given by its flat part, which should be the type I parents given in \cite{evans/pugh:2009ii}.

\subsection{Flat $A_2$-$C^{\ast}$-planar algebra from $SU(3)$ $\mathcal{ADE}$ subfactors} \label{sect:flat_Planar_algebras_give_subfactors}

We will now associate a flat $A_2$-$C^{\ast}$-planar algebra $P$ to a double sequence $(B_{i,j})$ of finite dimensional algebras with a flat connection.

\begin{figure}[htb]
\begin{center}
  \includegraphics[width=135mm]{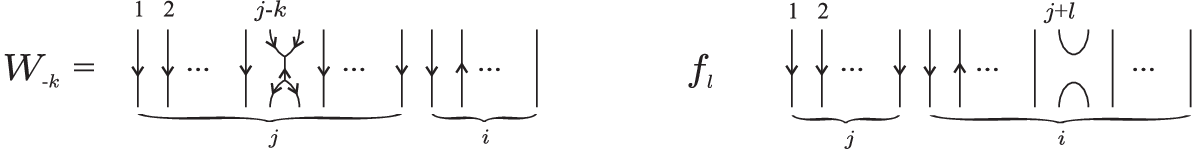}\\
 \caption{Tangles $W_{-k}$ and $f_l$} \label{fig:W_-k&f_l}
\end{center}
\end{figure}

We define the tangles $W_{-k}$, $k=0,\ldots,j-1$, and $f_l$, $l=1,\ldots,i$, in $\mathcal{P}_{i,j}(\varnothing)$ as in Figure \ref{fig:W_-k&f_l}, where the orientations of the strings without arrows depends on the parity of $i$ and $l$.

Let $\widetilde{\mathcal{G}}$ denote the graph $\mathcal{G}$ with all orientations reversed and let $P_{\sigma}$ be the space of closed paths on $\mathcal{G}$, $\widetilde{\mathcal{G}}$ which start at the distinguished vertex $\ast$, where a `$-$' denotes that an edge is on $\mathcal{G}$ and `$+$' denotes that an edge is on $\widetilde{\mathcal{G}}$. We will define a presenting map $Z:\mathcal{P}_{\sigma}(P) \rightarrow P_{\sigma}$ such that $P_{i,j} \cong B_{i,j}$, where we identify a path $\gamma_1 \cdot \gamma_2$ of length $2m$ in $P_{i,j}$ with the pair of paths $(\gamma_1,\widetilde{\gamma_2})$ of length $m$ (i.e. an element in $B_{i,j}$) by cutting the original path in half and reversing the path $\gamma_2$.
We define a $\ast$-operation on $P$ by $\gamma^{\ast} = \widetilde{\gamma} \in P_{\sigma^{\ast}}$ for $\gamma \in P_{\sigma}$. For $\gamma_1 \cdot \gamma_2 \in P_{i,j}$, $(\gamma_1 \cdot \gamma_2)^{\ast} = \widetilde{\gamma}_2 \cdot \widetilde{\gamma}_1$ which is mapped to $(\widetilde{\gamma}_2,\gamma_1) \in B_{i,j}$ under the isomorphism $P_{i,j} \cong B_{i,j}$. Note that $(\widetilde{\gamma}_2,\gamma_1) = (\gamma_1,\widetilde{\gamma}_2)^{\ast}$ in $B_{i,j}$, so the $\ast$-structure on $B$ is preserved under the isomorphism.

\begin{figure}[h]
\begin{center}
 \begin{minipage}[b]{13cm}
   \includegraphics[width=45mm]{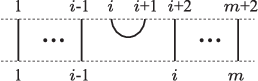}
  \hfill
   \includegraphics[width=45mm]{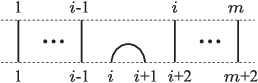}
 \end{minipage}
 \caption{Cup $\cup^{(i)}$ and cap $\cap^{(i)}$} \label{fig:cup(i)&cap(i)}
\end{center}
\end{figure}

\begin{figure}[h]
\begin{center}
 \begin{minipage}[b]{13cm}
   \includegraphics[width=45mm]{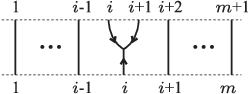}
  \hfill
   \includegraphics[width=45mm]{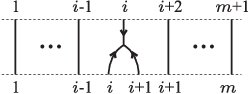}
 \end{minipage}
\caption{Y-forks $\curlyvee^{(i)}$ and $\overline{\curlyvee}^{(i)}$} \label{fig:Y-forks}
\end{center}
\end{figure}

\begin{figure}[h]
\begin{center}
 \begin{minipage}[b]{13cm}
   \includegraphics[width=45mm]{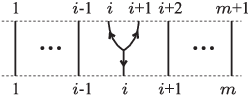}
  \hfill
   \includegraphics[width=45mm]{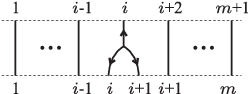}
 \end{minipage}
\caption{Inverted Y-forks $\curlywedge^{(i)}$ and $\overline{\curlywedge}^{(i)}$} \label{fig:inverted_Y-forks}
\end{center}
\end{figure}

Let $T$ be a labelled tangle in $\mathcal{P}_{\sigma}$ with $m$ internal discs $D_k$ with pattern $\sigma_k$ and labels $x_k \in P_{\sigma_k}$, $k=1, \ldots, m$. We define $Z(T)$ as follows. First, convert all the discs $D_k$ to rectangles (including the outer disc) so that its edges are parallel to the $x,y$-axes, and such that all the vertices on its boundary lie along the top edge of the rectangle.
Next, isotope the strings of $T$ so that each horizontal strip only contains one of the following elements: a rectangle with label $x_k$, a cup, a cap, a Y-fork, or an inverted Y-fork (see Figures \ref{fig:cup(i)&cap(i)}, \ref{fig:Y-forks} and \ref{fig:inverted_Y-forks}).
For a tangle $T \in \mathcal{P}_{\sigma}$ with $l$ horizontal strips $s_l$, where $s_1$ is the highest strip, $s_2$ the strip immediately below it, and so on, we define $\widetilde{Z}(T)=Z(s_1) Z(s_2) \cdots Z(s_l)$, which will be an element of $P_{\sigma}$. We then define $Z(T)$ by $Z(T) = \widetilde{Z}(T)$ if $\gamma$ is a path of odd length, and $Z(T) = \sqrt{\phi_{s(\gamma_1)}}/\sqrt{\phi_{r(\gamma_k)}} \widetilde{Z}(T)$ if $\gamma = \gamma_1 \cdot \gamma_2 \cdots \gamma_{2k}$ is a path of length $2k$, where the $\gamma_i$ are edges on $\mathcal{G}$ or $\widetilde{\mathcal{G}}$. Note that we have $P_{\sigma}^{\overline{0}} = P_{\sigma}^{\overline{1}} = P_{\sigma}^{\overline{2}}$.
This algebra is normalized in the sense that for the empty tangle $\bigcirc$, $Z(\bigcirc)=1$.
We will need to show that this definition only depends on $T$, and not on the decomposition of $T$ into horizontal strips.

Let $C$ be the set of all strips containing one of these elements except for a labelled rectangle. We will use the following notation for elements of $C$, as shown in Figures \ref{fig:cup(i)&cap(i)}, \ref{fig:Y-forks} and \ref{fig:inverted_Y-forks}: A strip containing a cup, cap will be $\cup^{(i)}$, $\cap^{(i)}$ respectively, where there are $i-1$ vertical strings to the left of the cup or cap. A strip containing an incoming Y-fork, inverted Y-fork will be $\curlyvee^{(i)}$, $\curlywedge^{(i)}$ respectively, where there are $i-1$ vertical strings to the left of the (inverted) Y-fork. A bar will denote that it is an outgoing (inverted) Y-fork.

For an element $c \in C$ we have sign strings $\sigma_1$, $\sigma_2$ given by the endpoints of the strings along the top, bottom edge respectively of the strip (we will call these endpoints vertices), where, along the top edge `$+$' is given by a sink and `$-$' by a source, and along the bottom edge `$+$' is given by a source and `$-$' by a sink.
The leftmost region of the strip $c$ corresponds to the vertex $\ast$ of $\mathcal{G}$, and each vertex along the top (or bottom) with downwards, upwards orientation respectively, corresponds to an edge on $\mathcal{G}$, $\widetilde{\mathcal{G}}$ respectively. Then the top, bottom edge of the strip is labelled by elements in $P_{\sigma_1}$, $P_{\sigma_2}$ respectively, which start at $\ast$. Then $Z(c)$ defines an operator $M_c \in \textrm{End}(P_{\sigma_2},P_{\sigma_1})$ as follows.

For a cup $\cup^{(i)}$, and paths $\alpha = \alpha_1 \cdot \alpha_2 \cdots \alpha_j$, $\beta = \beta_1 \cdots \beta_{j+2}$,
\begin{equation}\label{M_cup}
\left( M_{\cup^{(i)}} \right)_{\alpha, \beta} = \delta_{\alpha_1, \beta_1} \delta_{\alpha_2, \beta_2} \cdots \delta_{\alpha_{i-1}, \beta_{i-1}} \delta_{\alpha_i, \beta_{i+2}} \delta_{\alpha_{i+1}, \beta_{i+3}} \cdots \delta_{\alpha_m, \beta_{m+2}} \delta_{\widetilde{\beta}_i, \beta_{i+1}} \frac{\sqrt{\phi_{r(\beta_i)}}}{\sqrt{\phi_{s(\beta_i)}}}.
\end{equation}

For a cap $\cap^{(i)}$,
\begin{equation}\label{M_cap}
M_{\cap^{(i)}} = M_{\cup^{(i)}}^{\ast}.
\end{equation}

For an incoming (inverted) Y-fork $\curlyvee^{(i)}$ or $\curlywedge^{(i)}$,
\begin{eqnarray}
\left( M_{\curlyvee^{(i)}} \right)_{\alpha, \beta} & = & \delta_{\alpha_1, \beta_1} \cdots \delta_{\alpha_{i-1}, \beta_{i-1}} \delta_{\alpha_{i+1}, \beta_{i+2}} \cdots \delta_{\alpha_m, \beta_{m+1}} \frac{1}{\sqrt{\phi_{s(\alpha_i)}\phi_{r(\alpha_i)}}} W(\triangle^{(\widetilde{\alpha}_i, \beta_i, \beta_{i+1})}), \nonumber \\
& & \label{M_Yfork(in)} \\
\left( M_{\curlywedge^{(i)}} \right)_{\alpha, \beta} & = & \delta_{\alpha_1, \beta_1} \cdots \delta_{\alpha_{i-1}, \beta_{i-1}} \delta_{\alpha_{i+2}, \beta_{i+1}} \cdots \delta_{\alpha_{m+1}, \beta_m} \frac{1}{\sqrt{\phi_{s(\beta_i)}\phi_{r(\beta_i)}}} \overline{W(\triangle^{(\beta_i, \widetilde{\alpha}_{i+1}, \widetilde{\alpha}_i)})}, \nonumber \\
& &  \label{M_invYfork(in)}
\end{eqnarray}
where $W$ is a cell system on $\mathcal{G}$ satisfying (\ref{eqn:typeI_frame}) and (\ref{eqn:typeII_frame}).

For an outgoing (inverted) Y-fork $\overline{\curlyvee}^{(i)}$ or $\overline{\curlywedge}^{(i)}$,
\begin{eqnarray}
M_{\overline{\curlyvee}^{(i)}} & = & M_{\curlywedge^{(i)}}^{\ast}, \label{M_Yfork(out)} \\
M_{\overline{\curlywedge}^{(i)}} & = & M_{\curlyvee^{(i)}}^{\ast}. \label{M_invYfork(out)}
\end{eqnarray}

\begin{figure}[tb]
\begin{center}
  \includegraphics[width=90mm]{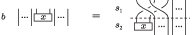}\\
 \caption{$Z(b)$ for horizontal strip $b$ containing a rectangle} \label{fig:Z(b)}
\end{center}
\end{figure}

For a strip $b$ containing a rectangle with label $x = \gamma$, where $\gamma$ is a single path in $P_{\sigma}$, we define the operator $M_{b} = Z(b)$ as follows. Let $p$, $p'$ be the number of vertical strings to the left, right respectively of the rectangle in strip $b$, with orientations given by the sign strings $\sigma^{(p)}$, $\sigma^{(p')}$ respectively. We attach trivial tails $\mu$ of length $p$ to $x$, where $\mu$ has edges on $\mathcal{G}$, $\widetilde{\mathcal{G}}$ as dictated by the sign string $\sigma^{(p)}$, so that we have a sum $\sum_{\mu} \gamma \cdot \mu$ of paths in the basis given by the sign string $\sigma \sigma^{(p)}$. We use the connection to transform this to a linear combination of paths in the basis given by the sign string $\sigma^{(p)} \sigma$.
By flatness of the connection on $\mathcal{G}$, this will be an element of the form $\sum_{\zeta, \mu} p_{\zeta} \mu \cdot \zeta$, where $p_{\zeta} \in \mathbb{C}$ are given by the connection, $\zeta$ are paths in $P_{\sigma}$, and $\mu$ are again paths in $P_{\sigma^{(p)}}$. We then add trivial tails $\nu$ of length $p'$ to this element, where $\nu$ has edges on $\mathcal{G}$, $\widetilde{\mathcal{G}}$ as dictated by the sign string $\sigma^{(p')}$. This gives an element $\sum_{\zeta, \mu, \nu} p_{\zeta} \mu \cdot \zeta \cdot \nu$, which is an element in $P_{\sigma^{(p)} \sigma \sigma^{(p')}}$. Then we define $Z(b) \in \mathrm{End}(P_{\sigma^{(p)} \sigma^{(p')}},P_{\sigma^{(p)} \sigma \sigma^{(p')}})$ to be $\sum_{\zeta, \mu, \nu} p_{\zeta} (\mu \cdot \zeta \cdot \nu, \mu \cdot \nu)$.
We extend this definition of $Z(b)$ linearly to strips $b$ where the label $x$ is a linear combination of paths in $P_{\sigma}$. This definition means that $Z(b)$ is defined as the product $Z(s_1)Z(s_2)$, where $s_1$, $s_2$ are the horizontal strips on the right hand side of Figure \ref{fig:Z(b)}. Thus we see that $P$ is a flat $A_2$-planar algebra.

The following theorem shows that for an $\mathcal{ADE}$ graph with a flat connection, $P = \bigcup_{\sigma} P_{\sigma}$ is a flat $A_2$-$C^{\ast}$-planar algebra, such that the subalgebra $\bigcup_{i,j} P_{i,j}$ is given by the subfactor double complex $(B_{i,j})$.

\begin{Thm} \label{thm:type_I_subfactors_give_flat_planar_algebras}
Let $\mathcal{G}$ be an $\mathcal{ADE}$ graph such that the connections (\ref{eqn:connection}), (\ref{eqn:inverse_graph_connection}) are flat. The above definition of $Z(T)$ for any $A_2$-planar tangle $T$ makes $P = \bigcup_{\sigma} P_{\sigma}$ into a flat $A_2$-$C^{\ast}$-planar algebra, such that $P_{i,j} \cong B_{i,j}$, with $\mathrm{dim} \left( P_0^{\overline{0}} \right) = \mathrm{dim} \left( P_0^{\overline{1}} \right) = \mathrm{dim} \left( P_0^{\overline{2}} \right) = 1$. This $A_2$-$C^{\ast}$-planar algebra has parameter $\alpha = [3]$ (the Perron-Frobenius eigenvalue for $\mathcal{G}$), and $Z(I_{\sigma\sigma^{\ast}}(x)) = x$, where $I_{\sigma\sigma^{\ast}}(x)$ is the tangle $I_{\sigma\sigma^{\ast}}$ with $x \in P_{\sigma\sigma^{\ast}}$ as the insertion in its inner disc. For $x \in P_{i,j}$, $i,j \geq 0$, we have \\
\begin{minipage}[b]{16cm}
 \begin{minipage}[t]{1cm}
  \parbox[t]{1cm}{$(i)$}
 \end{minipage}
 \begin{minipage}[t]{6cm}
  $Z(W_{-k}) = U_{-k}$, \qquad $k \geq 0,$
 \end{minipage}
\end{minipage} \\
\begin{minipage}[t]{16cm}
 \begin{minipage}[t]{1cm}
  \parbox[t]{1cm}{$(ii)$}
 \end{minipage}
  \begin{minipage}[t]{6cm}
  $Z(f_{l}) = \alpha e_l$, \qquad $l \geq 1,$
 \end{minipage}
\end{minipage} \\
\begin{minipage}[t]{16cm}
 \begin{minipage}[t]{1cm}
  \mbox{} \\
  \parbox[t]{1cm}{$(iii)$}
 \end{minipage}
 \begin{minipage}[t]{7cm}
  \mbox{} \\
   \includegraphics[width=60mm]{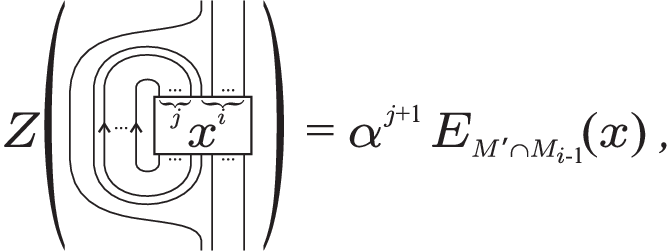}
 \end{minipage}
 \begin{minipage}[t]{7cm}
  \mbox{} \\
   \includegraphics[width=60mm]{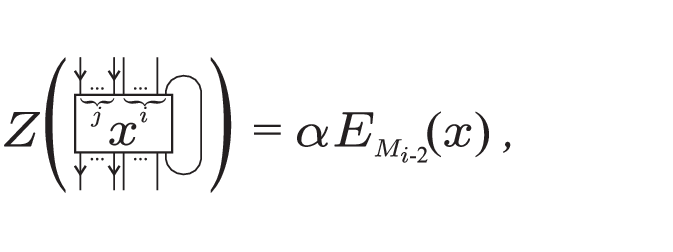}
 \end{minipage}
\end{minipage} \\
\mbox{} \\
\begin{minipage}[t]{16cm}
 \begin{minipage}[t]{1cm}
  \mbox{} \\
  \parbox[t]{1cm}{$(iv)$}
 \end{minipage}
 \begin{minipage}[t]{6cm}
  \mbox{} \\
   \includegraphics[width=45mm]{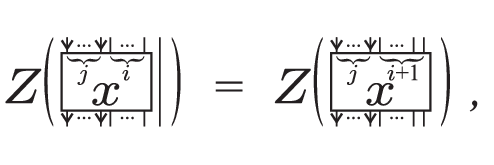}
 \end{minipage}
 \begin{minipage}[t]{6.5cm}
  \mbox{} \\
   \includegraphics[width=45mm]{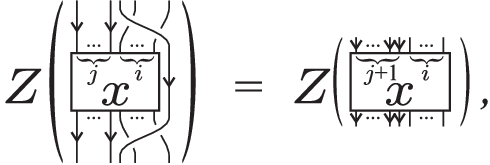}
 \end{minipage}
\end{minipage} \\
\mbox{} \\
\begin{minipage}[t]{16cm}
 \begin{minipage}[t]{1cm}
  \mbox{} \\
  \parbox[t]{1cm}{$(v)$}
 \end{minipage}
 \begin{minipage}[t]{7cm}
  \mbox{} \\
   \includegraphics[width=43mm]{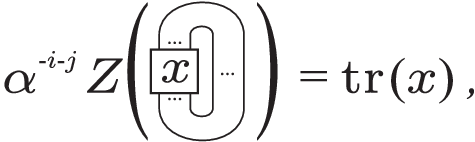}
 \end{minipage}
\end{minipage} \\ \\
In the first equation of $(iii)$ the first $j+1$ vertices along the top and bottom of the rectangle are joined by loops, and the second equation only holds for $i \neq 0$. In the first, second equation of $(iv)$ respectively, the $x$ on the right hand side is considered as an element of $P_{i+1,j}$, $P_{i,j+1}$ respectively.
\end{Thm}

\noindent
\emph{Proof:}
First we show that $Z(T)$ does not change if the labelled tangle is changed by isotopy of the strings. We use the following notation $\partial_{\alpha_i, \beta_j}^{\alpha_{i+k}, \beta_{j+k}} := \delta_{\alpha_i, \beta_j} \delta_{\alpha_{i+1}, \beta_{j+1}} \cdots \delta_{\alpha_{i+k}, \beta_{j+k}}$. The identities are simply a consequence of the identification in Figure \ref{fig:Oc-Kup} of the Ocneanu cells with trivalent vertices, and of cups and caps with the Perron-Frobenius weights. \\

\begin{figure}[b]
\begin{center}
 \begin{minipage}[b]{9cm}
   \includegraphics[width=90mm]{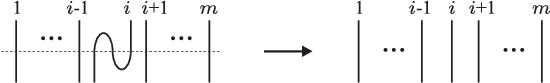}
 \end{minipage} \\
 \mbox{} \\
 \begin{minipage}[b]{9cm}
   \includegraphics[width=90mm]{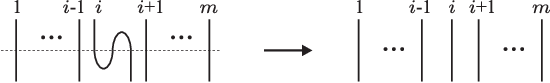}
 \end{minipage}
 \caption{Two cup-cap simplifications} \label{fig:cup-cap_simplifications}
\end{center}
\end{figure}

\noindent
\emph{Case (1)- Topological moves:}
We consider the cup-cap simplifications (which Kauffman calls Move Zero in \cite{kauffman:2001}) shown in Figure \ref{fig:cup-cap_simplifications}.
For the first cup-cap simplification of Figure \ref{fig:cup-cap_simplifications} we have \\
\begin{eqnarray}
\lefteqn{ \left( M_{\cup^{(i+1)}} M_{\cap^{(i)}} \right)_{\alpha, \beta} = \sum_{\gamma} \left( M_{\cup^{(i+1)}} \right)_{\alpha, \gamma} \left( M_{\cup^{(i)}} \right)_{\beta, \gamma} } \nonumber \\
& = & \sum_{\gamma} \partial_{\alpha_1, \gamma_1}^{\alpha_i, \gamma_i} \partial_{\alpha_{i+1}, \gamma_{i+3}}^{\alpha_m, \gamma_{m+2}} \delta_{\widetilde{\gamma_{i+1}}, \gamma_{i+2}} \frac{\sqrt{\phi_{r(\gamma_{i+1})}}}{\sqrt{\phi_{s(\gamma_{i+1})}}} \, \partial_{\beta_1, \gamma_1}^{\beta_{i-1}, \gamma_{i-1}} \partial_{\beta_i, \gamma_{i+2}}^{\beta_m, \gamma_{m+2}} \delta_{\widetilde{\gamma_i}, \gamma_{i+1}} \frac{\sqrt{\phi_{r(\gamma_i)}}}{\sqrt{\phi_{s(\gamma_i)}}} \;\; = \;\; \delta_{\alpha, \beta}. \qquad \label{cup-cap_simplification}
\end{eqnarray}
The second simplification in Figure \ref{fig:cup-cap_simplifications} follows from the first, since
\begin{equation}\label{cup-cap_simplification2}
M_{\cup^{(i)}} M_{\cap^{(i+1)}} = \left( M_{\cup^{(i+1)}} M_{\cap^{(i)}} \right)^T = \mathbf{1}.
\end{equation}

\begin{figure}[tb]
\begin{center}
  \includegraphics[width=110mm]{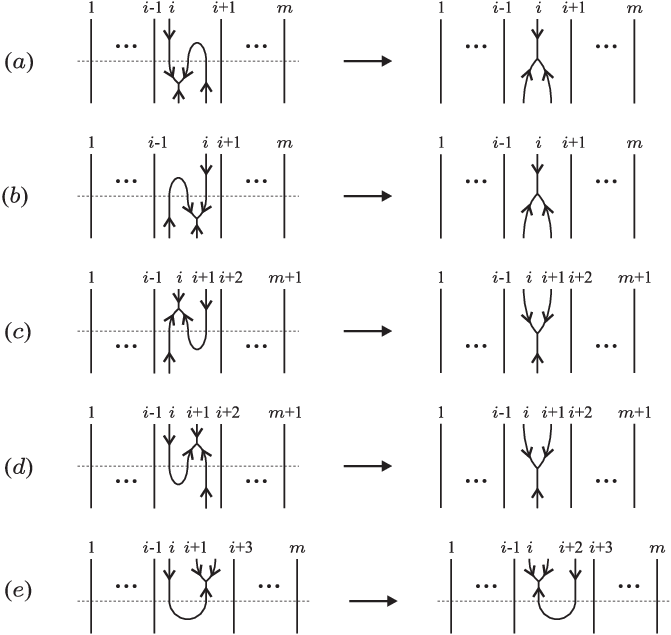}
 \caption{Isotopies involving an incoming trivalent vertex}\label{isotopies_inolving_trivalent_vertices}
\end{center}
\end{figure}

\begin{figure}[tb]
\begin{center}
  \includegraphics[width=90mm]{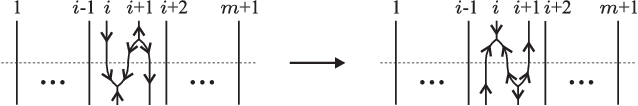}
 \caption{An isotopy involving an incoming and outgoing trivalent vertex} \label{fig:isotopy_inovlving_2_trivalent_vertices}
\end{center}
\end{figure}

\noindent
\emph{Case (2)- Isotopies involving incoming trivalent vertices:}
We require the identities of Figure \ref{isotopies_inolving_trivalent_vertices}.
For $(a)$ we verify that $\left(M_{\curlyvee^{(i)}} M_{\cap^{(i+1)}} \right)_{\alpha, \beta} = \left( M_{\curlywedge^{(i)}} \right)_{\alpha, \beta}$, and similarly for the identities $(b)$, $(c)$ and $(d)$. For $(e)$ we need to verify that $\left(M_{\cup^{(i-1)}} M_{\curlyvee^{(i)}} \right)_{\alpha, \beta} = \left(M_{\cup^{(i-1)}} M_{\curlyvee^{(i-1)}} \right)_{\alpha, \beta}$.
The corresponding identities for outgoing trivalent vertices hold in the same way. Then the identity in Figure \ref{fig:isotopy_inovlving_2_trivalent_vertices} follows from the cup-cap simplifications and identities (a)-(e) for incoming and outgoing trivalent vertices. \\

\noindent
\emph{Kuperberg relations:}
Before checking isotopies that involve rectangles, we will show that the Kuperberg relations K1-K3 are satisfied. For K1, a closed loop gives\\
\begin{equation}\label{proof_of_K1}
\left( M_{\cup^{(i)}} M_{\cap^{(i)}} \right)_{\alpha, \beta} = \sum_{\gamma} \left( M_{\cup^{(i)}} \right)_{\alpha, \gamma} \left( M_{\cup^{(i)}} \right)_{\beta, \gamma} = \delta_{\alpha, \beta} \sum_{\stackrel{\gamma_i:}{\scriptscriptstyle{s(\gamma_i)=r(\alpha_{i-1})}}} \frac{\phi_{r(\gamma_i)}}{\phi_{s(\gamma_i)}}  = \delta_{\alpha, \beta} [3],
\end{equation}
by the Perron-Frobenius eigenvalue equation $\Lambda x = [3] x$, where $x = (\phi_{v})_v$, and $\Lambda$ is $\Delta_{\mathcal{G}}$ or $\Delta_{\mathcal{G}}^T$ depending on whether the loop has anticlockwise, clockwise orientation respectively.
Relations K2 and K3 are essentially Ocneanu's type I, II formulas (\ref{eqn:typeI_frame}), (\ref{eqn:typeII_frame}) respectively. \\

\noindent
\emph{Property $(ii)$ and the connection:}
We obtain
$$Z(W_{-k}) = \sum_{\stackrel{|\zeta_1|=j-2-k, |\zeta'|= i}{\scriptscriptstyle{|\gamma_i|=|\eta_i|=1, |\zeta_2|=k}}} \mathcal{U}^{\gamma_2,\eta_2}_{\gamma_1,\eta_1}  \; \zeta_1 \cdot \gamma_1 \cdot \eta_1 \cdot \zeta_2 \cdot \zeta' \cdot \widetilde{\zeta'} \cdot \widetilde{\zeta_2} \cdot \widetilde{\eta_2} \cdot \widetilde{\gamma_2} \cdot \widetilde{\zeta_1},$$
and we identify the path
$\zeta_1 \cdot \gamma_1 \cdot \eta_1 \cdot \zeta_2 \cdot \zeta' \cdot \widetilde{\zeta'} \cdot \widetilde{\zeta_2} \cdot \widetilde{\eta_2} \cdot \widetilde{\gamma_2} \cdot \widetilde{\zeta_1} \in P_{\sigma \sigma^{\ast}}$ with the matrix unit $(\zeta_1 \cdot \gamma_1 \cdot \eta_1 \cdot \zeta_2 \cdot \zeta', \zeta_1 \cdot \gamma_2 \cdot \eta_2 \cdot \zeta_2 \cdot \zeta') \in \mathrm{End}(P_{\sigma},P_{\sigma})$.
The property $(ii)$ in the statement of the theorem follows from (\ref{eqn:HeckeRep}) and the definition of $U_{-k}$.
Since $U_{-k}$ is given by the tangle $W_{-k}$, we see that the partial braiding defined in (\ref{braiding2}) gives the connection, where (\ref{eqn:connection}) is given by \includegraphics[width=5mm]{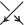} and (\ref{eqn:inverse_graph_connection}) is given by \includegraphics[width=5mm]{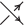}. For the latter connection, which involves the reverse graph $\widetilde{\mathcal{G}}$, if
$$
\begin{array}{ccc}
a & \rightarrow & b \\
\downarrow & & \downarrow \\
c & \rightarrow & d
\end{array}
$$
is a connection on the graph $\mathcal{G}$, then
\begin{center}
 \begin{minipage}[b]{10.5cm}
  \begin{minipage}{2.6cm}
   $$
   \begin{array}{ccc}
   c & \rightarrow & d \\
   \downarrow & & \downarrow \\
   a & \rightarrow & b
   \end{array} =
   $$
  \end{minipage}
  \begin{minipage}{1.1cm}
 \includegraphics[width=11mm]{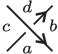}
  \end{minipage}
  \hspace{-0.5cm}
  \begin{minipage}{1cm}
   $$
   \begin{array}{c}
    \\
    \\
   \end{array} =
   $$
  \end{minipage}
  \begin{minipage}{2cm}
 \includegraphics[width=20mm]{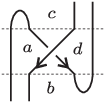}
  \end{minipage}
  \begin{minipage}{3.55cm}
   $$
   = \sqrt{\frac{\phi_a \phi_d}{\phi_b \phi_c}} \;
   \overline{\begin{array}{ccc}
   a & \rightarrow & b \\
   \downarrow & & \downarrow \\
   c & \rightarrow & d
   \end{array}}.
   $$
  \end{minipage}
 \end{minipage}
\end{center}

So we have that $Z(T)$ is invariant under all isotopies that only involve strings (and the partial braiding). This shows that the operators $U_{-k}$ do not change their form under the change of basis using the connection, since
\begin{center}
   \includegraphics[width=35mm]{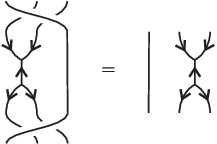}
\end{center}
Note that we have not used the fact that the connection is flat yet, so the operators $U_{-k}$ do not change their form under the change of basis for any of the $SU(3)$ $\mathcal{ADE}$ graphs. \\

\begin{figure}[tb]
\begin{center}
   \includegraphics[width=90mm]{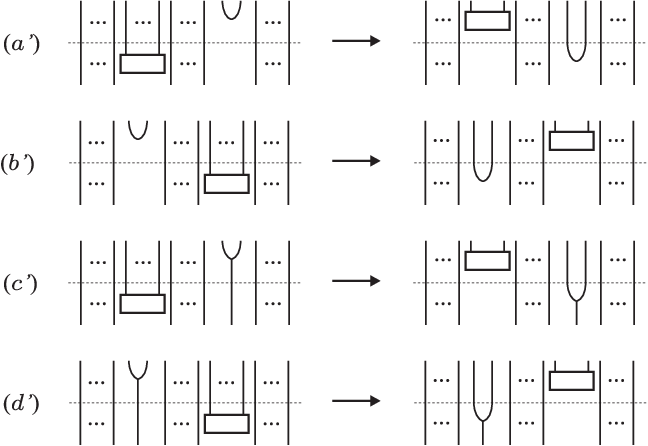}
 \caption{Isotopies involving rectangles} \label{fig:isotopies_involving_boxes}
\end{center}
\end{figure}

\begin{figure}[tb]
\begin{center}
  \includegraphics[width=100mm]{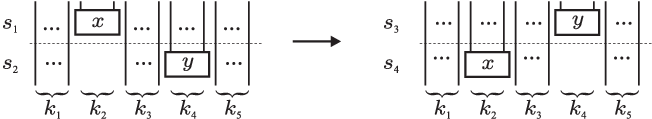}
 \caption{An isotopy involving two rectangles} \label{fig:commuting_boxes}
\end{center}
\end{figure}

\noindent
\emph{Case (3)- Isotopies that involve rectangles:}
We need to check invariance as in Figure \ref{fig:isotopies_involving_boxes}.
For $(a')$, pulling a cup down to the right of a rectangle $b$ is trivial since $M_{\cup}$ commutes with $M_b$ (since $b$, $\cup$ are localized on separate parts of the Bratteli diagram). Now consider $(b')$. We have for the left hand side
\begin{center}
   \includegraphics[width=140mm]{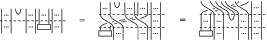}
\end{center}
where the first equality is the definition of $Z(s)$ for a horizontal strip $s$ containing a rectangle $b$ with through strings on its left, and the second equality follows since $Z$ is invariant under all isotopies that only involve strings and the partial braiding. Similarly, for the right hand side we obtain
\begin{center}
   \includegraphics[width=35mm]{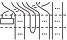}
\end{center}
and the result follows from $(a')$.
The situations for $(c')$, $(d')$ are similar to $(a')$, $(b')$. We also have the isotopy in Figure \ref{fig:commuting_boxes}.
Let $x \in P_{\sigma_1}$, $y \in P_{\sigma_2}$ be given by the paths $\alpha_1$, $\alpha_2$ respectively, of lengths $k_2$, $k_4$ respectively. The case for general elements $x \in P_{\sigma_1}$, $y \in P_{\sigma_2}$ follows by linearity.
Consider first the case where $k_1 = k_3 = k_5 = 0$. Then $Z(s_1)Z(s_2) = (\alpha_1 \cdot \alpha_2, \ast) \in \mathrm{End}(P_0,P_{\sigma_1\sigma_2})$. For the right-hand side we have $Z(s_4) = (\alpha_1, \ast) \in \mathrm{End}(P_0,P_{\sigma_1})$ and $Z(s_3) = \sum_{\mu, \alpha_2'} p_{\alpha_2'} (\mu \cdot \alpha_2', \mu) \in \mathrm{End}(P_{\sigma_1},P_{\sigma_1\sigma_2})$, where $p_{\alpha_2'} \in \mathbb{C}$ are given by the connection. By the flatness condition (\ref{flat_connection}), $p_{\alpha_2'} = \delta_{\alpha_2',\alpha_2}$, so we have $Z(s_3)Z(s_4) = \sum_{\mu} (\mu \cdot \alpha_2, \mu) (\alpha_1, \ast) = (\alpha_1 \cdot \alpha_2, \ast) = Z(s_1)Z(s_2)$. The cases where $k_1$, $k_3$, $k_5$ are non-zero follow similarly. \\

\begin{figure}[tb]
\begin{center}
  \includegraphics[width=80mm]{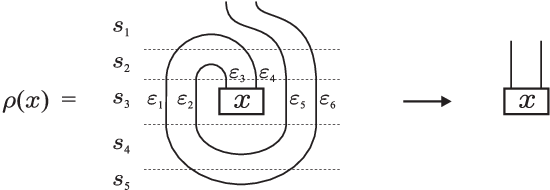}
 \caption{Rotation of internal rectangles by $2 \pi$} \label{fig:rotation_of_internal_boxes}
\end{center}
\end{figure}

\noindent
\emph{Case (4)- Rotational invariance:}
The other isotopy that needs to be checked is the rotation of internal rectangles by $2 \pi$. We illustrate the case where rectangle $b$ has $k_b = 2$ vertices along its top edge in Figure \ref{fig:rotation_of_internal_boxes}. We have divided $\rho(x)$ into horizontal strips $s_1,\ldots,s_5$.

Let $x \in P_{\sigma_b}$ be the label of the rectangle $b$, where $x$ is the single path $\gamma$ of length $k_b$. The strip $s_3$ containing the rectangle $b$ gives $Y = Z(s_3) = \sum_{\mu,\zeta,\nu} p_{\zeta} (\mu \cdot \zeta \cdot \nu, \mu \cdot \nu)$, where $p_{\zeta} \in \mathbb{C}$ are given by the connection, and $\mu$, $\nu$ are paths of length $k_b$ with edges on the graphs $\mathcal{G}$ or $\widetilde{\mathcal{G}}$ as dictated by the sign strings $\sigma_b^{\ast}$, $\sigma_b$ respectively.

\begin{figure}[tb]
\begin{center}
  \includegraphics[width=35mm]{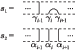}
 \caption{Horizontal strips $s_1$, $s_2$} \label{fig:fig36}
\end{center}
\end{figure}

For a horizontal strip $s_1$ and strip $s_2$ immediately below it, an entry in the operator $Z(s') = Z(s_1)Z(s_2)$ is only defined when the path corresponding to the bottom edge of the strip $s_1$ is equal to the path given by the top edge of $s_2$. So for example, for the two strips $s_1$, $s_2$ in Figure \ref{fig:fig36}, even though there are non-zero entries in $Z(s_2)$ for any path $\alpha = \alpha_1 \cdot \alpha_2 \cdots$, the entries in $Z(s')$ will be zero unless edge $\alpha_i$ is the reverse edge $\widetilde{\alpha_{i+1}}$ of $\alpha_{i+1}$ since the entries in $Z(s_2)$ are only non-zero for the paths $\gamma = \gamma_1 \cdot \gamma_2 \cdots$ such that $\gamma_i = \widetilde{\gamma_{i+1}}$.

Let $\varepsilon = \varepsilon_1 \cdot \varepsilon_2 \cdots \varepsilon_{3k_b}$, $\varepsilon' = \varepsilon_1' \cdot \varepsilon_2' \cdots \varepsilon_{k_b}' \cdot \varepsilon_{2k_b}' \cdot \varepsilon_{2k_b+1}' \cdots \varepsilon_{3k_b}'$ be two paths which label the indices for $Y$. For simplicity we consider the case $k_b = 2$ as in Figure \ref{fig:rotation_of_internal_boxes}. By considering the horizontal strip $s_3$ containing the rectangle, we see that $Y_{\varepsilon, \varepsilon'} = 0$ unless $\varepsilon_i = \varepsilon_i'$ for $i=1,2,5,6$. We see that in $\rho(x)$, $\varepsilon_1$ is the same string as $\varepsilon_4$ and $\varepsilon_6$, but that $\varepsilon_4$ has the opposite orientation to $\varepsilon_1$ and $\varepsilon_6$. We define the operator $\widehat{Y}$ by $\widehat{Y}_{\varepsilon, \varepsilon'} = 0$ unless $\varepsilon_1 = \widetilde{\varepsilon_4} = \varepsilon_6$ and $\varepsilon_2 = \widetilde{\varepsilon_3} = \widetilde{\varepsilon_5}$, and $\widehat{Y}_{\varepsilon, \varepsilon'} = Y_{\varepsilon, \varepsilon'}$ otherwise. Then \begin{eqnarray*}
\rho(x) & = & M_{\cup^{(1)}} M_{\cup^{(2)}} \cdots M_{\cup^{(k_b)}} Y M_{\cap^{(k_b)}} M_{\cap^{(k_b-1)}} \cdots M_{\cap^{(1)}} \\
& = & M_{\cup^{(1)}} M_{\cup^{(2)}} \cdots M_{\cup^{(k_b)}} \widehat{Y} M_{\cap^{(k_b)}} M_{\cap^{(k_b-1)}} \cdots M_{\cap^{(1)}}.
\end{eqnarray*}
For any two paths $\varepsilon$ and $\varepsilon'$ such that $\widehat{Y}_{\varepsilon, \varepsilon'}$ is non-zero, the caps contribute a scalar factor $\sqrt{\phi_{r(\varepsilon_2)}}/\sqrt{\phi_{s(\varepsilon_1)}} = \sqrt{\phi_{s(\varepsilon_1)}}/\sqrt{\phi_{s(\varepsilon_1)}} = 1$, and similarly we have a scalar factor of 1 from the cups. Now $\varepsilon_1$ is an edge on $\mathcal{G}$ (or $\widetilde{\mathcal{G}}$) with $s(\varepsilon_1)=\ast$, and hence $\rho (x)$ is only non-zero for paths $\varepsilon_3 \cdot \varepsilon_4$ such that $s(\varepsilon_3) = \ast$. By the flatness of the connection on $\mathcal{G}$, the only path $\zeta$ starting from $\ast$ for which $p_{\zeta} \neq 0$ is $\gamma$, i.e. the original element $x$. Then the resulting operator given by $\rho(x)$ will have all entries 0 except for that for $\gamma$, and we have $\rho(x) = \gamma = x$.

Then $Z(T)$ is invariant under all isotopies of the tangle $T$. \\

\noindent
\emph{Properties $(i)$-$(v)$:}
Property $(i)$ follows from the definition of $e_i$ in (\ref{def:e_i}), and
property $(ii)$ has already been shown. Now consider property $(iii)$. We start with the first equation.
For any $x \in P_{i,j}$, the left hand side is equal to $Z(EL^{i,0}_{i,0} ER^{i,1}_{i,0} ER^{i,2}_{i,1} \cdots ER^{i,j}_{i,j-1} (x))$, and so gives is the conditional expectation of $x$ onto $P^{(1)}_{i,0}$ (see Section \ref{Sect:Basic_planar_tangles}).
We now show that $P^{(1)}_{i,0} = M' \cap M_{i-1}$. Embedding the subalgebra $P_{i,0}^{(1)} \subset P_{i,0}$ in $P_{i,\infty}$ we see that it lives on the last $i-1$ strings, with the rest all vertical through strings. Then $P_{i,0}^{(1)}$ clearly commutes with $M$, since the embedding of $M = P_{1,\infty}$ in $P_{i,\infty}$ has the last $i-1$ strings all vertical through strings, so we have $M' \cap M_{i-1} \supset P_{i,0}^{(1)}$.
For the opposite inclusion, we extend the double sequence $(B_{i,j})$ to the left to get

$$
\begin{array}{cccccccccc}
& & B_{0,0} & \subset & B_{0,1} & \subset & B_{0,2} & \subset & \cdots \qquad \longrightarrow & B_{0,\infty} \\
& & \cap & & \cap & & \cap & & & \cap \\
B_{1,-1} & \subset & B_{1,0} & \subset & B_{1,1} & \subset & B_{1,2} & \subset & \cdots \qquad \longrightarrow & B_{1,\infty} \\
\cap & & \cap & & \cap & & \cap & & & \cap \\
B_{2,-1} & \subset & B_{2,0} & \subset & B_{2,1} & \subset & B_{2,2} & \subset & \cdots \qquad \longrightarrow & B_{2,\infty} \\
\cap & & \cap & & \cap & & \cap & & & \cap \\
\vdots & & \vdots & & \vdots & & \vdots & & & \vdots
\end{array}
$$
Note that $B_{1,-1} = B_{0,0} = \mathbb{C}$. Since the connection is flat, by Ocneanu's compactness argument \cite{ocneanu:1991} we have $B_{1,\infty}' \cap B_{i,\infty} = B_{i,-1}$. Let $x = (\alpha_1, \alpha_2)$ be an element of $B_{i,-1}$. We embed $x$ in $B_{i,0}$ by adding trivial horizontal tails of length one, and using the connection we can write $x$ as $x' = \sum_{\mu} p_{\beta_1, \beta_2} (\mu \cdot \beta_1, \mu \cdot \beta_2)$, where $p_{\beta_1, \beta_2} \in \mathbb{C}$.
We see that $x' \in B_{i,0} = P_{i,0}$ is summed over all trivial edges $\mu$ of length 1 starting at $\ast$, and hence is given by $Z(T)$ for some $T \in \mathcal{P}_{i,0}$ which has a vertical through string from the first vertex along the top to the first vertex along the bottom, i.e $x \in P_{i,0}^{(1)}$. So $M' \cap M_{i-1} = B_{i,-1} \subset P_{i,0}^{(1)}$.

For the second equation of $(iii)$, if $x \in P_{i,\infty}$ then $x \rightarrow Z(E_{i-1,\infty}^{i,\infty}(x))$ is the conditional expectation onto $P_{i-1,\infty} = M_{i-2}$, and the result for $x \in P_{i,j}$ follows by Lemma \ref{Lemma:commuting_square_of_Bij's}.

Property $(iv)$ is clear. Finally, for $(v)$ let $x$ be an element $(\alpha, \beta)$, where $\alpha$, $\beta$ are paths of length $k$ on $\mathcal{G}$. Then
$$[3]^{-k} Z(\widehat{x}) = [3]^{-k} \delta_{\alpha, \beta} \frac{\phi_{r(\alpha)}}{\phi_{\ast}} \phi_{\ast}^2 = [3]^{-k} \delta_{\alpha, \beta} \phi_{r(\alpha)} = \textrm{tr}((\alpha, \beta)),$$
since $\phi_{\ast}=1$, where $\widehat{x}$ is the tangle defined by joining the last vertex along the top of $T$ to the last vertex along the bottom by a string which passes round the tangle on the right hand side, and joining the other vertices along the top to those on the bottom similarly.

\begin{figure}[b]
\begin{center}
  \includegraphics[width=70mm]{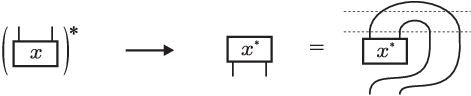}
 \caption{The $\ast$-structure on strips containing rectangles} \label{fig:fig86}
\end{center}
\end{figure}

To see the $\ast$-structure, note that under $\ast$ the order of the strips is reversed so that $(Z(s_1) Z(s_2) \cdots Z(s_l))^{\ast} = Z(s_l)^{\ast} Z(s_{l-1})^{\ast} \cdots Z(s_1)^{\ast}$.
For a strip containing a rectangle with label $x \in P_{\sigma}$ given by a path $\gamma$, $\ast$ sends the rectangle to the right hand side of Figure \ref{fig:fig86}. Since $s(\gamma) = r(\gamma)$, the caps contribute a coefficient $\sqrt{\phi_{r(\widetilde{\gamma})}}/\sqrt{\phi_{s(\widetilde{\gamma})}} = 1$ as required.
For $M_{\cup^{(i)}}$, the ratio $\sqrt{\phi_{r(\beta_i)}}/\sqrt{\phi_{s(\beta_i)}}$ does not change under reflection of the tangle and reversing the orientation, so that $(M_{\cup^{(i)}})^{\ast}$ is the conjugate transpose of $M_{\cup^{(i)}}$ as required, and similarly for $M_{\cap^{(i)}}$. Since the involution of the strip $\curlyvee^{(i)}$ containing an incoming trivalent vertex is $\overline{\curlywedge}^{(i)}$, whilst the involution of the strip $\curlywedge^{(i)}$ containing an incoming trivalent vertex is $\overline{\curlyvee}^{(i)}$, by (\ref{M_Yfork(out)}), $(M_{\overline{\curlyvee}^{(i)}})^{\ast}$ is the conjugate transpose of $M_{\curlywedge^{(i)}}$ and by (\ref{M_invYfork(out)}), $(M_{\overline{\curlywedge}^{(i)}})^{\ast}$ is the conjugate transpose of $M_{\curlyvee^{(i)}}$ as required.
To show that $P$ is an $A_2$-$C^{\ast}$-planar algebra we need to show that $P$ is non-degenerate, which is immediate from property $(v)$ in the statement of the theorem, Proposition \ref{Prop:non-degenerate=tr_pos_def} and the fact that $\mathrm{tr}$ is positive definite.
\hfill
$\Box$

\begin{Def}
We will say that an $A_2$-planar algebra $P$ is an $A_2$-planar algebra for the subfactor $N \subset M$ if $P_{0,\infty} = N$, $P_{1,\infty} = M$, $P_{n,\infty} = M_{n-1}$, the sequence $P_{0,0} \subset P_{1,0} \subset P_{2,0} \subset \cdots$ is the tower of relative commutants, and if conditions $(i)$-$(v)$ of Theorem \ref{thm:type_I_subfactors_give_flat_planar_algebras} are satisfied.
\end{Def}

Suppose $P$ is the $A_2$-planar algebra given by the double complex $(B_{i,j})$ for an $SU(3)$ $\mathcal{ADE}$ subfactor $N \subset M$. Then the $A_2$-planar subalgebra $P^{(1)} \subset P$ is the $A_2$-planar algebra given by the double complex for the subfactor $M \subset M_1$.

For the subalgebra $Q$ introduced in $\S$\ref{sec:general_A_2-planar_algebras}, we give an alternative proof of Jones' theorem that extremal subfactors give planar algebras \cite[Theorem 4.2.1]{jones:planar} in the finite depth case.
Jones' proof uses the bimodule setup- he works with the von Neumann algebras themselves, identifying the relative commutants with tensor powers of the von Neumann algebra $M$. The rotation tangle $\rho$ plays an important role in his proof of Theorem 4.2.1, as does the Pimsner-Popa basis.
In our setup we choose to work more directly with the finite dimensional relative commutants themselves. The rotation $\rho$ and the Pimsner-Popa basis do not appear in our proof.
The advantage of our proof is that it extends to our $A_2$ setting, whereas the bimodule setup seems difficult to adapt.

\begin{Cor}
Let $N \subset M$ be a finite depth type $II_1$ subfactor. For each $k$ let $Q_k = N' \cap M_{k-1}$. Then $Q = \bigcup_k Q_k$ has a spherical ($A_1$-)$C^{\ast}$-planar algebra structure (in the sense of Jones), with labelling set $Q$, for which $Z(I_k(x)) = x$, where $I_k(x)$ is the tangle $I_k$ with $x \in Q_k$ as the insertion in its inner disc, and \\
\begin{minipage}[t]{16cm}
 \begin{minipage}[t]{1cm}
  \parbox[t]{1cm}{$(i)$}
 \end{minipage}
  \begin{minipage}[t]{6cm}
  $Z(f_{l}) = \delta e_l$, \qquad $l \geq 1,$
 \end{minipage}
\end{minipage} \\
\begin{minipage}[t]{16cm}
 \begin{minipage}[t]{1cm}
  \mbox{} \\
  \parbox[t]{1cm}{$(ii)$}
 \end{minipage}
 \begin{minipage}[t]{6.5cm}
  \mbox{} \\
   \includegraphics[width=53mm]{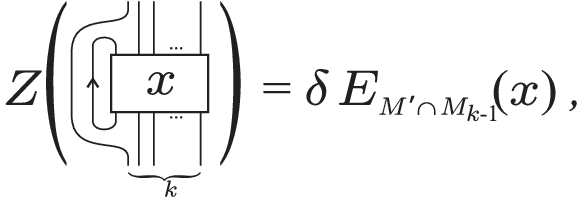}
 \end{minipage}
 \begin{minipage}[t]{6cm}
  \mbox{} \\
   \includegraphics[width=53mm]{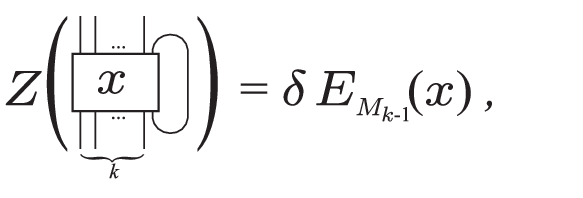}
 \end{minipage}
\end{minipage} \\
\mbox{} \\
\begin{minipage}[t]{16cm}
 \begin{minipage}[t]{1cm}
  \mbox{} \\
  \parbox[t]{1cm}{$(iii)$}
 \end{minipage}
 \begin{minipage}[t]{6cm}
  \mbox{} \\
   \includegraphics[width=53mm]{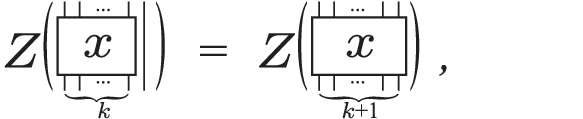}
 \end{minipage}
\end{minipage} \\
\mbox{} \\
\begin{minipage}[t]{16cm}
 \begin{minipage}[t]{1cm}
  \mbox{} \\
  \parbox[t]{1cm}{$(iv)$}
 \end{minipage}
 \begin{minipage}[t]{6cm}
  \mbox{} \\
   \includegraphics[width=53mm]{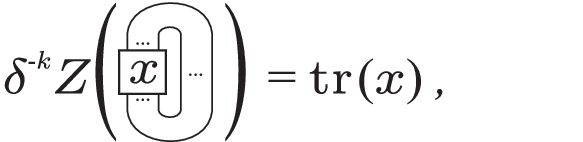}
 \end{minipage}
\end{minipage} \\ \\
for $x \in Q_k$, $k \geq 0$. In condition $(iii)$, the $x$ on the right hand side is considered as an element of $Q_{k+1}$.
Moreover, any other spherical planar algebra structure $Z'$ with $Z'(I_k(x)) = x$ and $(i)$, $(ii)$, $(iv)$ for $Z'$ is equal to $Z$.
\end{Cor}
\emph{Proof:}
We define $Z$ in the same way as above, by converting all the discs of a tangle $T$ to horizontal rectangles and isotoping the tangle so that in each horizontal strip there is either a labelled rectangle, a cup or a cap. Then we define $M_{\cup^{(i)}}$ and $M_{\cap^{(i)}}$ as in (\ref{M_cup}), (\ref{M_cap}). For strip $b_l$ containing a rectangle with label $x_l$, we define $M_{b_l}$ as in Theorem \ref{thm:type_I_subfactors_give_flat_planar_algebras}, using the connection on the principal graph $\mathcal{G}$ and its reverse graph $\widetilde{\mathcal{G}}$. The cup-cap simplification of Figure \ref{fig:cup-cap_simplifications} follows from (\ref{cup-cap_simplification}) and (\ref{cup-cap_simplification2}). The invariance of $Z$ under isotopies involving rectangles as in Figures \ref{fig:isotopies_involving_boxes}, \ref{fig:rotation_of_internal_boxes} follows as in the proof of Theorem \ref{thm:type_I_subfactors_give_flat_planar_algebras}. That closed loops give a scalar factor of $\delta$ follows from (\ref{proof_of_K1}), where the Perron-Frobenius eigenvalue now is $\delta$.
Properties $(i)$-$(iv)$ are proved in the same way as properties $(i)$, $(iii)$, $(iv)$, $(v)$ of Theorem \ref{thm:type_I_subfactors_give_flat_planar_algebras}, and uniqueness is proved as in \cite{jones:planar}.
\hfill
$\Box$

\subsection{Presentation of the Path Algebra for $\mathcal{A}^{(n)}$ as a $PTL$ Algebra} \label{sect:planar_alg_for_A=PTL}

We now show that each $B_{i,j}$ for the double sequence $(B_{i,j})$ defined above for $\mathcal{G} = \mathcal{A}^{(n)}$ also has a presentation as $\mathcal{PTL}_{i,j}$, where $\mathcal{PTL}_{i,j}$ is the quotient of $PTL_{i,j}$ by the subspace of zero-length vectors, as in Section \ref{sect:dimP}.

Now $B_{1,j} \cong \mathcal{PTL}_{1,j}$ by Lemma \ref{Lemma:Vn=Hn(q)}. Let $\psi:B_{1,j} \rightarrow \mathcal{PTL}_{1,j}$ be the isomorphism given by $\psi(U_{-k}) = W_{-k}$, $k=0,\ldots,j-1$. We define maps $\varrho_i$ for $i \geq 2$ by $\varrho_2 = \varphi$, $\varrho_3 = \omega \varphi$, $\varrho_4 = \varphi \omega \varphi$, $\varrho_5 = \omega \varphi \omega \varphi$, $\ldots$.
Let $x = \sum_{\gamma, \gamma'} \lambda_{\gamma, \gamma'} (\gamma, \gamma')$, $\lambda_{\gamma, \gamma'} \in \mathbb{C}$, be an element of $B_{i,j}$. Then $Z(\varrho_i^{-1}(x)) \in B_{1,i+j-1}$. We set $x_W \in \mathcal{PTL}_{1,i+j-1}$ to be the element $\psi(Z(\varrho_i^{-1}(x)))$, and since $Z(W_{-k}) = U_{-k}$ we have $Z(x_W) = Z(\varrho_i^{-1}(x))$. For any $x \in B_{i,j}$, $\varrho_i(x_W) \in PTL_{i,j}$ and $Z(\varrho_i(x_W)) = Z(\varrho_i(Z(x_W))) = Z(\varrho_i(Z(\varrho_i^{-1}(x)))) = Z(\varrho_i \varrho_i^{-1}(x)) = Z(I_{i,j}(x)) = x$. In fact, $\varrho_i(x_W) \in \mathcal{PTL}_{i,j}$, since if $\langle \varrho_i(x_W), \varrho_i(x_W) \rangle = 0$, then $\langle x,x \rangle = \langle \varrho_i^{-1}(x), \varrho_i^{-1}(x) \rangle = \langle x_W, x_W \rangle = \langle \varrho_i(x_W), \varrho_i(x_W) \rangle = 0$ as in Section \ref{sect:dimP}, so that $\varrho_i(x_W)$ is a zero-length vector if and only if $x$ is. Then for every $x \in B_{i,j}$ there exists a $y = \varrho_i(x_W) \in \mathcal{PTL}_{i,j}$ such that $Z(y) = x$, so that $Z$ is surjective. Since, by (\ref{eqn:dimP=dimP-2}), $\mathrm{dim}(\mathcal{PTL}_{i,j}) = \mathrm{dim}(\mathcal{PTL}_{1,i+j-1}) = \mathrm{dim}(B_{1,i+j-1}) = \mathrm{dim}(B_{i,j})$, this element $y$ is unique and $Z$ is a bijection. By its definition, $Z$ is linear and preserves multiplication. Then $Z:\mathcal{PTL}_{i,j} \rightarrow B_{i,j}$ is an isomorphism, and we have shown the following:
\begin{Lemma} \label{Lemma:B=PTL}
In the double sequence $(B_{i,j})$ defined above for $\mathcal{G} = \mathcal{A}^{(n)}$, each $B_{i,j}$ is isomorphic to $\mathcal{PTL}_{i,j}$
\end{Lemma}

In particular, there is a presentation of the path algebra for the 01-part $\mathcal{A}_{01}^{(n)}$ of $\mathcal{A}^{(n)}$ given by vectors of non-zero length, which are linear combinations of tangles generated by Kuperberg's $A_2$ webs, where $A(\mathcal{A}_{01}^{(n)})_k$ is the space of all such tangles on a rectangle with $k$ vertices along the top and bottom, with the orientations of the vertices alternating.

As a corollary to Lemma \ref{Lemma:B=PTL}, we thus obtain the following description of the ($A_1$-)planar algebra for the Wenzl subfactor \cite{wenzl:1988} which has principal graph $\mathcal{A}_{01}^{(n)}$. Let $Q$ be the ($A_1$-)planar algebra generated by $\; \bigcup_{i \geq 0} \mathcal{PTL}_{i,0} \;$ with relations K1-K3, and let $I^0$ be the ideal generated by the zero-length elements of $Q$, that is, $I^0 = \bigcup_{k > 0} I^0_k$ where $I^0_k = \{ x \in Q_k | \; tr(x^{\ast} x) = 0 \}$.
For a family $\{ U_m \}_{m \geq 1}$ of self-adjoint operators which generate an $A_2$-Temperley-Lieb algebra with parameter $\delta = q+q^{-1}$, where $q=e^{i \pi / n}$, let $M = \langle U_1, U_2, U_3, \ldots \; \rangle$ and $N = \langle U_2, U_3, U_4, \ldots \; \rangle$.
\begin{Cor}
The ($A_1$-)planar algebra $P$ corresponding to Wenzl's subfactor $N \subset M$ is the quotient $P = Q/I^0$.
\end{Cor}

\subsection{Comparison of $PTL_{i,0}$ with the Temperley-Lieb Algebra}

We will now compare $PTL_{i,0}$ with the Temperley-Lieb algebra. In particular we will write a basis for $PTL_{3,0}$, which will be given by the Temperley-Lieb diagrams $TL_3 = \textrm{alg}(\mathbf{1},f_1,f_2)$, and an extra diagram which contains trivalent vertices.

\begin{figure}[tb]
\begin{minipage}[t]{6cm}
 \begin{center}
  \includegraphics[width=22mm]{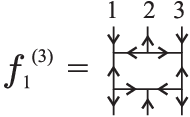}\\
 \caption{Element $f_1^{(3)}$} \label{fig:f1(3)}
 \end{center}
\end{minipage}
\hfill
\begin{minipage}[t]{9cm}
 \begin{center}
  \includegraphics[width=60mm]{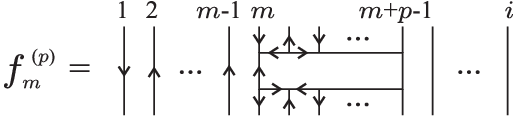}\\
 \caption{Element $f_m^{(p)}$ (with $m$ odd)} \label{fig:fm(p)}
 \end{center}
\end{minipage}
\end{figure}

Since $PTL_{1,2} = \textrm{alg}(\mathbf{1}_{1,2}, W_{-1}, W_0)$, we have $\varphi(W_{-1}) = q^{8/3} W_{-1}$ and $\varphi(W_0) = q^{5/3} \mathbf{1}_{2,1} - q^{-1/3} f_1$ so that $PTL_{2,1} = \textrm{alg}(\mathbf{1}_{2,1}, W_0, f_1)$. The action of $\omega$ on $\mathcal{PTL}_{2,1}$ is given by $\omega(f_1) = f_1$, $\omega(W_0) = f_1^{(3)} - q \alpha^2 f_1 f_2 - q^{-1} \alpha^2 f_2 f_1$ and $\omega(f_1 - q W_0 f_1 - q^{-1} f_1 W_0 + W_0 f_1 W_0) = f_2$, where $f_1^{(3)}$ is the tangle illustrated in Figure \ref{fig:f1(3)}.
We see that $PTL_{3,0}$ is generated by $\mathbf{1}, f_1, f_2$ and $f_1^{(3)}$. This new element $f_1^{(3)}$ cannot be written as a linear combination of products of $\mathbf{1}, f_1$ and $f_2$. Thus we see that $PTL_{3,0}$ is generated by $TL_3 = \textrm{alg}(\mathbf{1},f_1,f_2)$ and the extra element $f_1^{(3)}$. The following hold for $f_1^{(3)}$ (they can be easily checked by drawing pictures):
\begin{itemize}
\item[(i)] $\left( f_1^{(3)} \right)^2 = \delta f_1^{(3)} + \alpha (f_1 + f_2) + \alpha^2 (f_1 f_2 + f_2 f_1)$,

\item[(ii)] $f_1 f_1^{(3)} = \delta f_1 + \delta \alpha f_1 f_2$, \qquad $f_2 f_1^{(3)} = \delta f_2 + \delta \alpha f_2 f_1$,

\item[(iii)] $f_i f_1^{(3)} f_i = \delta^3 \alpha^{-1} f_i$, \quad $i=1,2$,

\item[(iv)] $f_1^{(3)} f_i f_1^{(3)} = \delta^2 (f_1 + f_2) + \delta^2 \alpha (f_1 f_2 + f_2 f_1)$, \quad $i=1,2$.
\end{itemize}

We define the operator $g_1^{(3)}$ to be $g_1^{(3)} = Z(f_1^{(3)})$. Then $A(\mathcal{A}^{(n)}_{01})_3 = \textrm{alg}(\mathbf{1}, e_1, e_2, g_1^{(3)})$.

For $n \geq 6$, with the rows and columns indexed by the paths of length 3 on $\mathcal{A}^{(n)}_{01}$ which start at vertex $(0,0)$, $g_1^{(3)}$ can be written explicitly as the matrix
$$
g_1^{(3)} = \left( {
\begin{array}{cccc}
                 [2]^3/[3] & \sqrt{[2]^3 [4]}/[3] & 0 & 0 \\
                 \sqrt{[2]^3 [4]}/[3]  & [4]/[3] & 0 & 0 \\
                 0 & 0 & [2] & 0 \\
                 0 & 0 & 0 & 0 \\
               \end{array} } \right).
$$
For $n=5$, $g_1^{(3)} = \alpha \mathbf{1} - e_1 - e_2 + \alpha e_1 e_2 + \alpha e_2 e_1$, so is a linear combination of $\mathbf{1}, e_1$ and $e_2$. This is not a surprise since $\mathcal{A}^{(5)}_{01}$ is just the Dynkin diagram $A_4$, and we know that $A(A_4)_3$ is generated by $\mathbf{1}$, $e_1$ and $e_2$. Note also that in this case we have $\alpha = \delta = \sin(2 \pi i/5)$.

It appears that $PTL_{i,j} = \textrm{alg}(\mathbf{1}_{i,j}, W_{-k}, f_l, f_m^{(3)} | \, k=0,\ldots,j-1; \, l=1,\ldots,i-1; \, m=1,\ldots,i-2)$, where $f_m^{(3)}$ is the tangle illustrated in Figure \ref{fig:fm(p)}, where $p=3$.
The more general elements $f_m^{(p)}$ illustrated in Figure \ref{fig:fm(p)} have an internal face with $2p$ edges. These elements are generated by $f_m^{(3)}$ and $f_l$: $f_m^{(3)} f_{m+1}^{(3)} = f_m^{(4)} + \textrm{ linear combination of } f_l$, and more generally, $f_m^{(3)} f_{m+1}^{(p)} = f_m^{(p+1)} + \textrm{ linear combination of } f_l, f_k^{(3)}$.

We know that $PTL_{0,j}$ is generated by $W_{-k}$, $k=1,\ldots,j-1$, by Lemma \ref{Lemma:tangles_as_Wi's}. As shown above using the maps $\omega$ and $\varphi$, any element in $PTL_{i,j}$ can be obtained from the $W_{-k}$ by using the braiding. Thus $PTL_{i,j}$ is generated by the $W_{-k}$ and the braiding.

\paragraph{Acknowledgements}

This paper is based on work in \cite{pugh:2008}. The first author was partially supported by the EU-NCG network in Non-Commutative Geometry MRTN-CT-2006-031962, and the second author was supported by a scholarship from the School of Mathematics, Cardiff University.
We thank the referees for their helpful comments. We are grateful to one of the referees for suggesting the use of the braiding in the proof of Lemma \ref{Lemma:tangles_as_Wi's}, which allowed for a simpler proof.

\end{document}